\documentclass[10pt]{article}

\usepackage[text={6in,650pt},centering]{geometry}

\usepackage{color}

\usepackage{comment}
\usepackage{datetime}

\usepackage{eucal}

\usepackage[all]{xy}

\usepackage{enumitem}

\usepackage{epsfig}
\usepackage{lpic}

\usepackage{amssymb, amsmath, amsthm, eucal, latexsym}

\usepackage{mathptmx}
\usepackage[scaled=.90]{helvet}
\usepackage{courier}

\usepackage[
color=orange!80,
bordercolor=black,
textwidth=.8in]
{todonotes}
\makeatletter \providecommand\@dotsep{5}
\makeatother

\usepackage[pdftex]{hyperref}

\hypersetup{
  colorlinks = true,
  linkcolor = black,
  urlcolor =  blue,
  citecolor = black,
}

\newcommand{\bfC}{\mathbf{C}}

\newcommand{\bfH}{\mathbf{H}}
\newcommand{\bfQ}{\mathbf{Q}}
\newcommand{\bfV}{\mathbf{V}}

\newcommand{\bSet}{\overline{\mathbf{Set}}}
\newcommand{\bCan}{\overline{\mathbf{Cantor}}}
\newcommand{\Set}{\mathbf{Set}}
\newcommand{\Can}{\mathbf{Cantor}}
\newcommand{\Exp}{\mathbf{Exp}}
\newcommand{\Lev}{\mathbf{Lev}}
\newcommand{\Tho}{\mathbf{Tho}}

\newcommand{\bbK}{\mathbb{K}}
\newcommand{\bbM}{\mathbb{M}}
\newcommand{\bbN}{\mathbb{N}}
\newcommand{\bbQ}{\mathbb{Q}}

\newcommand{\bbS}{\mathbb{S}}
\newcommand{\bbX}{\mathbb{X}}
\newcommand{\bbZ}{\mathbb{Z}}

\newcommand{\calE}{\mathcal{E}}
\newcommand{\calI}{\mathcal{I}}
\newcommand{\calL}{\mathcal{L}}
\newcommand{\calM}{\mathcal{M}}
\newcommand{\calP}{\mathcal{P}}

\newcommand{\rmB}{\mathrm{B}}
\newcommand{\rmC}{\mathrm{C}}
\newcommand{\rmD}{\mathrm{D}}
\newcommand{\rmE}{\mathrm{E}}
\newcommand{\rmF}{\mathrm{F}}

\newcommand{\rmH}{\mathrm{H}}
\newcommand{\rmK}{\mathrm{K}}
\newcommand{\rmL}{\mathrm{L}}
\newcommand{\rmM}{\mathrm{M}}

\newcommand{\rmP}{\mathrm{P}}
\newcommand{\rmS}{\mathrm{S}}
\newcommand{\rmT}{\mathrm{T}}
\newcommand{\rmU}{\mathrm{U}}
\newcommand{\rmV}{\mathrm{V}}
\newcommand{\rmW}{\mathrm{W}}

\newcommand{\bid}{{\bf id}}
\DeclareMathOperator{\id}{id}
\DeclareMathOperator{\op}{op}
\DeclareMathOperator{\sk}{sk}
\DeclareMathOperator{\ord}{ord}
\DeclareMathOperator{\Star}{Star}
\DeclareMathOperator{\Link}{Link}
\DeclareMathOperator{\Map}{Map}
\DeclareMathOperator{\Rep}{Rep}

\renewcommand{\phi}{\varphi}

\renewcommand{\emptyset}{\text{\upshape\O}}

\newcommand{\KO}{\mathbb{KO}}

\DeclareMathOperator*{\colim}{colim}
\DeclareMathOperator*{\hocolim}{hocolim}

\DeclareMathOperator{\Aut}{Aut}

\newcommand{\al}{\alpha}
\newcommand{\la}{\lambda}
\newcommand{\s}{\sigma}
\newcommand{\Si}{\Sigma}

\newcommand{\rar}{\longrightarrow}
\newcommand{\lar}{\leftarrow}
\newcommand{\inc}{\hookrightarrow}

\newcommand{\x}{\times}
\newcommand{\minus}{\setminus}

\renewcommand{\le}{\leqslant}
\renewcommand{\ge}{\geqslant}
\newcommand{\IV}{\rotatebox[origin=c]{90}{\ensuremath{\leqslant}}}
\newtheorem{theorem}{Theorem}
\newtheorem*{theorem*}{Main theorem}
\newtheorem{proposition}[theorem]{Proposition}
\newtheorem{corollary}[theorem]{Corollary}
\newtheorem*{corollary*}{Corollary}
\newtheorem{lemma}[theorem]{Lemma}
\newtheorem{prop}[theorem]{Proposition}
\newtheorem{cor}[theorem]{Corollary}
\newtheorem{lem}[theorem]{Lemma}

\theoremstyle{definition}
\newtheorem{definition}[theorem]{Definition}

\newtheorem{remark}[theorem]{Remark}
\newtheorem{Def}[theorem]{Definition}
\newtheorem{ex}[theorem]{Example}
\newtheorem{rem}[theorem]{Remark}
\newtheorem{notation}[theorem]{Notation}

\numberwithin{theorem}{section}
\numberwithin{equation}{section}
\numberwithin{figure}{section}

\newtheorem{Th}{Theorem}
\newtheorem{Co}[Th]{Corollary}%\renewcommand{\theCo}{\arabic{Co}}

\title{{The homology of the Higman--Thompson groups}}

\author{Markus Szymik and Nathalie Wahl}

\newdateformat{mydate}{\monthname~\twodigit{\THEYEAR}}
\date{\mydate\today}

\begin{document}

\maketitle

\renewcommand{\abstractname}{}

\begin{abstract}
\noindent  
We prove that Thompson's group~$\rmV$ is acyclic, answering a 1992 question of Brown in the positive. More generally, we identify the homology of the Higman--Thompson groups~$\rmV_{n,r}$ with the homology of the zeroth component of the infinite loop space of the mod~$n-1$ Moore spectrum. As~$\rmV=\rmV_{2,1}$, we can deduce that this group is acyclic. 
Our proof %, in fact, gives a space level identification, after plus construction, and 
involves establishing homological stability with respect to~$r$, as well as a computation of the algebraic K-theory of the category of finitely generated free Cantor algebras of any type~$n$. 
%As a special case we deduce that Thompson's group~$\rmV=\rmV_{2,1}$ is acyclic. 
%For all~\hbox{$a\geqslant3$}, the groups~$\rmV_{a,r}$ are only rationally acyclic. 
\vspace{\baselineskip}

\vspace{\baselineskip}

\end{abstract}

%%%

\section*{Introduction}

About half a century ago, Thompson introduced a group~$\rmV$ together with subgroups~$\rmF\leqslant\rmT\leqslant\rmV$ in order  to construct examples of finitely presented groups with  unsolvable word problem. 
Thompson's groups have since developed a life of their own, relating to many branches of  mathematics. %It is therefore  natural to wonder about the homological properties of these groups. 
The homology of the group~$\rmF$ was computed by Brown and Geoghegan~\cite{Brown+Geoghegan}; it is free abelian of rank~$2$ in all strictly positive degrees. The homology of the group~$\rmT$ was computed by Ghys and Sergiescu~\cite{Ghys+Sergiescu}; it is isomorphic to the homology of the free loop space on the~$3$-sphere.
%~(Of course, the latter is {\em not} a classifying space for the group~$\rmT$.) 
As for Thompson's group~$\rmV$ itself, Brown~\cite{Brown:Suggestion} proved that it is rationally acyclic and suggested that it might even be integrally so. 
In the present paper, we prove that~$\rmV$ is indeed integrally acyclic. 
%However, the existing partial results rationally or in low-dimensions could hardly be considered as evidence for a conjecture, and the problem remained open.

Thompson's group~$\rmV$ fits into the more general family of the Higman--Thompson groups~$\rmV_{n,r}$ for~\hbox{$n\geqslant 2$} and~$r\geqslant 1$, with~$\rmV=\rmV_{2,1}$ as the first case: %~(see~\cite{Higman}).
A Cantor algebra of type~$n$ is a set~$X$ equipped with a bijection~\hbox{$X^n\cong X$}, and the Higman--Thompson group~$\rmV_{n,r}$ is the automorphism group of the free Cantor algebra~$\rmC_{n}[r]$ of type~$n$ on~$r$ generators. The main result of this text is an identification of the homology of all of the groups~$\rmV_{n,r}$ in terms of a well-known object of algebraic topology: the mod~\hbox{$n-1$} Moore spectrum~$\bbM_{n-1}$.  

\begin{Th}%[Theorem~\ref{thm:main}]\label{homologyofVnr}
\label{thm:main}
For any~$n\geqslant2$ and~$r\geqslant1$ there is a map~$\rmB\rmV_{n,r} \to \Omega^\infty_0\bbM_{n-1}$ 
inducing an isomorphism 
\[
\rmH_*(\rmV_{n,r};M) \ \stackrel{\cong}\longrightarrow\ \rmH_*(\Omega^\infty_0\bbM_{n-1};M).
\] 
in homology for any coefficient system~$M$ on~$\Omega^\infty_0\bbM_{n-1}$. 
\end{Th}

Here the space~$\Omega^\infty_0\bbM_{n-1}$ is the zeroth component of the infinite loop space~$\Omega^\infty\bbM_{n-1}$ that underlies the mod~$n-1$ Moore spectrum~$\bbM_{n-1}$. 
Note that the target of the isomorphism does not depend on~$r$.

In the case~$n=2$, the spectrum~$\bbM_{n-1}$ is contractible and the above result answers Brown's question: 

\begin{Co}[{\bf Theorem~\ref{thm:Vacyclic}}]
Thompson's group~$\rmV=\rmV_{2,1}$ is acyclic.
\end{Co}

In \cite{Brown:Suggestion}, Brown indicates that his argument for the rational acyclicity of~$\rmV$ extends to prove rational acyclicity for all groups~$\rmV_{n,r}$. 
When~$n$ is odd, the group~$\rmV_{n,r}$ was known not to be integrally acyclic just from the computation of its first homology group, which is~$\bbZ/2$ in that case. 
Our main theorem applied to the case~$n\geqslant 3$ completes the picture, giving 
a proof of rational acyclicity for all~$\rmV_{n,r}$, and at the same time showing that integral acyclicity only holds in the special case~$n=2$\,: 

\begin{Co}[{\bf Theorem~\ref{thm:Qacyclic}}]
For all~$n\geqslant 3$, the group~$\rmV_{n,r}$ is rationally but not integrally acyclic.   
\end{Co}

We give in the last section of the paper some additional explicit
consequences of Theorem~\ref{thm:main}. In particular, we confirm and
complete the known information about the abelianizations and Schur
multipliers of the groups~$\rmV_{n,r}$~(Propositions~\ref{prop:H1}
and~\ref{prop:H2}), and compute the first non-trivial homology group
of~$\rmV_{n,r}$ for each $n$ and $r$~(Proposition~\ref{prop:non-trivial}). 
%We also deduce that all of the groups~$\rmV_{n,r}$ are rationally acyclic and show that, for~\hbox{$a\geqslant 3$}, the groups~$\rmV_{n,r}$ that are {\em not} integrally acyclic (Theorem~\ref{thm:Qacyclic}).
When~$n$ is odd, the commutator subgroups~$\rmV_{n,r}^+$ is an index two subgroup, and our methods can also be applied to study this group~(Corollary~\ref{cor:commutator}).

The proof of our main theorem %Theorem~\ref{thm:main} 
rests on two pillars. The first is homological stability: For any fixed~\hbox{$n\geqslant 2$}, the Higman--Thompson groups~$\rmV_{n,r}$ fit into a canonical diagram
\begin{equation}\label{eq:canonical_diagram}
\rmV_{n,1}
\longrightarrow\rmV_{n,2}
\longrightarrow\rmV_{n,3}
\longrightarrow\ \cdots\tag{$\star$}
\end{equation}
of groups and (non-surjective) homomorphisms, and we show that the maps~$\rmV_{n,r}\to\rmV_{n,r+1}$ induce isomorphisms in homology for large~$r$ in any fixed homological degree. 
The definition of Cantor algebras leads to isomorphisms~\hbox{$\rmC_{n}[r]\cong\rmC_{n}[r+(n-1)]$} for all~$n\geqslant 2$ and~$r\geqslant1$, giving isomorphisms%~$\rmV\cong\rmV_{2,r}$ for all~$r\geqslant1$ and more generally
~$\rmV_{n,r}\cong\rmV_{n,r+(n-1)}$ for all~$n\geqslant 2$ and~$r\geqslant1$. 
Using %additionally the isomorphisms~\hbox{$\rmV_{n,r}\cong\rmV_{n,r+(n-1)}$}
these isomorphisms, we obtain that the stabilization maps are actually
isomorphisms in homology in  {\em all} degrees, see Theorem~\ref{thm:HS}. 

%%%

To prove homological stability, we use the framework of \cite{Randal-Williams+Wahl}. 
%A Cantor algebra of arity~$n$ is a set~$X$ equipped with an isomorphism~$X\cong X^n$, and the group~$\rmV_{n,r}$ is the automorphism group of the free Cantor algebra~$\rmC_{n,r}$ of arity~$n$ on~$r$ generators. The groupoid~$\Can^\x_n$ of free Cantor alegras of arity~$n$ and their isomorphisms does not satisfy the hypotheses in \cite{Randal-Williams+Wahl}, but 
The main ingredient for stability is the proof of high connectivity of a certain simplicial complex of independent sets in the free Cantor algebra~$\rmC_{n}[r]$. 
It follows from \cite{Randal-Williams+Wahl} that homological stability also holds with appropriate abelian and polynomial twisted coefficients.

%%% 

Our stability theorem can be reformulated as 
saying that the map
\[
\rmV_{n,r}\longrightarrow\bigcup_{r\geqslant 1}\rmV_{n,r}=\rmV_{n,\infty}
\]
is a homology isomorphism, where the union is defined using the maps in the diagram ($\star$), and 
%the isomorphism 
%$$\rmH_*(\rmV_{n,r})\cong \rmH_*(\rmV_{n,\infty})$$
%for~$\rmV_{n,\infty}$ the colimit of the diagram ($\star$).
our second pillar is the identification of the homology
of~$\rmV_{n,\infty}$. This is achieved by the identification of the
K-theory of the groupoid of free Cantor algebras of type~$n$, as we describe now. 

Let~$\Can_n^\x$ denote the category of free Cantor algebras of
type~$n$ with morphisms their isomorphisms.  The category~$\Can_n^\x$
is symmetric monoidal, and hence has an associated
spectrum~$\bbK(\Can_n^\x)$, its algebraic K-theory. We denote by~$\Omega^\infty_0\bbK(\Can_n^\x)$ the zeroth component of its associated infinite loop space. 
Applying the group completion theorem, we get a map  
\[
\rmB\rmV_{n,\infty}\longrightarrow\Omega^\infty_0\bbK(\Can_n^\x),
\]
defined up to homotopy, that induces an isomorphism in homology with all local coefficient systems on the target 
%yields that the homology of~$\rmV_{n,\infty}$ identifies with the homology of the zeroth component of the infinite loop space of this spectrum: 
%$$\rmH_*(\rmV_{n,\infty})\cong \rmH_*(\Omega^\infty_0\bbK(\Can_n^\x))$$
%where~$\bbK(\Can_n^\x)$ is the spectrum associated to the symmetric monoidal category~$\Can^\x_n$ 
(see Theorem~\ref{thm:gpcomp}).  
%the homology of~$\rmV_{n,\infty}$ and that of the classifying space of the category of free Cantor algebras of arity~$n$. 
%Now, by their defining property, the free Cantor algebra of arity~$n$ on a set and that on the same set crossed with a set of size~$n$ are isomorphic.  What we show is that the category~$\Can_n^\x$ is a homotopy mapping cylinder of the functor, defined on the category of finite sets and bijections, that crosses with a set of size~$n$.  
%More formally, 
Using a model of Thomason, we identify the classifying space~$|\Can_n^\x|$ with that of a homotopy colimit~$\Tho_n$ in symmetric monoidal categories %describing that mapping cylinder: 
build out of the category of finite sets and the functor that takes the product with a set of cardinality~$n$:
\[
|\Can_n^\x|\simeq |\Tho_n|,
\]
where the equivalence respects the symmetric monoidal structure, see Theorem~\ref{back}. 
In particular, the two categories have equivalent algebraic K-theory spectra:~$\bbK(\Can_n^\x)\simeq \bbK(\Tho_n)$. The main theorem then follows from an identification
\[
\bbK(\Tho_n)\simeq \bbM_{n-1},
\]
see Theorem~\ref{thm:hocolim_is_Moore}. The idea behind the last two equivalences is as follows. 
%that, by their defining property, the free Cantor algebra of arity~$n$ on a set and that on the same set crossed with a set of size~$n$ are isomorphic. 
The category~$\Tho_n$ is a homotopy mapping torus of the functor, defined on the category of finite sets and bijections, that takes the product with a set of size~$n$. Thinking of the finite sets as the generating sets of free Cantor algebras of a given type~$n$,  this functor implements, for any~$r$, the identification of a Cantor algebra~$\rmC_{n}[r]$  with the Cantor algebra~\hbox{$\rmC_{n}[rn]=\rmC_{n}[r+r(n-1)]$}, which reflects the defining property of Cantor algebras.  Now the K-theory of the groupoid of finite sets is the sphere spectrum, by the classical Barratt--Priddy--Quillen theorem, and in spectra, this mapping torus equalizes multiplication by~$n$ with the identity on the sphere spectrum, which leads to the Moore spectrum~$\bbM_{n-1}$. 

%%%

The paper is organized as follows. In Section~\ref{sec:Cantor}, we introduce the Cantor algebras and the Higman--Thompson groups, and we give some of their basic properties that will be needed later in the paper. In Section~\ref{sec:spaces}, we show how the groups~$\rmV_{n,r}$ fit into the set-up for homological stability of \cite{Randal-Williams+Wahl} and construct the spaces relevant to the proof of homological stability, which is given in Section~\ref{sec:connectivity}. The following Section~\ref{sec:classifying} is devoted to the homotopy equivalence~\hbox{$|\Can_n^\x|\simeq |\Tho_n|$}, which is given as a composition of three homotopy equivalences. Section~\ref{sec:main} then relates the first of these two spaces to~$\rmV_{n,\infty}$ and the second to the Moore spectrum; the section ends with the proof of the main theorem. 
Finally, Section~\ref{sec:explicit} draws the computational consequences of our main theorem. Throughout the paper, the symbol~$A$ will denote a fixed finite set of cardinality~$n$. 

%%%

% !TEX root = allinone_main_v2.tex

\section{Cantor algebra}\label{sec:Cantor}

In this section, we recall some facts we need about Cantor algebras, and their groups of automorphisms,  the Higman--Thompson groups. We follow Higman's own account~\cite{Higman}. See also~\cite[Sec.~4]{Bro87} for a shorter survey.

%\Nnote{Higman works with~$n\ge 2$ so we need to exclude~$n=1$ early on. }

Let~$A$ be a finite set  of cardinality at least~$2$ and 
let
\[
A^*=\bigsqcup_{n\geqslant 0}A^n
\]
denote the word monoid on the set~$A$. This is the free monoid generated by~$A$, with unit~$\emptyset\in A^0$ and multiplication by juxtaposition. 
By freeness, a (right) action of the monoid~$A^*$ on a set~$S$ is uniquely determined by a map of sets~$S\times A\to S$. Such a map has an adjoint
\begin{equation}\label{eq:adjoint}
S\longrightarrow S^A=\Map(A,S).
\end{equation}

\begin{definition}\label{def:Cantor}
A {\em Cantor algebra of type~$A$} is an~$A^*$--set~$S$ such that the
adjoint structure map~\eqref{eq:adjoint} is a bijection. The morphisms of Cantor algebras of type~$A$ are the maps of~$A^*$--sets, that is the set maps that commute with the action of~$A$.  
\end{definition}
Such objects go also under the name {\em J\'onsson--Tarski algebras}.
%Ref to Smirnov? 73 or 92? or Jonsson-Tarski? probably not really necessary as it's in Wikipedia :). 

For any finite set~$X$, there exists a free Cantor algebra~$\rmC_A(X)$ of type~$A$ generated by~$X$. It can be constructed from the free~$A^*$--set 
with basis~$X$, namely the set
\[
\rmC_A^+(X):=X\x A^*=\bigsqcup_{n\geqslant 0} X\x A^n
\]
by formally adding elements so as to make the adjoint of the map defining the action bijective.~(See~\cite[Sec.~2]{Higman}; the set~$\rmC_A^+(X)$ is denoted~$X\langle A\rangle$ in \cite{Higman}.) Throughout the paper we will work with free Cantor algebras, but only the elements of the canonical free~$A^*$--set~$\rmC^+_A(X)\subset \rmC_A(X)$ will play a direct role.

\begin{figure}[ht]
\centering
\includegraphics[width=0.7\textwidth]{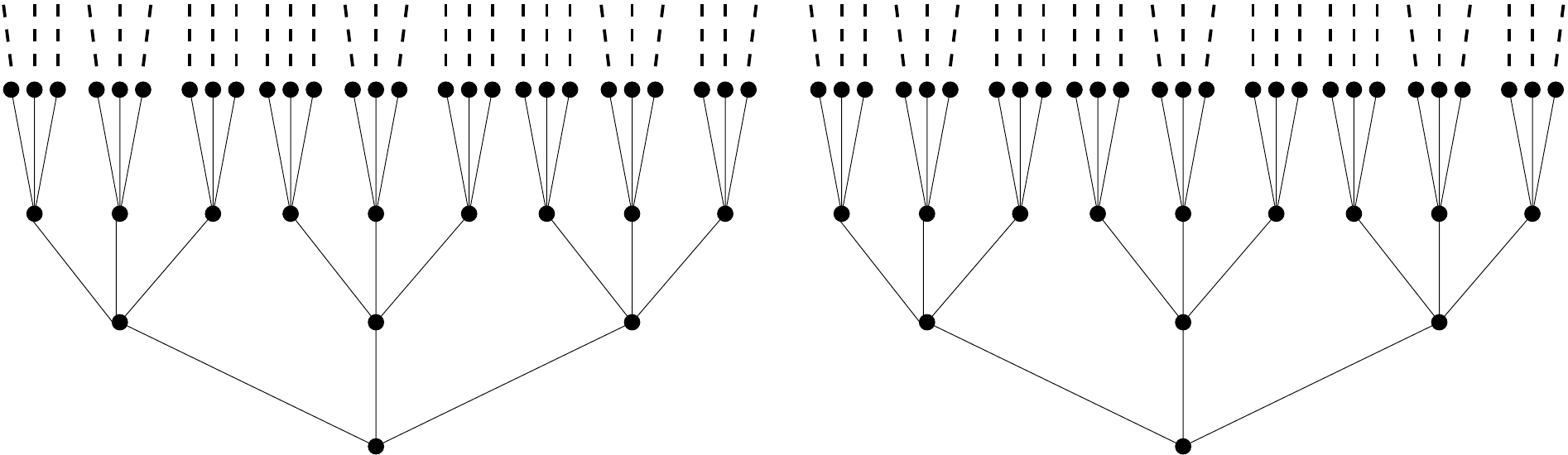}
\caption{The set~$\rmC_{A}^+[2]$ for~$|A|=3$ identifies with the vertices of a forest consisting of two infinite planar ternary trees.}\label{fig:forest}
\end{figure}

The set~$A$ will be fixed throughout the paper to be a set 
$$A=\{a_1,\dots,a_n\}$$ 
with~$n\ge 2$ ordered elements. %(Though some of our results do hold for~$n=1$, the case~$n\ge 2$ is the interesting case, so we will restrict to that. See also Remark~\ref{rem:n=1}.) 
The set~$A$ comes thus with a canonical isomorphism to the set~$[n]=\{1,\dots,n\}$, but we prefer giving it a different name to emphasize its special role. 
As for the generating set, we will be particularly interested in the case~$X=[r]=\{1,\dots,r\}$, though we will also make use of the Cantor algebras and~$A^*$--sets~$\rmC_A(E)$ and~$\rmC_A^+(E)$ generated by other sets build from~$[r]$ and~$A$. 
%in which case we will sometimes shorten the notation and just write~$\rmC_{n,r}$  for~$\rmC_{A}(X)$ and 
%$\rmC_{n,r}^+$ for~$\rmC_{A}^+(X)$. 
Note that the elements of~$\rmC_{A}^+[r]$ canonically identify with the vertices of a planar forest consisting of~$r$ rooted infinite~$n$--ary trees, as in Figure~\ref{fig:forest}. It can be useful to have this picture in mind in what follows, and interpret results using it. 

Our main object of study, the Higman--Thompson groups~$\rmV_{n,r}$, are the automorphism groups of the free Cantor algebras~$\rmC_A[r]:=\rmC_{A}(X)$ for~$X=[r]$: 

\begin{definition}
Let~$A=\{a_1,\dots,a_n\}$ with~$n\geqslant 2$ and let~$r\geqslant 1$. The {\em Higman--Thompson group}
\[
\rmV_{n,r}=\Aut(\rmC_{A}[r])
\]
is the automorphism group of the free Cantor algebra of type~$A$ on~$r$ generators. 
\end{definition}

\subsection{Bases, independent sets and expansions}\label{sec:exp}

We need to understand  isomorphisms of Cantor algebras. 
By freeness, a morphism of Cantor algebras from~$\rmC_A(X)$ to a Cantor algebra~$S$ is determined by its value on the generating set~$X$.  
For instance, the canonical map~$X\x A \to \rmC_A^+(X)\to \rmC_A(X)$ induces a map~$\rmC_A(X\x A)\to \rmC_A(X)$, which one can show is an isomorphism, using the Cantor algebra structure map of~$\rmC_A(X\x A)$. An isomorphism~$f\colon\rmC_A(X)\to \rmC_A(Y)$ of Cantor algebras  can in general take the generating set~$X$ to any generating set of~$Y$. In this section, we study generating sets of Cantor algebras, in particular those called expansions, preparing for Higman's simple description of  isomorphisms between free Cantor algebras, which is recalled in the following section. 

\begin{definition}
A subset~$S\subset \rmC_A(X)$ is called 
%an {\em independent set} if the induced map~$\rmC_A(S)\to\rmC_A(X)$ is injective, and we say that~$S$ is 
a {\em  basis} for the free Cantor algebra~$\rmC_A(X)$ if  the induced homomorphism~\hbox{$\rmC_A(S)\to\rmC_A(X)$} of Cantor algebras is an isomorphism. 
%By \cite[Lem.~2.5]{Higman}, every finite independent set of~$\rmC_A^+(X)$ is contained in a finite basis also inside~$\rmC_A^+(X)$.
\end{definition}

\begin{definition}
Given a subset~$Y$ of a free Cantor algebra~$\rmC_A(X)$, an {\em expansion} of~$Y$ is a subset of~$\rmC_A(X)$ obtained from~$Y$ by applying a {\em finite} sequence of
{\em simple expansions}, where a simple expansion replaces one element~\hbox{$y\in Y$} by the elements~\hbox{$\{y\}\x A\subset \rmC_A(X)$}, its ``descendants'' in~$\rmC_A(X)$.  For~\hbox{$Y\subset\rmC_A(X)$}, we denote by~$\calE(Y)$ the set of all expansions of~$Y$. 
\end{definition}

\begin{figure}
\centering
\includegraphics[width=0.7\textwidth]{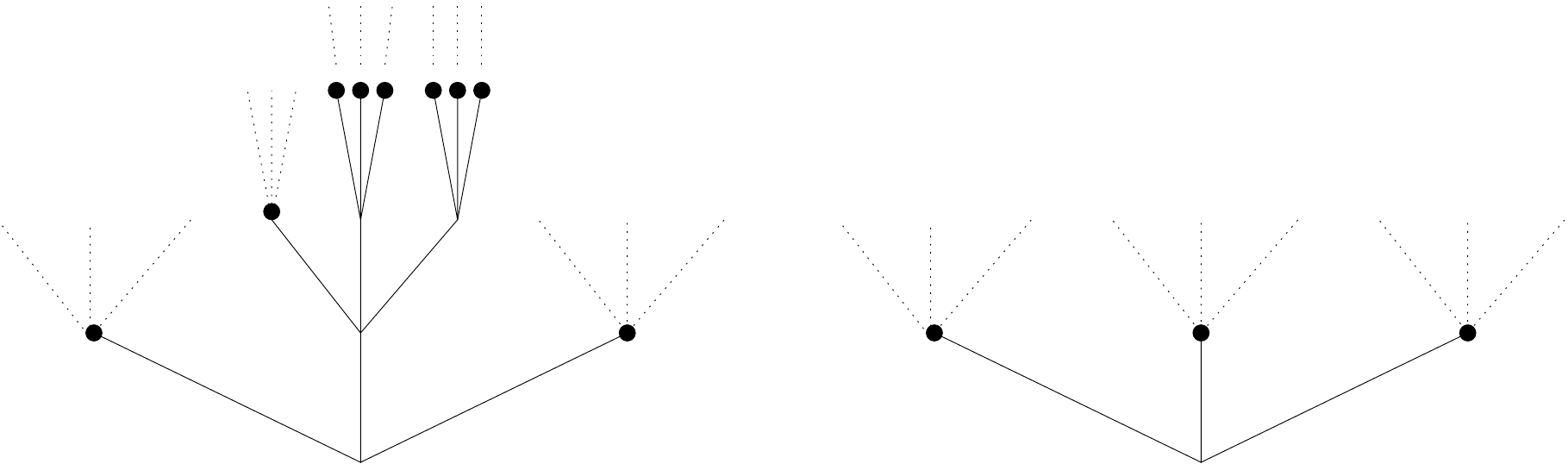}
\caption{Recall the identification of~$\rmC_{A}[2]^+$ for~$|A|=3$ with the set of vertices of the forest of  Figure~\ref{fig:forest}. 
Under this identification, the vertices marked by black dots here define an expansion of the standard basis~$[2]$, obtained from it by applying a sequence of five simple expansions. This expansion has cardinality~\hbox{$12=2+5(3-1)$}.}\label{fig:basis}
\end{figure}

If~$Y$ was a basis, then so is any of its expansions \cite[Lem.~2.3]{Higman}. In particular, all the expansions of the canonical basis~$X$ represent bases for~$\rmC_A(X)$, and these are the bases we will work with. 
Note that any such basis is a finite  subset of~$\rmC_A^+(X)$. 
If we think of~$\rmC_A^+(X)$ as the vertices of an infinite~$|A|$--ary forest with roots the elements of~$X$ (as in Figure~\ref{fig:forest}), an expansion of~$X$ is the set of leaves of finite~$|A|$--ary subforest~$F$ which~``generates'' in the sense that the infinite forest is the union of~$F$ and the infinite trees attached to the leaves of~$F$. (See 
Figure~\ref{fig:basis} for an example.) %and for instance~\cite[Sec.~4]{Bro87} for more details on this point of view.)

\begin{lem}\label{lem:expord}
The set of finite bases~$S$ of~$\rmC_A(X)$ satisfying that~$S\subset \rmC_A^+(X)$ identifies with the set~$\calE(X)$ of all expansions of~$X$. It is a partially ordered set with the relation 
$$Y\le Z \ \ \ \Longleftrightarrow\ \ \  Z \ \textrm{is an expansion of}\ \ Y.$$  
The poset~$\calE(X)$ has a least element, namely~$X$. 
\end{lem}

%Note that for~$Y,Z\in \calE(X)$, 
%$$Y\le Z \ \ \Longleftrightarrow\ \ C_A^+(Y)\supset Z \ \ \Longleftrightarrow\ \ C_A^+(Y)\supset C_A^+(Z).$$ 

\begin{proof}
We first check that the given relation~$\le$ defines a partial ordering on~$\calE(X)$: transitivity follows directly from the definition of expansion and antisymmetry follows using in addition the fact that if~$Z$ is an expansion of~$Y$, then~$|Z|\ge |Y|$, with strict inequality if the extension is non-trivial. We have that~$X\in \calE(X)$ and, by definition,~$X\le Y$ for any~$Y\in \calE(X)$, so~$X$ is a least element.  

The fact that all expansions of~$X$ are bases is given by Lemma~2.3 of \cite{Higman}, and these, by definition, lie inside~$\rmC_A^+(X)$. The fact that any finite basis that  is a subset of~$\rmC_A^+(X)$ is an expansion of~$X$ follows from Lemma~2.4 in \cite{Higman}: Suppose~$S$ is such a  basis, and let~\hbox{$U=\rmC_A^+(S)\subset\rmC_A^+(X)$}. Then~\hbox{$U=\rmC_A^+(S)\cap\rmC_A^+(X)$}, hence satisfies condition~(i) in the lemma, which is equivalent to condition~(iii) in the lemma, that is~\hbox{$U=\rmC_A^+(Z)$} for~$Z$ some expansion of~$X$. Hence~$\rmC_A^+(Z)=\rmC_A^+(S)$ as~$A^*$--subsets of~$\rmC_A^+(X)$, which is only possible if~$S=Z$ because~$S$ and~$Z$ both generate this free~$A^*$--set. 
%Then~$S\le Z$ and~$Z\le S$, and hence they must be equal. 
\end{proof}

In Section~\ref{sec:spaces} and \ref{sec:connectivity}, we will work with independent sets, which are subsets of expansions:

\begin{definition}\label{def:independent}
A finite subset~$P\subset \rmC_A^+(X)$ is called {\em independent} if there exists an expansion~$E$ of~$X$ such that~$P\subset E$. We denote by~$\calI(X)$ the set of all independent sets of~$\rmC_A(X)$. 
If~$P,Q$ are independent sets, we will say that~$P$ is {\em independent from~$Q$} if they are disjoint and~$P\sqcup Q$ is still independent. 

We have that~$\calE(X)\subset \calI(X)$. We call the elements of~$\calI_0(X):=\calI(X)\minus \calE(X)$ the {\em non-generating independent sets}, so that the elements of~$\calI_0(X)$ are precisely the independent sets that are not bases.
\end{definition} 

When $X=[r]$, we write $\calE[r],\, \calI[r]$ and $\calI_0[r]$ for $\calE(X),\, \calI(X)$ and $\calI_0(X)$.

\begin{lem}\label{lem:minexp}
Let~$P\in \calI(X)$ be an independent set. The subposet 
of~$\calE(X)$ of expansions of~$X$ containing~$P$ has a least element. 
\end{lem}

\begin{proof}
Suppose~$E=Q\sqcup P$ is an expansion of~$X$ containing~$P$. Consider the subalgebra of~$\rmC_A(X)$ generated by~$Q$.  
By \cite[Lem.~2.7(ii),(iii)]{Higman}, this subalgebra has a finite generating set~$G\subset \rmC_A^+(X)$ with the property that any element of its intersection with~$\rmC_A^+(X)$ is an expansion of an element of~$G$. In particular,~$Q$ is necessarily an expansion of~$G$, and~$G\sqcup P$ is the requested least element. 
\end{proof}

\begin{figure}
\centering
\includegraphics[width=0.7\textwidth]{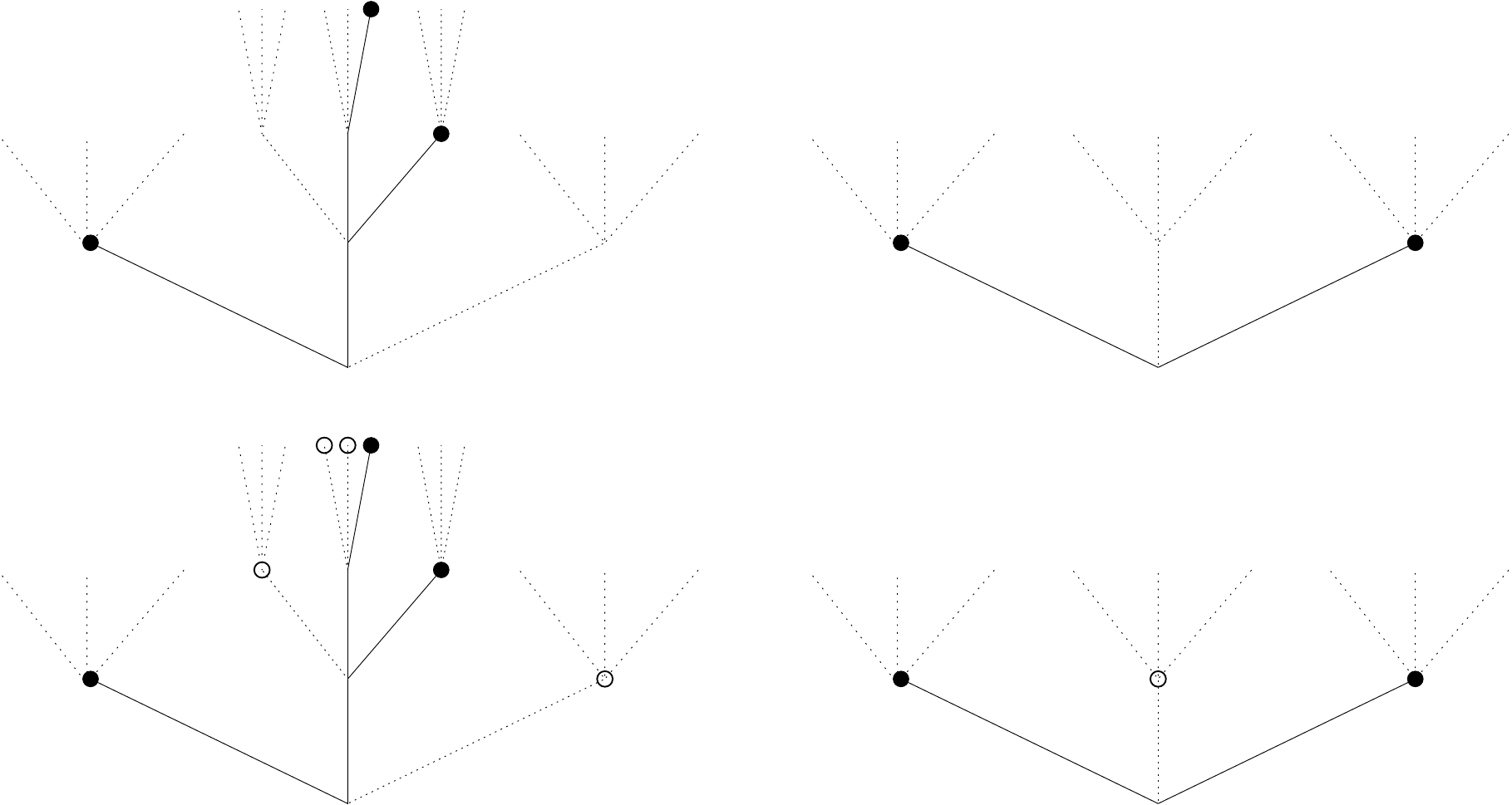}
\caption{The black dots define an independent set in~$\rmC_{A}[2]^+$ for~$|A|=3$, identified with the set of vertices of the forest of  Figure~\ref{fig:forest}, and the white dots complete it to the least expansion of~$[2]$ containing it.}\label{fig:indepset}
\end{figure}

Higman proved the following important fact about bases: 

\begin{lem}{\upshape \cite[Cor.~1]{Higman}}\label{lem:common} %p12 in Higman
Any two finite bases of a free Cantor algebra have a common expansion.
\end{lem}

In particular, for any~$X$, the poset~$\calE(X)$ is directed, that is any pair of elements in~$\calE(X)$ have a common upper bound. 

A consequence of the lemma is that the cardinality of a finite basis for~$\rmC_{A}(X)$ is congruent to~$|X|$ modulo~\hbox{$|A|-1$}. In fact, for two finite sets~$X$ and~$Y$, we have that~$\rmC_{A}(X)\cong\rmC_{A}(Y)$ if and only if~\hbox{$X=\emptyset=Y$} or if~$X$ and~$Y$ are both non-empty of cardinality congruent modulo~$|A|-1$; this is the condition that guarantees that the two Cantor algebras admit finite bases of the same cardinality.

If~$X=X_1\sqcup X_2$ is the disjoint union of two finite sets, we have a canonical isomorphism
$$\rmC_A^+(X) \ \cong\ \rmC_A^+(X_1)\sqcup \rmC_A^+(X_2)$$
and more generally, the set~$\rmC_A^+(X)$ splits as a disjoint union
$$\rmC_A^+(X)=\bigsqcup_{x\in X}\rmC_A^+(\{x\}).$$
In terms of the forests of Figure~\ref{fig:forest}, we see that~$\rmC_A^+(X)$ is a disjoint union of trees, one for each element of~$X$. 
Expansions of~$X$ are subsets of~$\rmC_A^+(X)$. The following result says that the property of being an expansion can be checked componentwise.

%Note that if~$Y$ is an expansion of~$X$, then expansions of~$Y$ are also expansions of~$X$, as~$\rmC_A^+(Y)$ is naturally a subset of~$\rmC_A^+(X)$.
\begin{lem}\label{lem:exp} Let~$S\subset \rmC_{A}^+(X\sqcup Y)$ be a finite subset. Then 
$S\in \calE(X\sqcup Y)$ if and only if~$S\cap \rmC_A^+(X)\in \calE(X)$ and~\hbox{$S\cap \rmC_{A}^+(Y)\in \calE(Y)$}, where we consider~$\rmC_A^+(X)$ and~$\rmC_A^+(Y)$ as subsets of~$\rmC_A^+(X\sqcup Y)$ through
 the identification~$\rmC_A^+(X\sqcup Y)\cong \rmC_A^+(X)\sqcup \rmC_A^+(Y)$. 
\end{lem}

\begin{proof}
Suppose first that~$S$ is an expansion of~$X\sqcup Y$, so~$S$ is obtained from~$X\sqcup Y$ by applying a finite sequence of simple expansions. Now each simple expansion is either expanding an element of~$\rmC_A^+(X)$ or of~$\rmC_A^+(Y)$, and we see that~$S=S_0\sqcup S_1$, with~$S_0$ the subset obtained by applying to~$X$ the expansions of the first type, and~$S_1$ the subset applying to~$Y$ the expansions of the second type. As~$S_0=S\cap \rmC_A^+(X)$ and~$S_1=S\cap \rmC_A^+(Y)$, this gives the first direction. The reverse direction is direct from the definition. 
%This follows from the fact that if~$E$ is an expansion of~$X$ and~$F$ an expansion of~$Y$, then~$E\cup F$ is an expansion of~$X\cup Y$, and it does not matter whether we first expand~$X$ and then~$Y$ or do it the other way around. Hence we can write any expansion of~$X\cup Y$ as the composition of an expansion of~$X$ and an expansion of~$Y$, and two such expansions are equal as expansions of~$X\cup Y$ if and only if they are equal as expansions of~$X$ and~$Y$ separately. 
\end{proof}

In fact, the lemma can be used to show that 
$$\calE(X)\cong \prod_{x\in X}\calE(\{x\})$$
as a poset. 

\subsection{Representing isomorphisms}\label{sec:isos}

Given a basis~$E$ of~$\rmC_A(X)$, a basis~$F$ of~$\rmC_A(Y)$ and a bijection~$\la\colon E\to F$, there is a unique isomorphism of Cantor algebras~$f\colon\rmC_A(X)\to \rmC_A(Y)$  that satisfies that~$f|_{E}=\la$. (Figure~\ref{fig:iso} gives an example in terms of forest.) 
But such a representing triple~$(E,F,\la)$ is far from unique. 
Indeed, any expansion~$E'$ of~$E$ defines a new such triple~$(E',F',\la')$ representing the same isomorphism~$f$ simply by taking 
$F'=f(E')$ and~$\la'=f|_{E'}$. We will now see that any isomorphism can be represented by a triple~$(E,F,\la)$ where~$E\in\calE(X)$ and~\hbox{$F\in\calE(Y)$}, and that there is in fact a canonical representative among such triples. To describe this canonical representative, we use the following partial ordering on the set of all such triples:

\begin{figure}
\centering
\includegraphics[width=0.4\textwidth]{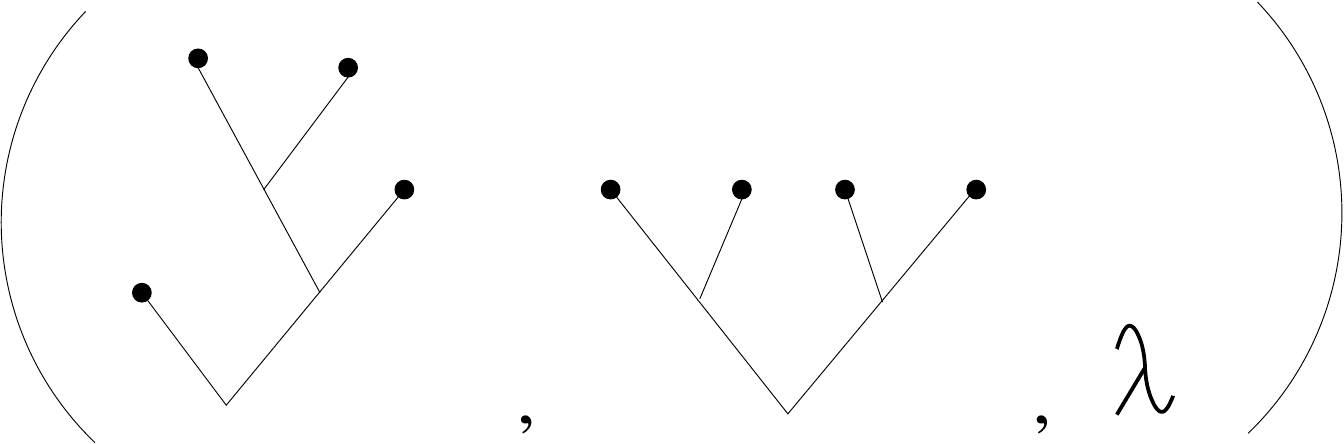}
\caption{Representation of an automorphism of~$\rmC_A[1]$ with~$|A|=2$ in terms of a triple~$(E,F,\la)$, where~$E$ is the set of black dots in the first tree,~$F$ the set of black dots in the second tree, and~$\la$ is some bijection between these two sets.}\label{fig:iso}
\end{figure}

\begin{Def}
For~$f\colon\rmC_A(X)\to \rmC_A(Y)$ an isomorphism of Cantor algebras, define
$$\Rep(f)=\{(E,F,\la) \ |\  E\in \calE(X), F\in \calE(Y) \ \textrm{and}\  \la=f|_E\colon E\stackrel{\cong}\rar F\}$$
to be the set of triples representing~$f$ with the property that~$E$ and~$F$ are expansions of~$X$ and~$Y$ respectively, with poset structure 
%$(E,F,\la)$,~$(E',F',\la')$ be two triple representing the same map~$f\colon\rmC_A(X)\to \rmC_A(Y)$, with~$E,E'\subset \rmC_A^+(X)$ and~$F,F'\in \rmC_A^+(Y)$. 
%We define the  
$$(E,F,\la)\le (E',F',\la') \ \  \ \Longleftrightarrow\ \ \ E\le E',$$ where the right hand side uses the partial order relation on the set~$\calE(X)$ of expansions of~$X$ 
from Lemma~\ref{lem:expord}.
\end{Def}

Note that this defines indeed a partial ordering on~$\Rep(f)$:  the transitivity follows from transitivity in the poset~$\calE(X)$, as does anti-symmetry  once one notices that~$E=E'$ forces~$F=F'$ and~$\la=\la'$ given that the triples represent the same morphism. 

Note also that having the relation~$(E,F,\la)\le (E',F',\la')$ in~$\Rep(f)$ forces the relation~$F\le F'$  in~$\calE(Y)$, because~\hbox{$F=f(E)$} and~\hbox{$F'=f(E')$}, and~$f$ takes expansions to expansions. Likewise, the map~$\la'$ under such a condition is necessarily the map induced by~$\la$ on~$E'$. 
Also, as~$F$ and~$\la$ are determined by~$E$ and~$f$, the forgetful map~$\Rep(f)\to \calE(X)$ that takes a triple~$(E,F,\la)$ to~$E$ is injective, so~$\Rep(f)$ canonically embeds as a sub-poset of~$\calE(X)$.

The main result of the section is the following: 

\begin{prop}\label{triples}
For any isomorphism of Cantor algebras~$f\colon\rmC_A(X)\to \rmC_A(Y)$, the poset~$\Rep(f)$ is non-emtpy and has a least element. 
\end{prop}

In other words, any isomorphism~$f\colon\rmC_A(X)\to \rmC_A(Y)$ can be represented in the above sense by a triple~$(E,F,\la)$ with~$E$ an expansion of~$X$,~$F$ an expansion of~$Y$, and with~$\la\colon E\to F$ a bijection, and  there is a unique minimal such representative, with~$E$ minimal with respect to the ordering on~$\calE(X)$.

The result is a mild generalization of~\cite[Lem.~4.1]{Higman} who gives the case~$X=Y$, and it follows from the same proof. We give it for completeness as the result is crucial for us. 

\begin{proof}
Consider the subset~$U$ of~$\rmC_A^+(X)$ defined by 
\begin{align*}U \ =\ \rmC_A^+(X)\ \cap\ f^{-1}(\rmC_A^+(Y)) %\ \cap\ \dots\cap\  f_1^{-1}\circ \dots\circ f_k^{-1}\rmC_A^+(X_k)\\
\ =\ \rmC_A^+(X)\ \cap\ \rmC_A^+(f^{-1}(Y)). %\ \cap\ \dots  \ \cap\ \rmC_A^+(f_1^{-1}\circ\dots\circ f_k^{-1}(X_k)).
\end{align*}
By \cite[Lem.~2.4]{Higman}, there exists an expansion~$E$ of~$X$ such that~$U=\rmC_A^+(E)$. As~$E$ lies inside~$U$, it satisfies that %$f_i\circ\dots\circ f_1(Y_0)$ 
$F=f(E)$ is an expansion of~$Y$. %for every~$i$. 
Taking~$(E,F,\la)$ with~$\la=f|_{E}$ yields a presentation of~$f$ with~$E$ and~$F$ expansions of~$X$ and~$Y$, which shows that~$\Rep(f)$ is non-empty. 

We also claim that the triple~$(E,F,\la)$ just constructed is in fact the least element of~$\Rep(f)$. Indeed, assume that~$E'$ is another expansion of~$X$ satisfying that~$f(E')$ 
%$f_i\circ\dots\circ f_1(Z)$
 is an expansion of~$Y$, and let~$(E',f(E'),\la')\in \Rep(f)$ be the associated triple. 
%for every~$i$. 
Then~$E'$ must lie in~$U=\rmC_A^+(E)$ and hence be an expansion of~$E$ by Lemma~\ref{lem:expord} as it is a basis of~$\rmC_A(X)$ and hence also of~$\rmC_A(E)$. 
It follows that~\hbox{$(E,F,\la)\le (E',f(E'),\la')$} in~$\Rep(f)$. 
%This proves the first part of the proposition. 
%The last statement then follows by picking~$Y_0$ to be the unique minimal expansion of~$X_0$ given by the first part of the proposition,~$Y_1$ its image by~$f$, and~$\la$ the restriction of~$f$ to~$X_0$. 
\end{proof}

\begin{ex}\label{ex:iso}
Any bijection~$\la\colon A\stackrel{\cong}\rar [n]$ induces an isomorphism~$\rmC_A[1]\stackrel{\cong}\rar \rmC_A[n]$ which is represented by the triple~$(\{1\}\x A,[n],\la)$. Its inverse~$\rmC_A[n]\to \rmC_A[1]$ is represented by the triple~$([n],\{1\}\x A,\la^{-1})$. More generally,  if $f$ is (minimally) represented by $(E,F,\la)$, then its inverse $f^{-1}$ is (minimally) represented by~$(F,E,\la^{-1})$. 
\end{ex}

\begin{rem}\label{rem:comptriples}
The composition of isomorphisms of Cantor algebras in terms of representatives can be computed as follows. If~$f\colon\rmC_A(X)\to \rmC_A(Y)$ is represented by~$(E,F,\la)$ and  
$g\colon\rmC_A(Y)\to \rmC_A(Z)$ is represented by~$(G,H,\mu)$, we have a diagram 
$$\xymatrix@R-1pc{\hat E \ar@{}[d]|\IV \ar[r]^-{\hat \la} & \hat F\ =\ \hat G \ar@{}[d]|{\IV\ \ \ \ \ \ \ \IV} \ar[r]^-{\hat \mu} & \hat H \ar@{}[d]|\IV \\
E \ar@{}[d]|\IV \ar[r]^-{\la} & F\ \ \ \ \ \  \  \  G \ar@{}[d]|{\IV\ \ \ \ \ \ \ \IV} \ar[r]^-{\mu} & H \ar@{}[d]|\IV \\
X & Y \ =\ Y  & Z
}$$
for~$\hat F=\hat G$ a common expansion of~$F$ and~$G$ (which exists by Lemma~\ref{lem:common}), and~$\hat \la$ and~$\hat \mu$ the maps induced by~$\la$ and~$\mu$, or equivalently the restrictions of~$f$ and~$g$ to~$\hat E$ and~$\hat G$, with~$\hat E:=f^{-1}(\hat F)$ and~$\hat H:=g(\hat G)$. Then the triple~$(\hat E,\hat H,\hat \mu\circ \hat \la)$ represents the composition~$g\circ f$.  Note that, even if the original representing triples were minimal, and~$\hat G$ is chosen minimally, the resulting triple representing the composition will in general not be minimal, as can be checked for instance in Example~\ref{ex:iso}. 
\end{rem}

\subsection{Categories of Cantor algebras}\label{sec:catcan}

We will in this paper work with {\em  permutative categories}, which are symmetric monoidal categories that are strictly associative and have a strict unit. By a {\em strict monoidal functor}, we will mean a monoidal functor~$F$ such that the morphisms~$F(x)\oplus F(y)\to F(x\oplus y)$ are the identity. 

Let~$\Set$ denote the category with objects the natural numbers, where we identify the integer~$r\ge 0$ with the set~\hbox{$[r]=\{1,\dots,r\}$}, and with morphisms the maps of sets. We denote by~$\Set^\x$ its subcategory of isomorphisms~(the bijections). The categories~$\Set$ and~$\Set^\x$ are both permutative categories with the monoidal structure~$\oplus$ defined using the sum on objects, and disjoint union on morphisms, using the canonical identification~$[r]\sqcup [s]\cong [r+s]$. The unit is the empty set~$[0]=\emptyset$ and the symmetry~$[r]\oplus [s]=[r+s]\to [r+s]=[s]\oplus [r]$ is given by the~$(r,s)$ block permutation. 

Let~$A=\{a_1,\dots,a_n\}$ as before. 
We now define our category~$\Can_A$ of finitely generated free Cantor algebras of type~$A$.  
To avoid set-theoretical issues, we will only consider the Cantor algebras freely generated by the sets~$[r]$ for~$r\ge 0$. So the objects of~$\Can_A$ are the natural numbers just like for the category~$\Set$, but with~$r$ now identified with the Cantor algebra~$\rmC_A[r]$. The morphisms in~$\Can_A$ are the morphisms of Cantor algebras as in Definition~\ref{def:Cantor}. We will denote by~$\Can_A^\x$ the subcategory  of isomorphisms in~$\Can_A$. 

%\begin{rem}\label{rem:n=1}
%We could have considered the case of a set~$A$ of cardinality~$1$. A Cantor algebra of type~$[1]$ is 
% \nnote{say something somewhere about the case~$A=[1]$} 
%\end{rem}

Taking the free Cantor algebra on a given set induces a  functor 
$$\rmC_A\colon\Set\rar \Can_A$$
which is the identity on objects.
Indeed, any map of sets~$[r]\to [s]$ induces a morphism~\hbox{$\rmC_A[r]\to \rmC_A[s]$} between the free Cantor algebras, because~$\rmC_A[r]$ is free on~$[r]$ and~$[s]$ canonically identifies with a subset of~$\rmC_A[s]$. This association is compatible with composition. 

\begin{Def}
For~$r<s$, we will denote by 
$$i_L\colon[r]\inc [s] \ \ \ i_R\colon[r]\inc [s]$$
the {\em left} and {\em right embeddings} of~$[r]$ into~$[s]$, i.e.~that  of~$[r]$ as the first~$r$ (resp.~last)~$r$ elements of~$[s]$. We will likewise denote by 
$$i_L,i_R\colon \rmC_A[r] \rar \rmC_A[s]$$
the corresponding induced maps~$\rmC_A(i_L)$ and~$\rmC_A(i_R)$. 
\end{Def}

We use the functor~$\rmC_A$ and the permutative structure of~$\Set$ to define a permutative structure, also denoted~$\oplus$, on 
$\Can_A$: 
On objects, we define 
\[
\rmC_A[r]\oplus \rmC_{A}[s]:=\rmC_{A}[r+s]
\]
and for morphisms~$f\colon\rmC_A[r]\to \rmC_A[r']$ and~$g\colon\rmC_A[s]\to \rmC_A[s']$, we define~$$f\oplus g\colon\rmC_A[r+s]\to \rmC_A[r'+s']$$  
to be the unique morphism defined on the basis~$[r+s]=i_L[r]\sqcup i_R[s]$ using the map 
$$[r]\ \stackrel{f|_{[r]}}\rar \ \rmC_A[r']\ \stackrel{i_L}\rar \ \rmC_A[r'+s']$$ on the first~$r$ elements and
$$[s]\ \stackrel{g|_{[s]}}\rar\  \rmC_A[s']\ \stackrel{i_R} \rar \ \rmC_A[r'+s']$$ on the last~$s$ elements. 
%where~$i_L\equiv \rmC_A(i_L)$ and~$i_R\equiv\rmC_A(i_R)$
%are the maps induced by the inclusions of sets~$i_L:[r']\inc [r'+s']$ and~$i_R:[s']\inc [r'+s']$. 
%the maps~$\rmC_A[r']\to \rmC_A[r'+s']$ and~$\rmC_A[s']\to \rmC_A[r'+s']$ are induced from the maps of sets~$[r'],[s']\to [r'+s']$. 
Finally, let 
\[
\s_{r,s}\colon\rmC_A[r+s]=\rmC_A([r]\oplus [s])\rar \rmC_A([s]\oplus [r])=\rmC_A[r+s]
\] 
be the image under the functor~$\rmC_A$ of the symmetry~$[r]\oplus [s]\to [s]\oplus [r]$ in the category~$\Set$.

\begin{prop}\label{injectiveAut}
The sum~$\oplus$ and symmetry~$\s_{r,s}$ defined above make~$\Can_A^\x$  into a  permutative category, with unit~$\rmC_A[0]=\emptyset$, with the property that 
 the free Cantor algebra functor~$\rmC_A\colon\Set^\x\to \Can_A^\x$  is a strict symmetric monoidal functor, and that the sum 
\[
\Can^\x_A(\rmC_A[r],\rmC_A[r'])\x \Can_A^\x(\rmC_A[s],\rmC_A[s'])\stackrel{\oplus}{\longrightarrow} \Can_A^\x(\rmC_A([r+s],\rmC_A[r'+s']))
\]
is injective. 
\end{prop}

One can likewise show that~$\Can_A$ is a permutative category. We restrict to~$\Can^\x_A$ for simplicity, as this is the part that is relevant to us. 

Before proving the result, we interpret the sum $\oplus$ in terms of representatives.

We can reinterpret Lemma~\ref{lem:exp} as saying that the sum of sets~$[r]\oplus [s]=[r+s]$ extends to a sum of expansions: 
\[
\oplus\colon\calE[r]\x \calE[s]\rar \calE[r+s], 
\]
formally defined by setting~$E\oplus F=i_L(E)\sqcup i_R(F)$. (By the same lemma, this sum is an isomorphism.) Using this sum operation, we have the following: 

\begin{lem}\label{lem:sumrep}
Let~$f\in \Can_A^\x(\rmC_A[r],\rmC_A[r'])$ and~$g\in \Can_A^\x(\rmC_A[s],\rmC_A[s'])$ be isomorphisms represented by~$(E,F,\la)\in \Rep(f)$ and~$(G,H,\mu)\in \Rep(g)$. Then~\hbox{$(E\oplus G,F\oplus H,\la\oplus \mu)$} represents the sum~\hbox{$f\oplus g$}. Moreover, if~$(E,F,\la)$ is the least element of~$\Rep(f)$ and~$(G,H,\mu)$ the least element of~$\Rep(g)$, then the least element of~$ \Rep(f\oplus g)$ is~\hbox{$(E\oplus G,F\oplus H,\la\oplus \mu)$}.
\end{lem}

\begin{proof}
The sum~$f\oplus g$ is defined as the unique map taking~$i_L[r]\subset [r+s]$ to~$\rmC_A[r'+s']$ using~$i_L\circ f$ and~\hbox{$i_R[s]\subset [r+s]$} using~$i_R\circ g$. 
Now suppose that~$f$ is represented by~$(E,F,\la)$ and~$g$ by~$(G,H,\mu)$. The map of Cantor algebras~$h\colon\rmC_A[r+s]\to \rmC_A[r'+s']$ represented by~$(E\oplus G,F\oplus H,\la\oplus \mu)$ by definition takes~$i_L(E)$ to~$i_L(F)$ using~$\la=f|_{E}$ and~$i_L(G)$ to~$i_L(H)$ using~$\mu=g|_{G}$. As~$i_L(E)$ is an expansion of~$i_L[r]$ with~$i_L(F)$ the corresponding expansion of~$i_L(f[r])$, and the induced map~$h$ is a map of Cantor algebras, we necessarily have that~$h|_{i_L[r]}=i_L\circ f$, and likewise,~$h|_{i_R[s]}=i_R\circ g$. Hence~$h=f\oplus g$.

The sum respects the poset structure by Lemma~\ref{lem:exp}. Now suppose that~$(E,F,\la)$ and~$(G,H,\mu)$ are  least elements, and that there exists~$(S,T,\nu)<(E\oplus G,F\oplus H,\la\oplus \mu)$ strictly smaller in~$ \Rep(f\oplus g)$.
By Lemma~\ref{lem:exp}, the subset~$S_0:=S\cap \rmC^+_A[r]$ is an expansion of~$[r]$. As~
$$S_0\subset \rmC^+_A[r]\subset\rmC^+_A[r+s] \ \cong \rmC_A^+[r]\sqcup \rmC_A^+[s]$$
is an expansion of~$i_L[r]$ inside~$\rmC_A[r+s]$, we have~$(f\oplus g)(S_0)=i_L\circ f(S_0)$ inside~$\rmC_A[r'+s']$, as~$f\oplus g$ restricts to~$f$ on this subset. It follows that~$(S_0,f(S_0),f|_{S_0})$ is also a representative of~$f$ which is smaller or equal to~$(E,F,\la)$. Likewise, setting~\hbox{$S_1:=S\cap \rmC^+_A[r]\subset\rmC^+_A[r+s]$}, we get a 
representative~$(S_1,g(S_1),g|_{S_1})$ of~$g$ which is smaller or equal to~$(G,H,\mu)$. Given that~\hbox{$S=S_0\oplus S_1$}, if~$S< E\oplus G$, we must have that either~$S_0<E$ or~\hbox{$S_1<G$}, contradicting the minimality assumption. Hence the sum~$\oplus$ takes least elements to least elements. 
\end{proof}

\begin{proof}[Proof of Proposition~\ref{injectiveAut}]
By the lemma, we can write the sum 
$$ \Can^\x_A(\rmC_A[r],\rmC_A[s]) \x \Can_A^\x(\rmC_A[r'],\rmC_A[s'])  \ \ \stackrel{\oplus}\rar\ \ \ \Can_A^\x(\rmC_A[r+r'],\rmC_A[s+s'])~$$
in terms of (minimal) representatives. 
Functoriality follows then from the fact that~$\rmC_A^+[r+r']\cong \rmC_A^+[r]\sqcup \rmC_A^+[r']$ and likewise for~$s,s'$, and the fact that composition can be computed componentwise. 
Associativity of the sum is likewise easily checked in this description of the sum, and~$\rmC_A[0]$ is a strict unit. For the symmetry, we need to check that 
$$\xymatrix{\rmC_A[r]\oplus \rmC_A[s] \ar[r]^-{f\oplus g} \ar[d]_{\s_{r,s}} & \rmC_A[r']\oplus \rmC_A[s'] \ar[d]_{\s_{r',s'}} \\
\rmC_A[s]\oplus \rmC_A[r] \ar[r]^-{g\oplus f} & \rmC_A[s']\oplus \rmC_A[r']   }$$
commutes for any morphism~$f,g$. In terms of representatives, one checks that if~$f$ is represented by~$(E,F,\la)$ and~$g$ by~$(G,H,\mu)$, then both compositions will be represented by 
$$(E\oplus G,H\oplus F,\widehat \s_{r',s'}\circ (\la\oplus \mu))= (E\oplus G,H\oplus F, (\mu\oplus \la)\circ \widehat \s_{r,s})$$
for~$\widehat \s_{r,s}$ the restrictions of the symmetry~$\s_{r,s}\in\Aut(\rmC_A[r+s])$ to~$E\oplus G$ and~$\widehat \s_{r',s'}$ the corresponding restriction of~$\s_{r',s'}\in\Aut(\rmC_A[r'+s'])$ to~$\la(F)\oplus \mu(G)$.  So the square commutes and 
$(\Can_A,\oplus,\s,\rmC_A[0])$ is a permutative category. 

Finally, injectivity of the sum follows, also using representatives,  from the fact that an automorphism~$f$ is uniquely represented by its least representative in~$\Rep(f)$ (Proposition~\ref{triples}), that sums of minimal representatives are minimal representatives (Lemma~\ref{lem:sumrep}),
and that the map is injective on representatives. 
\end{proof}

The following proposition will be useful in Section~\ref{sec:spaces}. 

\begin{prop}\label{prop:split}
Let~$g\colon\rmC_A[r+s]\to \rmC_A[r+s]$ be an isomorphism and suppose that~$g$ restricts to the identity on a given finite basis of~$\rmC_A[s]\equiv i_R\rmC_A[s]\subset \rmC_A[r+s]$. 
Then we have~$g=g_0\oplus \rmC_A[s]$ for some isomorphism~\hbox{$g_0\colon\rmC_A[r]\to \rmC_A[r]$}. 
\end{prop}

Here and in the following we often employ the Milnor--Moore notation and denote the identity morphism of an object by that object.

\begin{proof}
Let~$(F,G,\la)\in \Rep(g)$ be a representative of~$g$ and let~$E\subset i_R\rmC_A[s]$ be a finite basis of~$\rmC_A[s]$ fixed by~$g$. As~$F$ is an expansion of~$[r+s]$, we have that~$F_1:=F\cap i_R\rmC_A^+[s]$ is an expansion of~$[s]$ (Lemma~\ref{lem:exp}). Now~$F_1$ and~$E$ have a common expansion~$\hat E$ by Lemma~\ref{lem:common}, which is an expansion of~$[s]$ as~$F_1$ was an expansion of~$[s]$. Let~$\hat F=(F\cap i_L\rmC_A^+[r])\cup \hat E$ and~$\hat G=g(\hat F)$ and~$\hat \la=g|_{\hat F}$. 
Then~$(\hat F,\hat G,\hat \la)$ is also a representative of~$g$. Note now that~$g(\hat E)=\hat E$ and~$\hat\la|_{\hat E}=\id$ as~$g$ respects expansions and restricted to the identity on~$E$. 
%\Mnote{Is that a sentence? N: yes.}
It follows that~$\hat G=(\hat G\cap i_L\rmC_A^+[r])\cup \hat E$ for~$(\hat G\cap \rmC_A[r])$ an expansion of~$[r]$. Hence~$g=g_0\oplus \rmC_A[s]$ for~$g_0$ represented by~$(\hat F_0=F\cap i_L\rmC_A^+[r],\hat G_0=\hat G\cap i_L\rmC_A^+[r], \hat\la|_{\hat F_0})$, which proves the result.  
\end{proof}

\section{Spaces associated to Higman--Thompson groups}\label{sec:spaces}

This section and the  next one are concerned with the proof of homological stability for the Higman--Thompson groups~$\rmV_{n,r}$, the automorphism group of the free Cantor algebra~$\rmC_{A}[r]$, with respect to the number~$r$ of generators, where~$A=\{a_1,\dots,a_n\}$ as before. 
%For this part of the result, it is enough to consider the Cantor algebras~$\rmC_{n,r}$. Throughout both sections,~$n\geqslant 2$ will be fixed.  
Given a family of groups satisfying a few properties, the 
paper  \cite{Randal-Williams+Wahl}  yields  a sequence of spaces whose high connectivity implies homological stability for the family of groups. In this section, we will show how 
the groups~$\rmV_{n,r}$ fit in the framework of \cite{Randal-Williams+Wahl} and construct spaces relevant to the proof of homological stability. %, which will be given in the following section.
% as well as closely related spaces that will allow us, in the following section, to prove the necessary connectivity result. 
%Given that the variables~$n$ and~$r$ play a rather different role, we will keep on using~$$A=[n]=\{1,\dots,n\}$$ for the set with~$n$ elements. 

For a fixed type~$A$ we collect the Higman--Thompson groups into a groupoid~$\bfV_A$, which is a subgroupoid of the category~$\Can_A$, or in fact of its groupoid of isomorphisms~$\Can_A^\x$ defined in the previous section.  The objects of~$\bfV_A$ are the same as those of~$\Can_A^\x$, namely the natural numbers~$r$ as placeholders for the free Cantor algebras~$\rmC_A[r]$, and the morphism sets are defined by setting 
\begin{align*}
\bfV_A(\rmC_A[r],\rmC_A[s])=\left\{ \begin{array}{ll}  \Can^\x_A(\rmC_A[r],\rmC_A[s]) & r=s\\
\emptyset & r\neq s
\end{array}\right.
\end{align*}
In other words,~$$\bfV_A\cong \bigsqcup_{r\ge 0}\rmV_{n,r}$$ is the groupoid of {\em automorphisms} in~$\Can_A^\x$, with each group~$\rmV_{n,r}$ considered as a groupoid with one object. Recall that there are isomorphisms~$\rmC_{A}[r]\cong\rmC_A[r+(n-1)]$ for any~$r\geqslant1$ (with~$n=|A|$). These isomorphisms are morphisms in the groupoid~$\Can^\x_A$  but we do {\em not} include these isomorphisms into the groupoid~$\bfV_A$. 
Note that the symmetric monoidal structure of~$\Can_A^\x$ restricts to~$\bfV_A$, making it a permutative groupoid. 

\subsection{Homogeneous categories}\label{sec:homogeneous}

Recall from~\cite[Def.~1.3]{Randal-Williams+Wahl} that a monoidal category~$(\bfH,\oplus,0)$ is called {\em homogeneous} if the monoidal unit~$0$ is initial and if the following two conditions are satisfied for every pair of objects~$X, Y$ in~$\bfH$: The set~${\bfH}(X,Y)$ is a transitive~$\Aut_{\bfH}(Y)$--set under post-composition, and the homomorphism
\[
\Aut_{\bfH}(X)\to\Aut_{\bfH}(X\oplus Y)
\]
that takes~$f$ to~$f\oplus Y$ is injective with image~$\{\phi\in\Aut_{\bfH}(X\oplus Y)\ |\phi\circ (\iota_X\oplus Y)=\iota_X\oplus Y\}$. %, the stabilizer of the image of~$\iota_A$.
The main examples of homogeneous categories can be obtained by applying Quillen's bracket construction to a groupoid. We recall this construction here, and show that it yields a homogeneous category when applied to the groupoid~$\bfV_A$. 

Let~\hbox{$\bfQ_A=\langle\bfV_A,\bfV_A\rangle$} denote the category obtained by applying Quillen's bracket construction~(see~\cite[p.~219]{Grayson:QuillenII}) to the groupoid~$\bfV_A$.  The category~$\bfQ_A$ has the same objects as~$\bfV_A$ and there are no morphisms from~$\rmC_A[r]$ to~$\rmC_A[s]$ unless there exists a~$k$ such that~\hbox{$\rmC_A[k]\oplus\rmC_A[r]\cong\rmC_A[s]$} in~$\bfV_A$, i.e.~$r\leqslant s$ in our case, with~\hbox{$k=s-r$}. If this is the case, morphisms are equivalence classes~$[f]$ of morphisms~$f$ in~$\rmV_{n,s}=\Aut_{\bfV_A}(\rmC_A[s])$ with~$f\sim f'$ if there exists an element~$g$ in~$\rmV_{n,k}$ such 
that
\[
f'=f\circ (g\oplus\rmC_A[r])\colon\ \rmC_A[s]=\rmC_A[k]\oplus\rmC_A[r]\ \rar\ \rmC_A[s].
\] 
Note that the unit~$\emptyset=\rmC_A[0]$ is now an initial object in the category~$\bfQ_A$. We will write~\hbox{$\iota_{r}\colon\emptyset\to\rmC_A[r]$} for the unique morphism, which we can represent as the equivalence class~$[\rmC_A[r]]$ of the identity in~$\rmV_{n,r}$.

\begin{proposition}
The category~$\bfQ_A$ is a permutative and homogeneous category with maximal subgroupoid~$\bfV_A$.   
\end{proposition}

\begin{proof} This is a direct application of three results in~\cite{Randal-Williams+Wahl}: Because~$(\bfV_A,\oplus,\s,\emptyset)$ is a symmetric monoidal groupoid,~\cite[Prop.~1.7]{Randal-Williams+Wahl} gives that~$\bfQ_A$ (denoted~$U\bfV_A$ in that paper) inherits a symmetric monoidal structure, with its unit~$\emptyset$ initial. And given that~$\bfV_A$ was actually permutative, so is~$\bfQ_A$. 
We have that~$\Aut(\emptyset)=\{\id\}$ and that there are no zero divisors in~$\bfQ_A$: If there is an isomorphism~$\rmC_A[r]\oplus\rmC_A[s]\cong\emptyset$ in~$\bfV_A$, then we must have that~$r=s=0$. Then~\cite[Prop.~1.6]{Randal-Williams+Wahl} gives that~$\bfV_A$ is the maximal subgroupoid of~$\bfQ_A$. Also, the groupoid~$\bfV_A$ satisfies cancellation~(by construction): If there exists an isomorphism~$\rmC_A[r]\oplus\rmC_A[s]\cong\rmC_A[r]\oplus\rmC_A[s']$ in~$\bfV_A$, then~\hbox{$\rmC_A[s]\cong\rmC_A[s']$} in~$\bfV_A$; they are in fact equal, because we have~\hbox{$\rmC_A[r+s]\cong\rmC_A[r+s']$} in the groupoid~$\bfV_A$ if and only if~$r+s=r+s'$. Finally, the groupoid~$\bfV_A$ satisfies that the map~\hbox{$\Aut_{\bfV_A}(\rmC_A[r])\to\Aut_{\bfV_A}(\rmC_A[r+s])$} adding the identity on~$\rmC_A[s]$ is injective by Proposition~\ref{injectiveAut}. It then follows from~\cite[Thm.~1.9]{Randal-Williams+Wahl} that~$\bfQ_A$ is homogenous, which completes the proof. 
\end{proof}

%\begin{remark}
%In fact, we have forced cancellation in the groupoid~$\bfV_a$ by forgetting the isomorphisms displayed in Example~\ref{isoex}. Including all isomorphisms of Cantor algebras in the groupoid, we would have had that~$\rmC_1\oplus\rmC_{n-1}=\rmC_a\cong\rmC_1=\rmC_1\oplus\rmC_0$, giving a counter-example to cancellation, because~$\rmC_{n-1}$ is not isomorphic to~$\rmC_0$.
%\end{remark}

\begin{remark}\label{rem:comp}
Explicitly, the monoidal structure of~$\bfQ_A$ is defined as follows. On objects, it is as in~$\bfV_A$ induced by the sum of natural numbers. Given two morphisms~$[f]\in \bfQ_A(\rmC_A[r],\rmC_A[s])$ and~\hbox{$[g]\in \bfQ_A(\rmC_A[r'],\rmC_A[s'])$}, the equivalence class of the composition
\[
\xymatrix{\rmC_A[s+s']\ar@{=}[r]&\rmC_A[s-r]\oplus\rmC_A[s-r']\oplus\rmC_A[r]\oplus\rmC_A[r']\ar[d]^{\rmC_A[s-r]\oplus\sigma_{s-r',r}\oplus\rmC_A[r']} & &\\ 
&\rmC_A[s-r]\oplus\rmC_A[r]\oplus\rmC_A[s-r']\oplus\rmC_A[r']\ar@{=}[r]&\rmC_A[s]\oplus\rmC_A[s']\ar[d]^-{f\oplus g}\\
 &&\rmC_A[s]\oplus\rmC_A[s']\ar@{=}[r]&\rmC_A[s+s'].}
\]
defines~$[f]\oplus [g]\in \bfQ_A(\rmC_A[r+r'],\rmC_A[s+s'])$, where~$\sigma_{s-r',r}$ is the symmetry of~$\bfV_A$. (See the proof of Proposition~1.6 in \cite{Randal-Williams+Wahl}.) 
\end{remark} 

%%%

\subsection{Representing morphisms in the category~$\bfQ_A$}

Recall from Proposition~\ref{triples} that the isomorphisms of~$\Can_A$, and hence also the elements~$f$ of~$\rmV_{n,s}$, admit a unique minimal
presentation~$(E,F,\la)$ where~$E,F$ are expansions of the standard generating set~$[s]$
of~$\rmC_A[s]$, and~$\la\colon E\to F$ is a bijection. 
We will now define an analogous minimal representation for the morphisms in~$\bfQ_A(\rmC_A[r],\rmC_A[s])$ for~\hbox{$r<s$}.
This description makes use of the posets~$\calE[r], \, \calI[r]$ and $\calI_0[r]$ of  expansions, independent sets, and non-generating independent sets of~$[r]$ of Section~\ref{sec:exp}.

An element~$[f]\in \bfQ_A(\rmC_A[r],\rmC_A[s])$ is represented by some~$f\in \rmV_{n,s}$, which itself admits a representative~$(E,F,\la)$, with~$E,F\in \calE[s]$ expansions of~$[s]$ and~$\la\colon E\to F$ a bijection. Using the decomposition 
$$[s]=[s-r]\oplus [r],$$ 
with the associated embedding~$i_R\colon\rmC_A^+[r]\inc \rmC_A^+[s]$, 
we can consider the intersection~$E_0:=E\cap i_R\rmC_A^+[r]$, which is an expansion of~$[r]$ by Lemma~\ref{lem:exp}. Then~\hbox{$P:=\la(E_0)\subset \la(E)=F$} is an independent set of~$\rmC_A[s]$. This way we can associate to any morphism~$[f]\in \bfQ_A$ a triple~$(E_0,P,\la)$ with~$E_0\in \calE[r]$ and~$P\in \calI[s]$ and~$\la\colon E_0\to P$ a bijection. 
With this in mind, we extend the definition of~$\Rep(f)$ from Section~\ref{sec:isos} to all morphisms of~$\bfQ_A$: 
\begin{Def}
For~$[f]\in \bfQ_A(\rmC_A[r],\rmC_A[s])$, let 
$$\Rep[f]=\{(E,P,\la)\ |\ E\in \calE[r], P\in \calI[s], \la=[f]|_E\colon E\stackrel{\cong}\rar P\},$$
which we consider as a poset by setting~$(E,P,\la)\le (E',P',\la')$ if and only if~$E\le E'$ in~$\calE[r]$.  
\end{Def}

Note that~$[f]|_E$ is a well-defined map as~$E\subset \rmC_A[r]\equiv i_R\rmC_A[r]$ and any two representatives of~$[f]$ agree as maps on that subalgebra. 

If~$r=s$, we have that~$\Rep[f]=\Rep(f)$ as defined earlier. In particular, if~$(E,P,\la)\in \Rep[f]$ in this case, then~$P\in \calE[s]$. On the other hand, if~$r<s$ and~$(E,P,\la)\in \Rep[f]$, we must have that~$P\in \calI_0[s]$ as~$P$ is non-generating in that case. 
Proposition~\ref{triples} extends to the following result: 

\begin{prop}\label{triples2}
For any morphism~$[f]\in \bfQ_A(\rmC_A[r],\rmC_A[s])$, the set~$\Rep[f]$ is non-empty and has a least element. 
\end{prop}

\begin{proof}
We have already seen above that the set is non-empty. We will in this proof construct a least element. Pick a representative~$f\in \rmV_{n,s}$ of~$[f]$, and let~$(E,F,\la)$ be the least element of~$\Rep(f)$ given by Proposition~\ref{triples}. As above, we consider the set~$[s]$ as the sum~$[s-r]\oplus [r]$, which gives embeddings~$i_L$ and~$i_R$  of~$\rmC_A[s-r]$ and~$\rmC_A[r]$ inside~$\rmC_A[s]$. 
Consider the triple~$(E_0,P,\la|_{E_0})$ with~$E_0=E\cap i_R\rmC_A^+[r]$ and~$P=\la(E_0)$. We have that~$E_0\in \calE[r]$ by Lemma~\ref{lem:exp} and~$P\in \calI[s]$ because~$P\subset \la(E)=F$. Given that~$f|_E=\la$, we get that 
$(E_0,P,\la|_{E_0})\in \Rep[f]$. We claim that~$(E_0,P,\la|_{E_0})$ is the least element of~$\Rep[f]$. So suppose~$(E_0',P',\la')\in \Rep[f]$ is another element. Set~$$E'=(E\cap i_L\rmC_A^+[s-r])\cup E_0'.$$ By Lemma~\ref{lem:exp}, we have~$E'\in \calE[s]$, and hence this is a basis for~$\rmC_A[s]$. 
%Suppose that~$E'_0\le E_0$. Then we also have~$E'\le E$. Note also that 
Now~$f(E')=f(E\cap \rmC_A^+[s-r])\cup P'$ is a basis, as~$f$ is an isomorphism, and it is included in~$\rmC_A^+[s]$ as both~$f(E\cap \rmC_A^+[s-r])$  and~$P'$ are. Hence it is an expansion of~$[s]$ by Lemma~\ref {lem:expord} and~$(E',f(E'),f|_{E'})$ is an element of~$\Rep(f)$. By the minimality of~$(E,F,\la)$, we must have that~$E\le E'$, which gives~$E_0\le E_0'$ as requested. 
\end{proof}

Define 
$$\rmP_A(r,s)=\left\{ \begin{array}{ll} \{(E,P,\la)\ |\ E\in\calE[r], P\in\calI_0[s], \la\colon E\stackrel{\cong}\rar P\}  & \textrm{if}\ \ r<s\\
\{(E,P,\la)\ |\ E\in\calE[r], P\in\calE[s], \la\colon E\stackrel{\cong}\rar P\}  & \textrm{if}\ \ r=s. \end{array}\right.$$
%\end{itemize} 
%\begin{itemize}
%\item~$E\subset \rmC_{n,r}^+$ is an expansion of the canonical basis~$[r]$, 
%\item~$P\subset \rmC_{n,s}^+$ is an independent set which is not a basis, %which is a basis if and only if~$r=s$,  
%\item~$\la\colon E\to P$ is bijection
%\end{itemize} 
We partially order~$\rmP_A(r,s)$ by setting~\hbox{$(E,P,\la)\leqslant (E',P',\la')$} if and only if~$E\le E'$ and~$\la'$ is the restriction to~$E'$ of the map~$\rmC_A(E)\to \rmC_A(P)\to \rmC_A[s]$ induced by~$\la$. As~$P'=\la'(E')$, it follows that~$P'$ is an expansion of~$P$ in such a situation.
%with~$A=[n]$ as above. 

\begin{prop}\label{prop:morphisms}
For all~$0\leqslant r\leqslant s$, there is a poset isomorphism~$$\rmP_A(r,s)\cong \bigsqcup_{[f]\in \bfQ_A(\rmC_A[r],\rmC_A[s])}\Rep[f].$$
In particular, taking least representative in~$\Rep[-]$ defines an isomorphism between~$\bfQ_A(\rmC_A[r],\rmC_A[s])$ and the set of minimal elements in~$\rmP_A(r,s)$. 
%A morphism~$[f]\colon\rmC_r\to\rmC_s$ in~$\bfQ_a$ with~$r\leqslant s$ is determined by a triple~$(F=F_1\sqcup\dots\sqcup F_r,G=G_1\sqcup\dots\sqcup G_s;\lambda)$ for~$F$ an~$(n,r)$--forest,~$G$ an~$(n,s)$--forest, and~$\lambda\colon L(F)\to L(G)$ an injective map which is surjective if and only if~$r=s$. Moreover, every such triple represents a morphism in~$\bfQ_a$, and for each morphism of~$\bfQ_a$, there is a unique minimal representative under expansion. 
%
%Alternatively, we can represent a morphism from~$\rmC_r$ to~$\rmC_s$ as a pair~$(B,\la)$ where~$B$ is a basis for~$\rmC_{n,r}$ and~$\la\colon B\to\rmC_{a,s}$ is an injective map with image an {\color{red}independent} set, which is a basis if and only if~$r=s$.  
\end{prop}

\begin{proof}
Forgetting which~$[f]$ a triple represents defines  a  map 
$$\al\colon\bigsqcup_{[f]\in \bfQ_A(\rmC_A[r],\rmC_A[s])}\Rep[f] \ \rar\ \rmP_A(r,s)$$
which is a map of posets as the order relation is defined in the same way in both cases. We want to show that this map is a poset bijection, that is a set bijection such that any two elements that are related in the target, are also related in the source, or equivalently, that any two elements that are related  in~$\rmP_A(r,s)$ represent the same morphism~$[f]$. 

When~$r=s$,  any element~$(E,P,\la)\in \rmP_A(s,s)$ determines a unique automorphism~$f$ of~$\rmC_A[s]$, as~$E$ and~$P$ are bases in this case. This proves that~$\al$ is bijective
in that case. 
It is also immediate that, if~\hbox{$(E,P,\la)\le (E',P',\la')$} in~$\rmP_A(s,s)$, then the triples represent the same morphism~$f$, proving the result in the case~$r=s$. 

We now assume that~$r<s$. To show that~$\al$ is injective, 
suppose that~$[f]$ and~$[f']$ have the same image~\hbox{$(E,P,\la)\in\rmP_A(r,s)$}. 
%Pick representatives~$f,f'$ for the two morphisms, and consider their minimal representatives~$\al(f)=(\bar E,F,\bar \la)$ and~$\al(f')=(\bar E',F',\bar \la')$. 
Consider the isomorphism 
\[
g:=(f')^{-1}\circ f\colon\rmC_A[s]=\rmC_A[r-s]\oplus \rmC_A[r]\rar \rmC_A[r-s]\oplus \rmC_A[r]=\rmC_A[s].
\]
%Pick a representative~$(F,G,\mu)\in \Rep(g)$ with the property that~$E\subset F$: indeed, let~$\hat E$ be a common expansion of~$F\cap~$ 
We have that~$g$ restricts to the identity on~$E$ as~$f|_E=\lambda=f'|_E$. By Proposition~\ref{prop:split},  it follows that~\hbox{$g=g_0\oplus\rmC_A[r]$} for some~$g_0\in\rmV_{n,s-r}$. %represented by~$(F\minus E,G\minus E,\mu|_{F\minus E})$. 
Hence~$f=f'\circ (g_0\oplus \rmC_A[r])$, which proves that~$[f]=[f']$.  

To show surjectivity 
in the case~$r<s$, suppose that~$(E,P,\la)\in \rmP_A(r,s)$. 
%We are left to show that~$\al$ has image the set of minimal elements of~$\rmP_A(r,s)$. 
%Let~$(E,P,\la)\in \rmP_A(r,s)$ be a minimal element. 
Let~$F$ be an expansion  of~$[s]$ containing~$P$. We have that~\hbox{$|E|=r+a(n-1)=|P|$} for some~$a\geqslant0$ and~$|F|=s+b(n-1)$ for some~\hbox{$b\geqslant0$}, with~\hbox{$|F|-|P|=s-r+(b-a)(n-1)>0$} under our assumption that~$(E,P,\la)\in \rmP_A(r,s)$ with~\hbox{$r<s$}. If~\hbox{$b-a<0$}, 
replace~$F$ by an expansion of~$F$ still containing~$P$ by expanding an element of~\hbox{$F\setminus P$}~(which is non-empty by the same assumption) at least~$(b-a)$ times. 
After doing this, we can assume moreover  that~\hbox{$b-a\geqslant0$}. Now let~$G$ be an expansion of~$[s-r]$ of cardinality~$s-r+(b-a)(n-1)$ and pick a bijection~$\mu\colon G\to F\setminus P$. 
Then~$G\cup E$ is an expansion of~$[s]$ and~$(G\cup E,\mu(G)\cup P,\mu\cup\la)$ represents an element~$f\in\rmV_{n,s}$. Let~$[f]\in \bfQ_A(\rmC_A[r],\rmC_A[s])$ be its equivalence class. 
By construction, we have~\hbox{$(E,P,\la)\in \Rep[f]$} which gives surjectivity 
also in this case.

Finally, still assuming~$r<s$, we need to check that~$(E,P,\la)\le (E',P',\la')$ in~$\rmP_A(r,s)$ can only happen if both triples represent the morphism~$[f]\in \bfQ_A(\rmC_A[r],\rmC_A[s])$, with~$\la$ and~$\la'$ necessarily the restriction of a representative~$f$ of~$[f]$ to~$E$ and~$E'$ respectively. By the surjectivity of~$\al$, we know that~$(E,P,\la)$ represents some~$[f]\in \bfQ_A(\rmC_A[r],\rmC_A[s])$. Let~$(G\cup E,\mu(G)\cup P,\mu\cup \la)\in \Rep(f)$ be the representative of such an~$f$ constructed above. Now~$(G\cup E,\mu(G)\cup P,\mu\cup \la)\le (G\cup E',\mu(G)\cup P',\mu\cup \la')$, showing that both triples necessarily represent the same morphism~$f$.  Hence~$(E,P,\la)$ and~$(E',P',\la')$ are both in~$\Rep[f]$ and were already comparable there. 
\end{proof}

%%%

\subsection{Spaces constructed from~$\bfQ_A$}

In the general context of the paper~\cite{Randal-Williams+Wahl}, given a pair~$(B,X)$ of objects in a homogeneous category, a sequence of semi-simplicial sets~$\rmW_r(B,X)$ is defined, and the main theorem in that paper says that homological stability holds for the automorphism groups of the objects~\hbox{$B\oplus X^{\oplus r}$} 
as long as the associated semi-simplicial sets are highly connected. In good cases, the connectivity of the semi-simplicial sets~$\rmW_r(B,X)$ can be computed from the connectivity of closely related simplicial complexes~$\rmS_r(B,X)$. 

We are here interested in the pair of objects~$(B,X)=(\emptyset,\rmC_A[1])$ in the homogeneous category~$\bfQ_A$. Indeed, the automorphism group of~$\emptyset\oplus\rmC_A[1]^{\oplus r}=\rmC_A[r]$ in the category~$\bfQ_A$ is the Higman--Thompson group~$\rmV_{n,r}$. We will therefore begin by describing the semi-simplicial sets~\hbox{$\rmW_r=\rmW_r(\emptyset,\rmC_A[1])$} and the simplicial complexes~$\rmS_r=\rmS_r(\emptyset,\rmC_A[1])$ from Definitions~2.1 and~2.8 in~\cite{Randal-Williams+Wahl}, and show that we are in a situation where we can use the connectivity of the latter to compute the connectivity of the former. In the following Section~\ref{sec:connectivity} we will estimate that connectivity. 

\begin{definition}
Given an integer~$r\geqslant1$, let~$\rmW_r=\rmW_r(\emptyset,\rmC_A[1])$ be the semi-simplicial set where the set of~$p$--simplices is the set~$\bfQ_A(\rmC_A[p+1],\rmC_A[r])$ of morphisms $\rmC_A[p+1]\to\rmC_A[r]$ in~$\bfQ_A$ and the~$i$--th boundary map~\hbox{$\bfQ_A(\rmC_A[p+1],\rmC_A[r])\to \bfQ_A(\rmC_A[p],\rmC_A[r])$} is defined by precomposing 
with~\hbox{$\rmC_A[i]\oplus\iota_{1}\oplus\rmC_A[p-i]$}.
\end{definition}

%Note that~$\rmW_r$ has dimension~$r-1$, and a top dimensional simplex in~$\rmW_r$ is an automorphism of~$\rmC_r$, which is an element of~$\rmV_{n,r}$.  
%For~$p<r-1$, Proposition~\ref{prop:morphisms} shows that~$p$--simplices of~$\rmW_r$ have a unique minimal presentation by a triple~$(E,P,\la)$ where~$E$ is an expansion of~$[p+1]$, and~$P$ is an independent set of~$\rmC_r^+$ that is not a basis, and~$\la\colon E\to P$ is a bijection. 

To the semi-simplicial set~$\rmW_r$ we associate the following simplicial complex, of the same dimension~\hbox{$r-1$}. 

\begin{definition} 
Given~$r\geqslant1$, let~$\rmS_r=\rmS_r(\emptyset,\rmC_A[1])$ be the simplicial complex with the same vertices as~$\rmW_r$, namely the set of maps~$\bfQ_A(\rmC_A[1],\rmC_A[r])$, and where~$p+1$ distinct vertices~$[f_0],\dots,[f_p]$ form a~$p$--simplex if there exists a~$p$--simplex of~$\rmW_r$ having them as vertices.
%such that~$f\circ i_j=f_j$ for each~$j=0,\dots,p$, where~$i_j=\iota_{j}\oplus\id_{\rmC_1}\oplus\iota_{{p-j}}\colon\rmC_1\to\rmC_{p+1}$ for~$\iota_{k}\colon 0=\rmC_0\to\rmC_k$ the unique morphism coming from the fact that 0 is initial in~$\bfQ_A$. 
\end{definition}

Given a simplicial complex~$S$, one can build a semi-simplicial set~$S^{\ord}$ that has a~$p$--simplex for each ordering of the vertices of each~$p$--simplex of~$S$.  

\begin{proposition}\label{prop:S=>W}
 There is an isomorphism of  semi-simplicial sets~$\rmW_r\cong\rmS_r^{\ord}$. Moreover, if~$\rmS_r$ is~$(r-3)$--connected, then so is~$\rmW_r$.  
\end{proposition}

\begin{proof}
This follows from~\cite[Prop.~2.9, Thm.~2.10]{Randal-Williams+Wahl}, given that~$\bfQ_A$ is symmetric monoidal, once we have checked that it is {\em locally standard} at~$(\emptyset,\rmC_A[1])$ in the sense of~\cite[Def.~2.5]{Randal-Williams+Wahl}. This means two things: Firstly, the morphisms~$\rmC_A[1]\oplus\iota_{1}$ and~$\iota_{1}\oplus\rmC_A[1]$ are distinct in~$\bfQ_A(\rmC_A[1],\rmC_A[2])$, and secondly, for all~$r\geqslant1$, the map~$\bfQ_A(\rmC_A[1],\rmC_A[r-1])\to \bfQ_A(\rmC_A[1],\rmC_A[r])$ that takes a morphism~$[f]$ to~$[f]\oplus\iota_{1}$ is injective.

For the first statement, 
we need to describe the two morphisms~$\rmC_A[1]\oplus\iota_{1}$ and~$\iota_{1}\oplus\rmC_A[1]$ precisely. Here~\hbox{$\rmC_A[1]\in \bfQ_A(\rmC_A[1],\rmC_A[1])=\bfV_A(\rmC_A[1],\rmC_A[1])$} is the identity, while~$\iota_1=[\rmC_A[1]]\in\bfQ_A(\emptyset,\rmC_A[1])$ is also represented by the identity on~$\rmC_A[1]$, but this time only up to an automorphism of~$\rmC_A[1]$. 
Using the explicit definition of the composition given in Remark~\ref{rem:comp}, we compute that~$\rmC_A[1]\oplus\iota_{1}=[\s_{1,1}]$ for~\hbox{$\s_{1,1}\colon\rmC_A[2]=\rmC_A[1]\oplus \rmC_A[1]\to \rmC_A[1]\oplus\rmC_A[1]=\rmC_A[2]$} the symmetry, while~$\iota_{1}\oplus\rmC_A[1]=\rmC_A[2]$ 
is represented by the identity (because the symmetry needed for that composition is~$\s_{0,0}$, which is the identity). 
To check that these morphisms are distinct in~$\bfQ_A(\rmC_A[1],\rmC_A[2])$, we use the minimal presentations (Propositions~\ref{triples2} and~\ref{prop:morphisms}). 
We have that~$[\s_{1,1}]$ has minimal presentation~$(\{1\},\{2\},\la)$ 
while the second has presentation~$(\{1\},\{1\},\mu)$, for~$\la,\mu$ the unique maps, showing that they are indeed distinct.

For the second statement, we have
\[
[f]\oplus\iota_{1}=[(f\oplus\rmC_A[1])\circ (\rmC_A[r-2]\oplus\s_{1,1})].
\]
If~$[f]$ is minimally presented by~$(E,P,\la)$, then one can check that~$[f]\oplus\iota_{1}$ is minimally 
presented by~$(E,i_L(P),i_L\circ \la)$ for~$i_L\colon\rmC^+_A[r-1]\to \rmC^+_A[r]$ the left embedding. As~$i_L$ is injective, the result follows.

By~\cite[Prop.~2.9]{Randal-Williams+Wahl} we now know that~$\rmW_r$ ``satisfies condition (A)'', which means that it is isomorphic to~$\rmS_r^{\ord}$, and~\cite[Thm.~2.10]{Randal-Williams+Wahl} of the same paper gives the second part of the statement.  
\end{proof}

\subsection{Variations}

Using the minimal
representatives of morphisms of~$\bfQ_A$ given by Proposition~\ref{prop:morphisms}, we can represent the vertices of~$\rmW_r$ and~$\rmS_r$ as triples~$(E,P,\la)$ where~$E\in \calE[1]$ is an expansion of~$[1]$ and~$P\in \calI[r]$ is  an independent set~(non-generating if~$r>1$ and an expansion if~$r=1$). 
Note that~$E$, and hence also~$P$, necessarily has cardinality~\hbox{$1+a(n-1)$} for some~$a\ge 0$. 
To study the connectivity of the simplicial complex~$\rmS_r$, we will use variants~$\rmU_r$,~$\rmU_r^{\infty}$ and~$\rmT_r^\infty$ of this complex, where only the independent sets~$P$ are remembered. 

Recall that~$\calI_0[r]=\calI[r]\minus \calE[r]$ denotes the set of non-generating independent sets. 

\begin{definition}\label{def:US}
Let~$\rmU_1=\calE[1]$ to be the simplicial complex of dimension 0 consisting of all the expansions of the set~$[1]$. 
For~$r\geqslant2$, let~$\rmU_r$ be the simplicial complex of dimension~\hbox{$r-1$} 
with vertices the non-generating independent subsets~$P\in \calI_0[r]$ of cardinality congruent to~$1$ modulo~$n-1$. A set of~\hbox{$p+1$} vertices~$P_0,\dots,P_q$ forms a~$q$--simplex of~$\rmU_r$ if the sets~$P_i$ are pairwise disjoint and 
\begin{itemize}[noitemsep,topsep=0pt]
\item~$q<r-1$ and~$P_0\sqcup\dots\sqcup P_q\in \calI_0[r]$, or
\item~$q=r-1$ and~$P_0\sqcup\dots\sqcup P_q\in\calE[r]$.
\end{itemize} 
\end{definition}

%Using Proposition~\ref{prop:morphisms}, we see that 
%there is a forgetful map~$\rmS_r\to\rmU_r$ that takes a vertex~$(E,P,\la)$ of~$\rmS_r$ to the independent set~$P$. We in fact have the following: 

Recall from~\,\upshape{\cite[Def.~3.2]{Hatcher+Wahl}} that a {\em complete join complex} over a simplicial complex~$X$ is a simplicial complex~$Y$ with a projection~$\pi\colon  Y\to X$, which is surjective, injective on each simplex, and with the property that, for each simplex~$\s=\langle x_0,\dots,x_q\rangle$ of~$X$, the subcomplex~$\pi^{-1}(\s)$ of~$Y$ of simplices projecting down to~$X$ is the join complex~$\pi^{-1}(x_0)*\dots*\pi^{-1}(x_q)$. 

\begin{prop}\label{join} 
The simplicial complex~$\rmS_r$ is a  complete join complex over~$\rmU_r$. 
\end{prop}

To prove the proposition, we will use the following lemma, which describes the vertices of any simplex of~$\rmW_r$. We introduce first a notation: for any~$1\le j\le r$, let 
$$i_j\colon[1]\to [r]$$
denote the map defined by~$i_j(1)=j$. We will also write 
$$i_j\colon\rmC_A[1]\rar \rmC_A[r]$$
for the induced map of Cantor algebras. 

\begin{lem}
Let~$[f]\in \bfQ_A(\rmC_A[q+1],\rmC_A[r])$ be a~$q$--simplex of~$\rmW_r$, and let~$(E,P,\la)\in \Rep[f]$ be a triple representing~$[f]$. Then the vertices of~$[f]$ in~$\rmW_r$ are represented by the triples 
$(E_0,P_0,\la_0),\dots,(E_q,P_q,\la_q)$ where $$E_j=E\cap i_j\rmC_A[1], \ \ \ P_j=\la(E_j) \ \ \ \textrm{and} \ \ \la_j=\la|_{E_j}.$$ Moreover,~$(E,P,\la)$ is the  minimal representative of~$[f]$ if and only if 
each~$(E_j,P_j,\la_j)$ is likewise minimal.
\end{lem}

\begin{proof}
The~$j$--th vertex of the simplex~$[f]$ is given by~\hbox{$[f]\circ (\iota_{j}\oplus \rmC_A[1]\oplus \iota_{q-j})\in \bfQ_A(\rmC_A[1],\rmC_A[r])$}. This composition is represented by the following automorphism of~$\rmC_A[r]$: 
\begin{equation}\label{equ:vert}
\xymatrix@R=11pt{
\rmC_A[r-q-1]\oplus \big(\rmC_A[j]\oplus \rmC_A[q-j]\oplus\rmC_A[1]\big)\ar[d]^{\id\oplus\s_{q-j,1}}\\
\rmC_A[r-q-1]\oplus \big(\rmC_A[j]\oplus \rmC_A[1]\oplus \rmC_A[q-j]\big)\ar[d]^\id\\
\rmC_A[r-q-1]\oplus \rmC_A[q+1]\ar[d]^f\\
\rmC_A[r]
}
\end{equation}
Now suppose~$[f]$ is represented by~$(E,P,\la)$ with~$f\in \Aut(\rmC_A[r])$ a representative of~$[f]$. This means that we have~\hbox{$E\in \calE[q+1]$},~$P\in \calI[r]$ and~\hbox{$\la=f|_{E}$}. Also, there exists~$F\in\calE[r-q-1]$ such that~$f$ is represented by 
$$(F\oplus E,f(F)\oplus P,f|_F\oplus \la)=(F\cup i_R(E),f(F)\cup i_R(P),f|_F\cup\la).$$ 
On~$\rmC_A^+[r]=\rmC_A^+[1]\sqcup\dots\sqcup\rmC_A^+[1]$, we see that the composition (\ref{equ:vert}) takes the last component~$\rmC_A^+[1]$ to the~\hbox{$(r-j-1)$}--st position and then applies~$f$.  Thus it is represented by   the triple 
\[
(F\oplus E_0\oplus\dots\oplus E_{j-1}\oplus E_{j+1}\oplus\dots\oplus E_q\oplus E_j,f(F)\oplus P,(f|_F\oplus\la)\circ (\id\oplus \s_{q-j,1}))
\]
where~\hbox{$E_j=E\cap i_j\rmC_A[1]\subset\rmC_A[r]$}. 
This proves the first part of the result, as the representative of the corresponding morphism in the category~$\bfQ_A(\rmC_A[1],\rmC_A[r])$ is obtained by restricting the source to~$i_R\rmC_A[1]$. 

Finally, the statement about minimality follows as in the proof of Proposition~\ref{triples2} by showing that if the~$j$--th vertex has a smaller representative~$(E_j',P_j',\la_j')$, we can use it to construct a smaller representative for~$[f]$ simply by replacing~$E\cap i_j\rmC_A[1]$ with~$i_j(E_j')$, and vice versa. 
\end{proof}

\begin{proof}[Proof of Proposition~\ref{join}]
Using Proposition~\ref{prop:morphisms}, we identify the vertices of the simplicial complex~$\rmS_r$ with the minimal
elements of~$\rmP_A(1,r)$. Using  the lemma, we see that, under this identification,  a collection of vertices~$(E_0,P_0,\la_0),\dots,(E_q,P_q,\la_q)$ forms a~$q$--simplex if and only if there exists a minimal element~\hbox{$(E,P,\la)\in \rmP_A(q,r)$} 
%having~$(E_0,P_0,\la_0),\dots (E_q,P_q,\la_q)$ as vertices in~$\rmW_r$. And this is the case precisely if 
such that 
$E=E_0\sqcup\dots\sqcup E_q$ for~$E_j=E\cap i_j\rmC_A^+[1]$,  and for each~$j$,~$\la_j=\la|_{E_j}$ and~$P_j=\la(E_j)$. 
It follows that~$P=P_0\cup\dots\cup P_q$ with each~$P_j$ pairwise disjoint in~$\rmC_A^+[r]$, and~$\la=\la_0\cup \dots\cup \la_q$. 
Note that by Lemma~\ref{lem:exp}, each~$E_j$ being an expansion of~$[1]$ gives that~$E$ is an expansion of~$[q+1]$.  
%Also, again by the lemma, minimality of the~$E_i$'s implies that~$E$ is minimal. 
%Given that we can also take~$E$ to be the union of the~$E_i$'s and~$\la$ the associated map, we see that the only condition for such a collection of vertices to form a simplex is in fact that the subsets~$P_0,\dots,P_q$ of~$\rmC_A[r]^+$ are pairwise disjoint and that their union is a 
It follows that~$P_1\cup\dots\cup P_q$ is an independent set that is  non-generating if~$q+1<r$, and an expansion of~$[r]$ if~$q+1=r$.  

We define the projection~$\pi\colon\rmS_r\to\rmU_r$ on vertices by setting~$\pi(E,P,\la)=P$. From the above discussion, we see that~$\pi$ is a map of simplicial complexes. 
The map is surjective as for any~$q$--simplex~$\langle P_0,\dots,P_q\rangle$ of~$\rmU_r$, we can choose expansions~$E_j$ of~$[1]$ with~$|E_j|=|P_j|$, by the condition on the cardinality of the~$P_j$'s, and bijections~$\la_j\colon E_j\to P_j$. Then the collection~$(E_0,P_0,\la_0),\dots,(E_q,P_q,\la_q)$ will define a simplex of~$\rmS_r$ projecting down to the given simplex. The map is injective on individual simplices as the~$P_j$'s of distinct vertices in a simplex of~$\rmS_r$ are by definition distinct. Finally we check the join condition. For a vertex~$P$ of~$\rmU_r$, we have that~$\pi^{-1}(P)$ is the set of~$(E,P,\la)$ where~$E\in \calE[1]$ has the same cardinality as~$P$ and~$\la\colon P\to E$ is a bijection. As any choice of~$(E_j,\la_j)$ for each vertex~$P_j$ of a simplex~$\s=\langle P_0,\dots,P_q\rangle$ of~$\rmU_r$ defines a lift~$\langle (E_0,P_0,\la_0),\dots,(E_q,P_q,\la_q)\rangle$ to~$\rmS_r$, we see that~$\pi^{-1}(\s)$ is the join~$\pi^{-1}(P_0)*\dots*\pi^{-1}(P_q)$, finishing the proof of the result. 
\end{proof}

In contrast to the simplicial complexes~$\rmS_r$ and~$\rmU_r$, the variants~$\rmU_r^{\infty}$ and~$\rmT_r^\infty$ that we will now introduce are both infinite-dimensional (as the notation suggests).

\begin{definition}\label{def:infty}
For~$r\geqslant1$, let~$\rmU_r^{\infty}$ be the simplicial complex whose vertices are the non-generating independent sets~$P\in \calI_0[r]$ of cardinality congruent to~$1$ modulo~$n-1$. Distinct vertices~$P_0,\dots,P_p$ form a~$p$-simplex in~$\rmU_r^{\infty}$ if the subsets are pairwise disjoint and their union~$P_0\sqcup\dots\sqcup P_p$ is also an element of~$\calI_0[r]$, that is the union is  a non-generating independent set. 
%even when~$p=r-1$, and also consequently allowing simplices of any dimension
\end{definition} 

Note that~$\rmU_r^\infty$ and~$\rmU_r$ have the same set of vertices when~$r\geqslant2$ but not when~\hbox{$r=1$}. In fact we have the following: 

\begin{lemma}\label{lem:skeleton}
The simplicial complexes~$\rmU_r$ and~$\rmU_r^\infty$ share the same~\hbox{$(r-2)$}--skeleton. 
%and so do~$\rmS_r^0$ and~$\rmS_r^{0,\infty}$. 
\end{lemma}

\begin{proof}
This follows immediately from the definitions. 
\end{proof}

\begin{definition}\label{def:infty,0}
For~$r\geqslant1$, let~$\rmT_r^\infty$ denote the full subcomplex of~$\rmU_r^\infty$ on the vertices that have cardinality~$1$.
%are ({\color{red}independent}) subsets 
%%of~$\rmF_{n,r}$ 
%of cardinality~1 (that are not bases).
\end{definition} 

The diagram 
\[
\xymatrix@R=10pt{
&&\rmU_r^\infty&\rmT_r^\infty\ar[l]_{\supseteq}\\
\rmS_r\ar[r]&\rmU_r\\
&\sk^{r-2}\rmU_r\ar[u]^\subseteq\ar@{=}[r]&\sk^{r-2}\rmU_r^\infty\ar[uu]_\supseteq}
\]
indicates the relations between the simplicial complexes that we have introduced so far, with the height reflecting their dimension.
In the next section, we will prove high connectivity of~$\rmS_r$ using this sequence of maps. 
This sequence of maps is chosen using the following ideas:
\begin{itemize}[noitemsep,topsep=0pt]
\item Going from~$\rmS_r$ to~$\rmU_r$: only the sets~$P_i$ are essential for determining the connectivity; 
\item Going from~$\rmU_r$ to~$\sk^{r-2}\rmU_r$: removing the top simplices means that all simplices correspond now to non-generating sets, and hence there is always something, and hence a lot, independent from any given simplex---this will allow coning off; 
\item Going from~$\rmU_r$ to~$\rmU_r^\infty$: it is often easier to show that an infinite complex is contractible; 
\item Going from~$\rmU_r^\infty$ to~$\rmT_r^\infty$: we have more control over the vertices of the latter complex. 
\end{itemize}

%%%

% !TEX root = allinone_main_v2.tex

\section{Homological stability}\label{sec:connectivity}

In this section, we estimate the connectivity of the simplicial complexes defined in the previous section, and use these to deduce homological stability for the canonical diagrams~\eqref{eq:canonical_diagram} of the Higman--Thompson groups. As before, the chosen set~$A=\{a_1,\dots,a_n\}$ is of cardinality~$n\ge 2$.

Recall the simplicial complexes~$\rmT_r^\infty$ and~$\rmU_r^\infty$ from the previous section. 
All the results of this section  will be based on the following two results. 

\begin{prop}\label{contract1}
Suppose~$|A|\ge 2$. Then for each integer~$r\geqslant1$ the simplicial complex~$\rmT_r^\infty$ is contractible. 
\end{prop}

\begin{prop}\label{contract2}
Suppose~$|A|\ge 2$. Then for each integer~$r\geqslant1$ the simplicial complex~$\rmU^\infty_r$ is contractible. 
\end{prop}

Before we give proofs of these two propositions in Section~\ref{sec:proofs}, let us state and prove some of their consequences.

Recall from~\cite[Def.~3.4]{Hatcher+Wahl} that a simplicial complex is called {\em weakly Cohen--Macaulay of dimension~$n$} if it is~$(n-1)$--connected and the link of every~$p$--simplex of it is~$(n-p-2)$--connected. 

\begin{cor}\label{wCM}
For all~$r\geqslant2$ the simplicial complexes~$\rmU_r$ are weakly Cohen--Macaulay of dimension~\hbox{$r-2$}. 
\end{cor}

\begin{proof}
A simplicial complex is~$(r-3)$--connected if and only if its~$(r-2)$--skeleton is. Since the simplicial complexes~$\rmU_r$ and~$\rmU_r^\infty$ share the same~$(r-2)$--skeleton by Lemma~\ref{lem:skeleton}, we see that the simplicial complex~$\rmU_r$ is~\hbox{$(r-3)$}--connected if the simplicial complex~$\rmU_r^\infty$ is, and the latter is even contractible by Proposition~\ref{contract2}. 

Let now~$\sigma$ be a~$p$--simplex of~$\rmU_r$  with vertices~$P_0,\dots,P_p$. We need to check that the link of~$\s$ is~\hbox{$(r-p-4)$}--connected. 
For~$p\geqslant r-2$, there is nothing to check as~$(r-2)-p-2\leqslant-2$ and~$(-2)$--connected is a non-condition. So assume~$p\leqslant r-3$. Then~$P_0\cup\dots\cup P_p\in \calI_0[r]$
is a non-generating independent set. By Lemma~\ref{lem:minexp}, we can find a least~$Q\in \calI_0[r]$ 
(disjoint from the~$P_i$'s) such that~$Q\cup P_0\cup \dots\cup P_p$ is an expansion of~$[r]$.   Note that~$Q$ is non-empty as~$\cup_iP_i$ is non-generating. 

To show the~$(r-p-4)$--connectivity of the link, it is enough to consider its~$(r-p-3)$-skeleton. We claim that there is an isomorphism~$\sk_{r-p-3}(\Link(\s))\cong\sk_{r-p-3}(\rmU_{q}^\infty)$ for~$q=|Q|$.  

To write a map~$\sk_{r-p-3}(\rmU_{q}^\infty)\to\sk_{r-p-3}(\Link(\s))$, we use the injection~\hbox{$i\colon\rmC_A^+[q]\to \rmC_A^+[r]$} induced by a chosen bijection~\hbox{$[q]\cong Q$ and the inclusion~$Q\subset\rmC_A^+[r]$}.
(As we will see in Section~\ref{sec:classifying}, there is actually a canonical such bijection as every~$Q$ has a canonical ordering induced by that of~$[r]$ and that of~$A$, but we can use any bijection here.) 
This map induces a map~$i\colon\calI_0[q]\to\calI_0[r]$ because~$Q$ was an independent set. Define 
\[
\al\colon\sk_{r-p-3}(\rmU_{q}^\infty)\rar \sk_{r-p-3}(\Link(\s))
\]
by setting~$\al\langle Q_0,\dots,Q_q\rangle=\langle i(Q_0),\dots,i(Q_q)\rangle$. This is an injective map of simplicial complexes because~$i$ is injective and preserves the independence condition. To show surjectivity, suppose that~$\langle P_{p+1},\dots,P_s\rangle$ is a simplex of~$\sk_{r-p-3}(\Link(\s))$. Then~$\langle P_0,\dots,P_p, P_{p+1},\dots,P_s\rangle$ is a simplex of~$\rmU_r$ and there must exist~$Q'$ disjoint from the~$P_i$'s such that~$Q'\cup P_0\cup \dots\cup P_s$ is an expansion of~$[r]$. Because~$Q$ was chosen so that~\hbox{$Q\cup P_0\cup \dots\cup P_p$} is the least expansion containing 
$P_0\cup \dots\cup P_p$, we must have that~\hbox{$Q'\cup P_{p+1}\cup \dots\cup P_s$} is an expansion of~$Q$. Now~$s\le p+(r-p-2)=r-2$ because~$\langle P_{p+1},\dots,P_s\rangle$ is in the~$(r-p-3)$--skeleton of the link. Hence the~$P_i$'s together form a non-generating independent set and~$|Q'|>0$.  It follows that~$\langle i^{-1}(P_{p+1}),\dots,i^{-1}(P_s)\rangle$ is a simplex of~$\sk_{r-p-3}(\rmU_{q}^\infty)$. 
This proves that the map~$\al$ is an isomorphism. 

The required connectivity of the link follows from the contractibility of~$\rmU_q^\infty$. 
\end{proof}

Note that~$q=1$ may occur in the proof above. This is the reason why we need to know that~$\rmU_r^\infty$ is contractible also when~$r=1$, even though~$\rmU_1^\infty$ has no relationship to~$\rmU_1$.  

\begin{cor}\label{r-3}
For each~$r\geqslant2$ the spaces~$\rmS_{r}$ and~$\rmW_r$ are~$(r-3)$--connected. 
\end{cor}

\begin{proof}%[Proof of Corollary~\ref{r-3}]
By Proposition~\ref{join}, the simplicial complex~$\rmS_r$ is a complete join complex over~$\rmU_r$, and~$\rmU_r$ is 
weakly Cohen--Macaulay of dimension~$r-2$ by Corollary~\ref{wCM}.  It thus follows from~\cite[Prop.~3.5]{Hatcher+Wahl} that~$\rmS_r$ is also weakly Cohen--Macaulay of that dimension, and so in particular~$(r-3)$--connected. 
Then, by Proposition~\ref{prop:S=>W}, it follows that the semi-simplicial sets~$\rmW_r$ are also~$(r-3)$--connected. 
\end{proof}

\begin{cor}\label{slope2}
The semi-simplicial set~$\rmW_{r+1}$ is~$(\frac{r-2}{2})$--connected for all~$r\geqslant0$. 
\end{cor}

\begin{proof}
There is a morphism~$\rmC_A[1]\to\rmC_A[r+1]$ in the category~$\bfQ_A$ as soon as~$r+1\geqslant1$, or equivalently~$r\geqslant0$. This shows that the semi-simplicial set~$\rmW_{r+1}$ is non-empty for all~$r\geqslant 0$. That gives the cases~\hbox{$r=0,1$}. For~\hbox{$r\geqslant2$}, Corollary~\ref{r-3} gives that~$\rmW_{r+1}$ is~$(r-2)$--connected, which, under the assumption on~$r$, implies that it is at least~$(\frac{r-2}{2})$--connected. 
\end{proof}

%%%

%%%

\subsection{Stability theorem}\label{sec:stabsubsec}

Define the stabilization homomorphism~$s_r\colon\rmV_{n,r}\to\rmV_{n,r+1}$ as the map that  takes an element~$f$ to~$f\oplus\rmC_A[1]$. 
The maps~$s_r$ for all~$r\ge 1$ together give the diagram~\eqref{eq:canonical_diagram} of groups of the introduction.

\begin{theorem}\label{thm:HS}
Suppose~$n\ge 2$. The stabilization homomorphisms induce isomorphisms
\[
s_r\colon \rmH_d(\rmV_{n,r};M)\longrightarrow\rmH_d(\rmV_{n,r+1};M)
\]
in homology in all dimensions~$d\geqslant0$, for all~$r\geqslant 1$, and for all~$\rmH_1(\rmV_{n,\infty})$--modules~$M$.
\end{theorem}

\begin{proof}
First, we can apply the stability result of~\cite{Randal-Williams+Wahl} to the category~$\bfQ_A$ with the chosen pair of objects~$(\rmC_A[1],\rmC_A[1])$: our Corollary~\ref{slope2} shows that the complexes~\hbox{$\rmW_r(\rmC_A[1],\rmC_A[1])\cong\rmW_{r+1}(\emptyset,\rmC_A[1])$} are~$(\frac{r-2}{2})$--connected, which implies, by~\cite[Thm.~3.4]{Randal-Williams+Wahl}, that~$s_r$ is an isomorphism in the range of dimensions~$d\leqslant \frac{r-4}{3}$, that is a range that increases with the `rank'~$r$. Recall now that the group~$\rmV_{n,r}$ is isomorphic to the group~$\rmV_{n,r+(n-1)}$ as soon as~\hbox{$r\geqslant 1$}. We choose an isomorphism between the two groups as follows. Let 
$h\colon \rmC_A[1]\to \rmC_A[n]$ be the isomorphism of Example~\ref{ex:iso} with minimal presentation~$(\{1\}\x A,[n],\la)$ for~$\la\colon A\to [n]$ a bijection (for example the canonical bijection coming from the picked ordering of~$A$). Then~$h_r=h\oplus \rmC_A[r-1]\colon \rmC_A[r]\to \rmC_A[n+r-1]$ is also an isomorphism. 
Let~$\gamma(h_r)$ denote conjugation with~$h_r$. Then we get a commutative diagram 
\[
\xymatrix{\rmV_{n,1+(r-1)}\ar[r]^-{s_r}\ar[d]_{\gamma(h_r)} &\rmV_{n,1+(r-1)+1}\ar[d]^{\gamma(h_{r+1})}\\
\rmV_{n,n+(r-1)}\ar[r]^-{s_{r+1}} &\rmV_{n,n+(r-1)+1}, 
}
\]
as
\[
(h\oplus \rmC_A[r])^{-1}\circ (f\oplus \rmC_A[1])\circ (h\oplus \rmC_A[r])=\big((h\oplus \rmC_A[r-1])^{-1}\circ f\circ (h\oplus \rmC_A[r-1])\big)\oplus \rmC_A[1]
\]
for all~$r\geqslant 1$. 
Given that the vertical maps are isomorphisms and increase the rank of the group (under the assumption~\hbox{$n\ge 2$}), 
while the horizontal maps induce an isomorphism in homology when the rank of the group is large enough,
we get that the horizontal maps must always induce an isomorphism in homology. 
\end{proof}

\begin{remark}
For the purposes of the present paper, we will only need homological stability with respect to trivial, or potentially abelian, coefficients, in the form  stated. 
Applying the more general Theorem~A of \cite{Randal-Williams+Wahl} instead of Theorem~3.4 in that paper, one obtains that stability also holds with finite degree coefficient systems. 
%There are, however, often good reasons to consider more general coefficients, and 
We refer the interested reader to~\cite{Randal-Williams+Wahl} for the definitions. %and results that apply in our situation as well.
\end{remark}

%\Nnote{are there natural examples of twisted coefficients? (question from the referee)}

%\Nnote{new}
%\begin{rem} The definitions we have made all make sense when the set~$A=\{a_1\}$ has cardinality~$1$, but some of the proofs used that the cardinality of~$A$ was at least two. 
%A Cantor algebra of type~$A$ if~$A$ has cardinality~$1$ is a set~$X$ together with a bijection~$f\colon X\to X$. The free Cantor algebra~$\rmC_A[r]$ in this case is~$X=\bbZ\times [r]$ with~$f(n,j)=(n+1,j)$. Its automorphism group is the wreath product~$V_{1,r}=\bbZ\wr\Sigma_r$. Homological stability is known for these groups (see eg.~\cite[Prop.1.6]{Hatcher+Wahl}), though in this case stability only holds 
%with a slope 2 stability bound:~$\rmH_d(\bbZ\wr\Sigma_r;M)\longrightarrow\rmH_d(\bbZ\wr\Sigma_{r+1};M)$ is an isomorphism when~$r\ge 2d+1$. The simplicial complex~$S_r$ is a complete join complex over the simplex~$\Delta^{r-1}$ in that case.
%\end{rem}

%%%

\subsection{Proofs of Propositions~\ref{contract1} and~\ref{contract2}}\label{sec:proofs}

In the rest of this section we give proofs of Propositions~\ref{contract1} and~\ref{contract2}.

Recall that a vertex of~$\rmT_r^\infty$ is a non-generating independent set of cardinality~$1$ in~$\rmC_A^+[r]$. The set
\[
\rmC_A^+[r]=\bigsqcup_{h\geqslant 0}[r]\x A^h
\]
is canonically graded by the {\em height}~$h$ of its elements. If~$r\neq 1$, any element of~$\rmC_A^+[r]$ is non-generating and independent, so the set of vertices of~$\rmT_r^\infty$ is the set of elements of 
$\rmC_A^+[r]$, whereas for~$r=1$ it is all the elements of height at least~$1$, because the bottom element (of height~$0$) is generating. 
%To do so, we first make sure that there is enough space in the target: 

The following lemma makes precise that there are many elements~$y\in \rmC_A^+[r]$ that are independent of any given non-generating independent set~$P$, as long as we go high enough in~$\rmC_A^+[r]$. 

\begin{lemma}\label{lem:height}
For any~$r,k,N\ge 1$, there exists a height~$h_{r,k,N}\ge 0$, such that for any 
non-generating independent set~$P\in \calI_0[r]$ of cardinality~$k$, and any~$h\ge h_{r,k,N}$, 
there are at least~$N$ elements~$y$ of height~$h$ in~$\rmC_A^+[r]\minus P$ such that~$P\cup \{y\}$ is still independent. 
%\nnote{version 2} For any~$k\ge 1$, there exist~$h_k\ge 1$ and~$w_k:\bbN\to \bbN$  strictly increasing for~$h\ge h_k$, such that for any 
%non-generating independent set~$P\in \calI_0[r]$ of cardinality~$k$, 
%there are at least~$w_k(h)$ elements~$x$ of height~$h$ in~$\rmC_A^+[r]\minus P$ such that~$P\cup \{x\}$ is still independent. 
%is at least~$w(h)$ for~$w$ a function depending only on~$k$
\end{lemma}

\begin{proof}
If we were only concerned with one particular non-generating independent set~$P$, we could produce such a height simply by picking one element~$y\in \rmC_A^+[r]$ independent of~$P$, which exists by the non-generating assumption, and considering the~$A^*$--subset~$\rmC_A^+(y)\subset\rmC_A^+[r]$ generated by~$y$. If~$y$ has height~$h_0$, then at height~$h=h_0+\ell$, this~$A^*$--subset has~$n^\ell$ elements, all independent of~$P$. And the number~$n^\ell$ can be chosen as big as we like by increasing~$\ell$ appropriately. To prove the lemma, we need to show that we can find a height that works for {\em any} given~$P$ of some fixed cardinality~$k$. 
%To do this, we will compute the number of elements that are independent of a fixed independent set~$P$ as a function of the height of its elements, and see that, as is intuitively true, the worst thing that can happen is if~$P$ has its vertices as low as possible. 

Note that it is enough to find a height $H$ such that for any independent set $P$ of cardinality $k$, there is at least~1 element $y$ of height $H$ independent of $P$. (Here $y$ will depend on $P$.)  Indeed, given such an $H$, just like in the case of a fixed set $P$, for all $h\ge H$, there will then be at least at least $n^{h-H}$ elements independent of $P$, namely the descendents of $y$ of height $H$. Setting $h_{r,k,1}=H$ and more generally $h_{r,k,N}=H+\ell$ for $\ell$ such that~\hbox{$n^\ell\ge N$}, will then prove the result.

We claim that we can take 
\[
H=\left\lceil\frac{k-r+1}{n-1}\right\rceil. %\min\{h\ |\ k\le (r-1)+h(n-1)\}.
\] 
Indeed, an element~$y\in \rmC_A^+[r]$ of height~$H$ is dependent of $P$ if there is a~$p\in P$ such that either~$y\in \rmC_A^+(p)$ or~\hbox{$p\in \rmC_A^+(y)$}. There are~$rn^{H}$ elements of height~$H$ in~$\rmC_A^+[r]$.
Suppose~$P=\{p_1,\dots,p_k\}$ with~$p_i$ of height~$h_i$, where we have ordered the elements of~$P$ so that~$h_1\le \dots\le h_i<H\le h_{i+1}\le \dots\le h_k$. Then 
there are~\hbox{$rn^{H}-(n^{H-h_1}+\dots+n^{H-h_i}+k-i)$} elements of height~$H$ that are independent of~$P$. This number is lowest when the heights of the elements of~$P$ are lowest, where lowest means as many elements as possible of height 0, then as many as possible of height one, and so on. 
It therefore remains to be shown that, if we choose a set $P$ with lowest heights, then there is still an element of height $H$ that is independent of $P$.
As~$P$ is assumed to be a non-generating independent set, it has at most~$(r-1)$ elements of height 0, and if~$P$ indeed has~$(r-1)$ elements of height 0, then it can have at most~$(n-1)$ elements of height 1, etc. (A set~$P$ with lowest height is displayed in  Figure~\ref{fig:indep}.) 
This shows that there will always be an element left in height $H$ when we have
\[
k\leqslant(r-1)+(n-1)H.
\]
By definition of~$H$, this inequality is satisfied, and we are done.
%such a lowest~$P$ will have at least one vertex independent from it at height~$H$. Any other~$P$ of the same cardinality will have at least as many elements of height~$H$ independent from it. 
\begin{figure}[ht]
\centering
\includegraphics[width=0.5\textwidth]{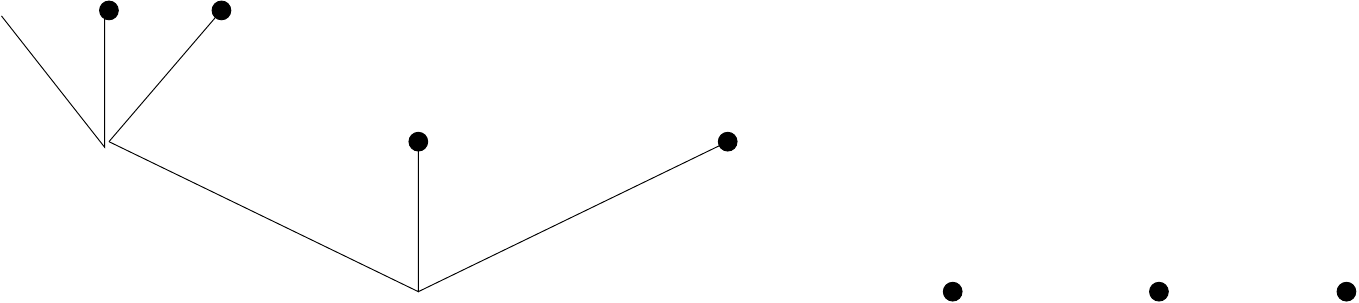}
\caption{Non-generating independent set~$P\subset\rmC^+_{A}[4]$ with~$|A|=3$ and~$|P|=7$ with most possible elements of height 0, then of height 1, etc. In this case, we see that~$h_{4,7,1}=2$ while~$h_{4,7,2}=3$.}
\label{fig:indep}
\end{figure}
\end{proof}

Note that we have defined the height~$h_{r,k,N}$ as the height such that there are at least~$N$ elements~$y$ of that height  with~$P\cup \{y\}$ still an independent set of cardinality~$k+1$. If we want~$P\cup \{y\}$ to be {\em non-generating}, and hence again defining a simplex of~$\rmT^\infty_r$, we need~$N\ge 2$.

\begin{proof}[Proof of Proposition~\ref{contract1}]
We will show that, for all~$k\ge 0$, all maps~$\rmS^k\to\rmT_r^\infty$ from spheres into the space~$\rmT_r^\infty$ are null-homotopic. 
Let us be given a map~$f\colon \rmS^k\to\rmT_r^\infty$. From \cite{Zeeman}, we can assume that the sphere~$\rmS^k$ comes with a triangulation such that the map~$f$ is simplicial. 
Let~$v_i$ be the number of simplices of dimension~$i$ is this triangulation of the sphere, and 
let~$v=v_0+\dots+v_k$ be the total number of simplices of all dimensions of that triangulation. In particular, this triangulation has~$v_0\le v$ vertices. Let \hbox{$N=v+2$} and~$h=h_{r,k,N}$ be the corresponding height obtained in Lemma~\ref{lem:height}. We will start by showing that the map~$f$ is homotopic to a map whose image only contains vertices of height at least~$h$ in~$\rmT_r^\infty$, and at most~$v$ of them: we will progressively retriangulate the sphere, and the new triangulation at all time will still have at most~$v$ vertices. 

We call a simplex of the sphere~$\rmS^k$~{\it bad} for~$f$ if all of its vertices are mapped to vertices in~$\rmT_r^\infty$ that have height less than~$h$.
%If there is a bad simplex for~$f$ with respect to the given triangulation of the sphere, then we will see that we can find a homotopy from~$f$ to a map that is simplicial with respect to another triangulation of the sphere that still has at most~$v$ vertices and that has fewer~bad simplices. 
We will modify~$f$ by removing the bad simplices inductively starting by those of highest dimension. 
So let~$\sigma$ be a bad simplex of maximal dimension~$p$ among all bad simplices. We will modify~$f$ and the triangulation of the sphere in the star of that simplex in a way that does not add any new bad simplex. In the process, we will increase the number of vertices by at most~$1$, and not at all if~$\sigma$ was a vertex. 
This implies that, after doing this for all bad simplices, we will have increased the number of vertices of the triangulation of the sphere by at most~$v_1+\dots+v_k$. As the sphere originally had~$v_0$ vertices, at the end of the process its new triangulation will  have at most~$v=v_0+v_1+\dots+v_k$ vertices (as we never introduced any new bad simplices, so only bad simplices of the original triangulation affect that number). There are two cases:

{\em Case~$p=k$.} If the bad simplex~$\sigma$ is of the dimension~$k$ of the sphere~$\rmS^k$, then its image~$f(\s)$ is a non-generating independent subset of~$\rmC_{A}^+[r]$. Because it is non-generating, we can choose~$y\in \rmC_A^+[r]$ that has height at least~$h$ and that is not a descendant of any vertex of~$f(\s)$ and still, together with~$f(\s)$, gives a non-generating independent subset. As the union~$f(\sigma)\cup\{y\}$ is again a simplex of~$\rmT_r^\infty$, we can add a vertex~$a$ in the center of~$\s$, replacing~$\s$ by~$\partial \s*a$ and replace~$f$ by the map~$(f|_{\partial\s})*(a\mapsto y)$ on~$\partial\sigma*a$. This map is homotopic to~$f$ through the simplex~$f(\s)\cup \{y\}$. We have added a single vertex to the triangulation. Because~$y$ has height~$h$, we have not added any new bad simplex, and we have removed one bad simplex, namely~$\s$. 

{\em Case~$p<k$.} If the bad simplex~$\sigma$ is a~$p$--simplex for some~$p<k$, by maximality of its dimension, the link of~$\sigma$ is mapped to vertices of height at least~$h$ in the complement of the free~$A^*$--set~$\rmC_A^+(f(\s))$ generated by~$f(\sigma)$. %These must lie in the the~$A^*$-invariant subset generated by the complement of the minimal basis containing~$f(\s)$. 
The simplex~$\s$ has~$p+1$ vertices whose images form an independent subset of~$\rmC_{A}^+[r]$ of cardinality at most~$p+1\leqslant k$.
By our choice of $h$, there are at least~$N= v+2$ vertices~$y_1,\dots,y_N$ of height~$h$ %in the complement of~$\rmC_{A}^+(f(\s))$, that is vertices 
such that each~$f(\s)\cup\{y_i\}$ still forms an independent set. 
Consider the~$A^*$--set 
$$\rmC_A^+(y_1)\sqcup\dots\sqcup\rmC_A^+(y_N)\ \subset \rmC_A^+[r].$$
As there are fewer vertices in the link than in the whole sphere, and the whole sphere has at most~$v$ vertices, by the pigeonhole principle, the vertices of~$\Link(\s)$ are mapped to at most~$v$ of the subsets~$\rmC_A^+(y_i)$. As~$N=v+2$, 
there are at least two of the above vertices~$y_i$ and~$y_j$ of height~$h$  such that no vertex of the link is mapped to any of their descendants~$\rmC_A^+(y_i)\sqcup\rmC_A^+(y_j)$. As these vertices are 
to start with independent of~$f(\s)$, it follows that for any simplex~$\tau$ of~$\Link(\s)$, the union~$f(\tau)\cup f(\s)\cup\{y_i\}$ is an independent set. Moreover this set is also non-generating because~$y_j$ is still independent of it. Hence  each~$f(\tau)\cup f(\s)\cup\{y_i\}$ forms a simplex of~$\rmT_r^\infty$. We can then replace~$f$ inside the star
\[
\Star(\s)=\Link(\s)*\s\simeq \rmS^{k-p-1}*\rmD^p
\]
by the map~$(f|_{\Link(\s)})*(a\mapsto y_i)*(f|_{\partial \s})$ on 
\[
\Link(\s)*a*\partial\s\simeq \rmS^{k-p-1}*D^0*\rmS^{p-1}.
\]
which agrees with~$f$ on the boundary~$\Link(\s)*\partial\s$ of the star, and is homotopic to~$f$ 
through the map~$(f|_{\Link(\s)})*(a\mapsto y_i)*(f|_{\s})$ defined on 
\[
\Link(\s)*a*\s\simeq \rmS^{k-p-1}*D^0*\rmD^p.
\] 
Now~$\Link(\s)*a*\partial(\s)$  has exactly one extra vertex compared to the star of~$\s$, unless~$\s$ was just a vertex, in which case its boundary is empty, and it has the same number of vertices.
As~$y_i$ has height~$h$, we have not added any new bad simplices. Hence we have reduced the number of bad simplices by one because~$\s$ was removed. 

By induction, we can now assume that there are no bad simplices for~$f$ with respect to a triangulation with at most~$v$ vertices. With this assumption, we can 
cone off~$f$ just as we coned off the links in the above argument:  We have more than~$N=v+2$ vertices of height~$h$ in~$\rmT_r^\infty$, and at most~$v$ vertices in the sphere. These vertices are mapped to vertices of height at least~$h$, that is to descendants of the vertices of height~$h$. By the pigeonhole principle, we know that there are at least two vertices~$y_i$ and~$y_j$, of height~$h$ such that no vertex of the sphere is mapped to any of their descendants. Hence we can cone off the sphere using~$\{y_i\}$. Indeed, this~$\{y_i\}$ is disjoint and independent from the set~$f(\sigma)$ for every simplex~$\s$ of the sphere, ensuring that the union~$f(\s)\cup\{y_i\}$ still forms an independent set, and any such set is non-generating because~$y_j$ is independent of it. Hence every~$f(\s)\cup\{y_i\}$ defines a simplex of~$\rmT_r^\infty$ and we can cone off the sphere by adding a single vertex mapped to~$y_i$. 
\end{proof}

\begin{proof}[Proof of Proposition~\ref{contract2}]
Consider a map~$f\colon\rmS^k\to\rmU_r^{\infty}$. Again we can assume that it is simplicial for some triangulation of the sphere~$\rmS^k$. We will show that there is a homotopy from~$f$ to a map that lands inside~$\rmT_r^\infty$. This will prove the result by Proposition~\ref{contract1}. 

We will, just like in the previous proof, modify~$f$ by a homotopy on the stars of the {\it bad} simplices in~$\rmS^k$, namely those whose vertices are all mapped to~$\rmU_r^\infty\backslash\rmT_r^\infty$. We will show how to reduce their number one by one, so that the result follows by induction.

Let~$\s$ be a bad simplex of maximal dimension, say~$p$, in the sphere~$\rmS^k$. By maximality, the link of~$\s$ is mapped to
%$\rmT_r^\infty\cap\operatorname{Star}(f(\s))$, which itself equals
~\hbox{$\rmT_r^\infty\cap\Link(f(\s))\subset \rmU_r^\infty$}.  
%as the vertices of~$f(\s)$ are mapped to~$\rmU_r^\infty\setminus\rmT_r^\infty$. 
We now argue that the intersection~\hbox{$\rmT_r^\infty\cap\Link(f(\s))$} is isomorphic to~$\rmT_{q}^\infty$ for some~$q\geqslant1$.

The simplex~$f(\sigma)$ is a collection of disjoint subsets of~$\rmC_{A}^+[r]$ that together form a non-generating independent subset~$P\in \calI_0[r]$. By Lemma~\ref{lem:minexp}, there exists a least expansion~$E$ of~$[r]$ containing~$P$, and because~$P$ is non-generating, the set~$E\minus P=Q$ has cardinality~$q\ge 1$. 
We claim that~$\rmT_r^\infty\cap\Link(f(\s))$ is isomorphic to~$\rmT_q^\infty$. 
To write a map, we first pick a bijection~$\la\colon[q]\to Q$.  Then~$\la$, together with the inclusion~$Q\subset\rmC_A^+[r]$ induces a map 
$$\hat\la\colon \rmC_A^+[q] \rar \rmC_A^+[r]$$
which in turn induces a map 
$$T(\hat \la)\colon \rmT_q^\infty \rar \rmT_r^\infty$$
defined on vertices by~$T(\hat\la)(y)=\hat\la(y)$. 
Indeed, the map~$\hat\la$ respects independency of subsets, which shows that~$T(\hat\la)$ is simplicial. 
Now note that any simplex of~$\rmT_r^\infty$ in the image of~$T(\hat \la)$ lies in the link of~$f(\s)$ because it necessarily is an independent strict subset of some expansion of~$Q$, and hence is independent of~$f(\s)$, and, together with $f(\s)$, non-generating. On the other hand, any simplex~$\tau$ of~$\rmT_r^\infty\cap\Link(f(\s))$ is defined by an independent set~$R$ such that~$R\cup P$ is (non-generating) independent. Let~$E'$ be some expansion of~$[r]$ containing~$R\cup P$. By minimality of~$E=Q\cup P$, we must have that~$E'$ is an expansion of~$E$ and hence that~$R$ lies in some expansion of~$Q$. It follows that~$\tau$ was actually in the image of~$T(\hat \la)$. 
Hence~$T(\hat\la)$ defines an isomorphism~\hbox{$\rmT_q^\infty\cong \rmT_r^\infty\cap\Link(f(\s))$}.

Let us consider the restriction of the map~$f$ to the star of the simplex~$\s$. Since we are working inside a~$k$-sphere, we have
\[
\Star(\sigma)\cong\Link(\sigma)*\sigma\cong \rmS^{k-p-1}*\sigma.
\]
The simplicial complex~$\rmT_{q}^\infty$ is contractible by Proposition~\ref{contract1}. As the link of~$f(\s)$ is mapped into the intersection~\hbox{$\rmT_r^\infty\cap\Link(f(\s))$}, which is isomorphic to~$\rmT_q^\infty$, we can extend the map from the link to a map~\hbox{$g\colon\rmD^{k-p}\to \rmT_r^\infty\cap\Link(f(\s))$}. 
Using such a map, the restriction of the map~$f$ to the star of the simplex~$\s$ can be extended to get a map~\hbox{$\hat f=g*f\colon\rmD^{k-p}*\s\to\rmU_r^\infty$} on the ball~$\rmD^{k-p}*\s$. 
%with~$D^{k-p}$ mapped to~$T_{q}^\infty$. 
This extension defines a homotopy on the star of~$\s$ in~$\rmS^k$, relative to the boundary link, from the restriction of~$f$ to a map with strictly fewer bad simplices: 
The new map~$g*f\colon\rmD^{k-p}*\partial\s\to\rmU_r^\infty$ has fewer old bad simplices, because~$\s$ has been removed, and 
there are no new bad simplices, because new simplices are joins of faces of~$\partial\s$ and simplices in the disc~$\rmD^{k-p}$, but the latter is mapped to good vertices. The result follows by induction.  
\end{proof}

%%%

% !TEX root = allinone_main_v2.tex

\section{The classifying space for the Cantor groupoid}\label{sec:classifying}

As before, we fix a finite set~$A=\{a_1,\dots,a_n\}$ of cardinality~$n\ge 2$. 
Recall from Section~\ref{sec:catcan} the groupoids~$\Set^\times$  of finite sets~$[r]=\{1,\dots,r\}$, with~$r\ge 0$, and their bijections,  and the groupoid~$\Can_A^\times$  of free Cantor algebras~$\rmC_A[r]$ of type~$A$, with~$r\ge 0$, and morphisms their isomorphisms. The groupoids~$\Set^\times$ and~$\Can_A^\times$ are both permutative categories, and taking the free Cantor algebra on a set defines a symmetric monoidal functor
\[
\rmC_A\colon\Set^\x\to \Can_A^\x.
\]
%The purpose of this section is to describe the homotopy type of the classifying space~$|\Can_A^\times|$ of the symmetric monoidal groupoid~$\Can_A^\times$ and its associated infinite loop space 
%$$\Omega^\infty \bbK(\Can_A^\times)$$
%in terms of the of the classifying spaces~$|\Set^\x|$.
The unit of the permutative category~$\Set^\x$ is the empty set~$[0]=\emptyset$ and that of~$\Can_A^\x$ is the empty Cantor algebra~$\rmC_A[0]$. Both have no non-trivial automorphisms and behave like disjoint basepoints:
\[
\Set^\x\ = \ \bSet^\x \sqcup \{[0]\} \ \ \ \textrm{and}\ \ \ \Can_A^\x\ =\ \bCan_A^\x\sqcup \{\rmC_A[0]\}
\]
for~$\bSet^\x$ the full subgroupoid of~$\Set^\x$ on the objects~$[r]$ with~$r>0$, and~$\bCan_A^\x$ the full subgroupoid of~$\Can_A^\x$ on the objects~$\rmC_A[r]$ with~$r>0$. 
Following Thomason \cite{Thomason}, we will in this section discard these units and work with non-unital permutative categories like~$\bSet^\x$ and~$\bCan^\x_A$,  
inserting units back again as disjoint basepoints only at the very end, just before taking spectra, as we explain now. 

There are algebraic K-theory ``machines'' that produce for every unital symmetric monoidal category~$\bfC$ a spectrum~$\bbK(\bfC)$ whose underlying infinite loop space
$\Omega^\infty\bbK(\bfC)$
is a group completion of the classifying space~$|\bfC|$ of~$\bfC$. (See~\cite[App C]{Thomason} for the particular machine that we will be using here.) 
If~$\bfC$ is a non-unital symmetric monoidal category, we will, following Thomason \cite[App C]{Thomason}, define 
$$\bbK(\bfC):=\bbK(\bfC_+)$$
for~$\bfC_+$ the category build from~$\bfC$ by adding a disjoint unit. In particular, when~$\bfC=\bSet^\x$, we have that~\hbox{$\bfC_+\cong \Set^\x$} and similarly for~$\bCan_A^\times$ and~$\Can_A^\x$. The Barratt--Priddy--Quillen theorem says that the spectrum
\[
\bbK(\bSet^\times)=\bbK(\Set^\times)\simeq \bbS
\]
is the sphere spectrum, so that there is a group completion~$|\Set^\times|\to\Omega^\infty\bbS$ (see \cite{BP72}). Our goal here is to compute~$\bbK(\Can_A^\x)$. We will do this using a permutative category~$\Tho_A$, build 
from~$\bSet^\times$ and the set~$A$, using a homotopy colimit construction of Thomason. Our main result for the section is the following: 

\begin{theorem}\label{back}
There is a strict symmetric monoidal functor~$\Tho_A \rar \bCan_A^\times$ of non-unital permutative categories that induces a homotopy equivalence~$|\Tho_A|\stackrel{\simeq}\rar |\bCan_A^\times|$ between the classifying spaces. 
\end{theorem}

The category~$\Tho_A$ will be build as a non-unital permutative category, and can be made unital by adding a disjoint basepoint. By a strict symmetric monoidal functor between non-unital permutative categories~$\bfC$ and~${\bf D}$, we will mean a functor~$F\colon\bfC\to {\bf D}$ such that~$F(X\oplus Y)=F(X)\oplus F(Y)$ for every objects~$X,Y$ of~$\bfC$ and likewise for morphisms. Adding disjoint units will give an associated strict symmetric monoidal functor~$F_+\colon\bfC_+\to {\bf D}_+$ in the usual sense. 

Using~\cite[Lem.~2.3]{Thomason}, we get the following corollary: 

\begin{corollary}\label{cor:KCan=KTho}
There is an equivalence~$\bbK(\Tho_A) \stackrel{\simeq}\rar \bbK(\bCan_A^\times)=\bbK(\Can_A^\x)$ of spectra.
\end{corollary}

The idea behind the permutative category~$\Tho_A$ is as follows. 
%First of all, there is a symmetric monoidal functor from the symmetric monoidal category~$\Set$ to~$\Can_A$ that sends a set~$[r]$ to the free Cantor algebra~$\rmC_A[r]$. It restricts to a symmetric monoidal functor~\hbox{$\rmC_A\colon\Set^\times\to\Can_A^\times$}. 
There is a symmetric monoidal endofunctor
\[
\Sigma_A\colon\bSet^\times\rar\bSet^\times
\]
that ``crosses with~$A$" (see Section~\ref{def:Sigma_A} for a precise definition). %takes a set~$[r]$ to the set~$[r]\times A$ and similarly for maps. %(Recall that~$A$ was any choice of set with~$a$ elements.) 
Now the functor~$\rmC_A\colon\bSet^\x\to \bCan_A^\x$ has the property that it is insensitive to pre-composition with~$\Sigma_A$: as seen in Example~\ref{ex:iso}, there are isomorphisms~\hbox{$\rmC_A(X\times A)\cong\rmC_A(X)$}, and these isomorphisms are essentially the defining property of Cantor algebras of type~$A$. There is in fact a natural isomorphism~\hbox{$\rmC_A\circ\Sigma_A\cong\rmC_A$} of functors. This suggests that the functor~$\rmC_A$ extends over a ``mapping torus of~$\Sigma_A$'' and the category~$\Tho_A$, which we define below, will be precisely such a device. 

%%%

Given a diagram~$\bfC\colon\lambda\mapsto\bfC_\lambda$ of non-unital permutative categories, indexed
on a small category~$\Lambda$, Thomason~\cite[Const.~3.22]{Thomason} defines a new non-unital  permutative category, denoted here by~$\hocolim_{\Lambda}\bfC_\lambda$, with the property that
\begin{equation}\label{eq:thmTho}
\bbK(\hocolim_{\Lambda}\bfC_\lambda)
\simeq
\hocolim_{\Lambda}\bbK(\bfC_\lambda).
\end{equation}
(See~\cite[Thm.~4.1]{Thomason}.) 
In order to construct~$\Tho_A$, we take the diagram that is indexed by the monoid~$\bbN$ of natural numbers, thought of as a category with one object. Given a pair~$(\bfC,F)$ consisting of a non-unital permutative category~$\bfC$ together with a symmetric monoidal endo-functor~$F$, we get a diagram on~$\bbN$ that associates to the unique object the category~$\bfC$ and to the morphism~$k\in \bbN$ the functor~$F^k$. We  define
\begin{equation}\label{eq:defTho}
\Tho_A=\hocolim_{\bbN}(\bSet^\times,\Sigma_A)
\end{equation}
for $(\bSet^\times,\Sigma_A)$ the corresponding diagram on $\bbN$. 
By the universal property of Thomason's construction (see~\cite[Prop 3.21]{Thomason}), %or~\cite[p.~337]{Thomason1979}), 
symmetric monoidal functors from~$\Tho_A$ to~$\bCan_A^\times$ can be defined by the following data:~a symmetric monoidal functor~$F\colon\bSet^\times\to\bCan_A^\times$ and a symmetric monoidal natural transformation~\hbox{$F\circ\Sigma_A\to F$}. Taking~$F=\rmC_A$ will thus give rise to a symmetric monoidal functor~\hbox{$\Tho_A\to\bCan_A^\times$}. We will construct this functor explicitly, and show that it induces an equivalence on the level of classifying spaces by showing that it fits inside a diagram
of categories and functors 
\[
\xymatrix@C=30pt{\Tho_A \ar[rrr] & & &  \bCan_A^\times\\
& \Lev_A \ar[ul]^\sim\ar[r]_\sim & \ar[ur]_\sim \Exp_A &}
\]
where all the other arrows induce homotopy equivalences on classifying spaces, and such that the diagram  commutes up to homotopy. 
The idea behind the intermediate categories~$\Exp_A$ and~$\Lev_A$ is as follows. Considering the functor~$\Tho_A\to\bCan_A^\times$, one sees that morphisms of~$\bCan_A^\times$ differ from morphisms of~$\Tho_A$ in two ways. Firstly, the morphisms of~$\Tho_A$ that map to the canonical isomorphism $\rmC_A[1]\to \rmC_A[n]$ are not invertible in~$\Tho_A$. Secondly, only certain simple types of expansions occur directly as morphisms of~$\Tho_A$, those that we will call~``level expansions''. We will take care of these two issues one at a time, writing the homotopy equivalence in three steps, studying the homotopy fibers each time.

The section is organized as follows: In Section~\ref{sec:Exp}, we define the category $\Exp_A$ and show that it is equivalent to $\bCan_A^\x$. In Section~\ref{levexpsec}, we define the category $\Lev_A$ and show that it is equivalent to $\Exp_A$. In Section~\ref{def:Sigma_A} we defined the category $\Tho_A$ and show that it is equivalent to $\Lev_A$. Finally, in Section~\ref{sec:square} we define the functor $\Tho_A\to \bCan_A^\x$ and show that the resulting diagram commutes.

\subsection{Expansions}\label{sec:Exp}

We will define the non-unital permutative category~$\Exp_A$ as a subcategory of~$\bCan_A^\x$, and show that passing from the groupoid~$\bCan_A^\times$ to its subcategory~$\Exp_A$ preserves the homotopy type of the classifying space. Here~$\bCan_A^\x$ denotes, as above, the category of free Cantor algebras~$\rmC_A[r]$ with~$r>0$ with morphisms the isomorphisms of Cantor algebras.  

Recall from Section~\ref{sec:isos} that 
a morphism from~$\rmC_A[r]$ to~$\rmC_A[s]$ in~$\bCan_A^\times$ can be described by a triple~$(E,F,\lambda)$, with the set~$E$ an expansion of~$[r]$, the set~$F$ an expansion of~$[s]$, and~$\lambda\colon E\to F$ a bijection. Recall also that there is a minimal such representative, where the minimality is defined using the poset structure of the set~$\calE[r]$ of expansions of $[r]$. 
Note that if either $E=[r]$ or if $F=[s]$, the representative is necessarily the minimal one. 
This allows for the following definition:

%\begin{definition}
%The {\em expansion category}~$\Exp_A$ has objects the finite sets and morphism~$X\to Y$ given by pairs~$(E,\lambda)$ of an expansion~$E$ of~$X$ and a bijection~$\lambda\colon E\to Y$. This expansion category can be identified with the subcategory of~$\Can_A^\times$, with the same objects, and with morphisms that can be presented as~$(E,Y,\lambda)$ with~$F=Y$. This description also makes the composition obvious.
%\end{definition}

\begin{definition}
The {\em expansion category}~$\Exp_A$ is the subcategory of~$\bCan_A^\times$, with the same objects, and with morphisms~$\rmC_A[r]\to\rmC_A[s]$ the morphisms of~$\bCan_A^\x$ that can be represented by a triple~$(E,F,\lambda)$ with~$F=[s]$. 
We will write~$(E,\lambda)$ for such a morphism. The symmetric monoidal structure of~$\bCan_A^\x$ restricts to one on~$\Exp_A$, making it a non-unital permutative category. 
\end{definition}

Note that $\Exp_A$ is indeed a subcategory: if $f:\rmC_A[r]\to \rmC_A[s]$ and $g:\rmC_A[s]\to \rmC_A[t]$ in $\bCan_A^\x$  are represented by $(E,[s],\la)$ and $(F,[t],\mu)$ respectively, then their composition $g\circ f$ is represented by $(\hat E,[t],\mu\circ\hat\la)$ for~\hbox{$\hat E=\hat\la^{-1}(F)$} the expansion of $E$ corresponding to $F$ under $\la$.

\begin{rem}
Recall from Section~\ref{sec:exp} that an expansion~$E$ of~$[r]$ can be thought of as a planar~$n$--ary forest on~$r$ roots corresponding to the elements of~$[r]$, whose leaves identify with the elements of~$E$. A bijection~\hbox{$\la\colon E\to[s]$} can then be interpreted as a labeling of the leaves of this forest. 
Identifying the objects of~$\Exp_A$ with the natural numbers and the the morphisms with labeled planar forests in this way, we thus see that~$\Exp_A$ can be thought of as a certain~``cobordism category'' of sets where the cobordisms are planar~$n$--ary forests. 
\end{rem}

Our main result in this section is the following: 

\begin{proposition}\label{prop:expcan}
The inclusion~$I\colon\Exp_A\to\bCan_A^\times$ is a symmetric monoidal
functor of non-unital permutative categories that induces an equivalence
\[
|\Exp_A|\simeq|\bCan_A^\times|
\]
on classifying spaces. 
%of~$\rmE_\infty$--monoids.
\end{proposition}

\begin{rem}
The above result can be seen as a special case of  \cite[Prop 2.13]{Thumann} by interpreting the Higman--Thompson groups as {\em operad groups} and applying the results in Section 3 of that paper. 
As the preparatory work for our proof will be useful in the following sections, we decided to keep the proof as is rather than explaining this alternative approach. 
\end{rem}

\begin{notation}[Induced maps~$\hat \la$]\label{not:hatla}
Given a bijection~$\la\colon X\stackrel{\cong}\to Y$ between finite sets, we get an induced isomorphism  
$$\rmC_A(\la)\colon\rmC_A(X)\stackrel{\cong}\rar \rmC_A(X)$$ of Cantor algebras, which in turn induces an isomorphism of posets 
$$\calE(\la)\colon\calE(X)\stackrel{\cong}\rar \calE(Y)$$ as~$\rmC_A(\la)$ takes expansions to expansions. And if~$E$ is an expansion of~$X$ and~$F=\calE(\la)(E)$ is the expansion of~$Y$ corresponding to~$E$ under~$\la$,  then restricting 
$\rmC_A(\la)$ to~$E$ induces a bijection~$E\to F$. In what follows, to ease notations, we will write~$\hat \la$ for these three types of maps induced by~$\la$: the map of Cantor algebra~$\rmC_A(\la)$, its restriction to subsets of~$\rmC_A(X)$, and the poset map~$\calE(\la)$. 
\end{notation}

Recall that~$A=\{a_1,\dots,a_n\}$ is an ordered set and that any expansion~$E$ of~$[r]$ is a subset 
$$E\subset \bigsqcup_{n\ge 0}[r]\x A^n\ =\ \rmC_A^+[r].$$
Hence~$E$ is also canonically ordered, using the lexicographic ordering on the words~$\bigsqcup_{n\ge 0}[r]\x A^n$.  (In terms of forests, this corresponds to the naturally induced planar ordering.) 
In the proof of Proposition~\ref{prop:expcan}, as well as later in Section~\ref{sec:classifying}, we will use this lexicographic ordering to identify any expansion~$E$ with the set~$[e]$ of the same cardinality as~$E$. 
The following result shows that this chosen identification is compatible with taking further expansions. 

\begin{lem}\label{lem:lexi}
Let~$A=\{a_1,\dots,a_n\}$. For any expansion~$E$ of~$[r]$ of cardinality~$e$, denote by~$\la_E\colon E \stackrel{\cong}\rar [e]$ the bijection defined by the lexicographic ordering of~$E$. 
Let~$F$ be an expansion of~$[e]$, with~$E'=\hat\la_E^{-1}(F)$ the corresponding expansion of~$E$ under~$\la_E$. Then~$\la_F\circ \hat\la_E=\la_{E'}$, i.e.~the upper triangle in the following diagram commutes:  
\begin{equation}\label{equ:EF}
\xymatrix@R-1pc{E'\ar@{}[d]|{\IV} \ar[rr]^-\cong_-{\hat \la_E} \ar@/^2pc/[rrrr]^-\cong_-{\la_{E'}} && \hat\la_E(E')=F \ar@{}[d]|{\IV} \ar[rr]^-\cong_-{\la_F} && [e']\\
 E \ar@{}[d]|{\IV}\ar[rr]^-\cong_-{\la_E} &&  [e]& & \\
[r] &  & &&
}\end{equation}
\end{lem}

\begin{proof}
The compatibility property follows from the fact that the map~$\la_E$ also induces a bijection of the set~$\rmC_A^+(E)\subset \rmC_A^+[r]$ with~$\rmC_A^+[e]$, and that this map respects the lexicographic order: this is true by definition of~$\la_E$ on the generating set~$E$, and by definition of the lexicographic order on the remaining elements. So the lexicographic order of~$E'\subset \rmC_A^+(E)\subset \rmC_A^+[r]$ agrees with that of~$\hat\la_E(E')\subset \rmC_A^+[e]$. 
%Through the map~$\rmC_A^+(\la_E)$, this induces a map of posets
%$$\hat \la_E:\calE(E)\rar \calE[e].$$
%\nnote{should we use such a notation for it?}
%The map~$\hat \la_E$ is compatible with the lexicographic orderings in the sense that, if~$E'\ge E\ge [r]$ is mapped to~$F\ge [e]$ by~$\hat \la_E$, then the composition 
%\begin{equation}\label{equ:EF}
%E'\ \stackrel{\rmC_A^+\la_E}\rar\ \hat\la_E(E')=F\ \stackrel{\la_F}\rar\  [e']
%\end{equation}
%equals the canonical map~$\la_{E'}:E'\to [e']$. This is because, whether we are working within expansions of~$[r]$ to compute~$\la_{E'}$ or expansions of~$[e]$ to compute~$\la_F$, an expansion replaces an element~$x$ by the~$n$ elements~$xa_1,\dots,xa_n$, and these are inserted instead of~$x$ in the ordering, wherever this expansion is considered. 
\end{proof}

\begin{proof}[Proof of Proposition~\ref{prop:expcan}]
We have already seen that the categories are (non-unital) permutative and that the inclusion respects this structure. 
The result will follow Quillen's Theorem A~\cite[\S1]{Quillen} if we show that, for any~$r>0$,  the fiber of~$I$ under~$\rmC_A[r]$ is equivalent to the poset~$\calE[r]$. Indeed, the poset~$\calE[r]$ has~$[r]$ as least element, and hence is contractible.

The fiber~$\rmC_A[r]\setminus I$ of the functor~$I$ under an object~$\rmC_A[r]$ of~$\bCan_A^\times$  has objects the pairs~$(\rmC_A[s],f)$ with~$f\colon\rmC_A[r]\to\rmC_A[s]$ a morphism in~$\bCan_A^\x$, and morphisms~$(\rmC_A[s],f)\to (\rmC_A[s'],f')$ given by morphisms~$(F,\lambda)\colon \rmC_A[s]\to \rmC_A[s']$ in~$\Exp_A$ such that
%~\hbox{$g\circ I(E,\lambda)=f$}.
the diagram 
$$\xymatrix{\rmC_A[s]\ar[rr]^-{I(F,\lambda)} &&\rmC_A[s'] \\
& \ar[ul]^{f} \rmC_A[r] \ar[ur]_{f'} & }$$ 
commutes in~$\bCan_A^\x$. 
 
We define a functor~$\Lambda\colon\calE[r]\to \rmC_A[r]\setminus I$ on objects by~\hbox{$\Lambda(E)=(\rmC_A[e],f_E)$}, where~$e=|E|$ is the cardinality of~$E$, and the map~\hbox{$f_E\colon\rmC_A[r]\to \rmC_A[e]$} is represented by the triple~$(E,[e],\la_E)$ 
with~$\la_E$ as in Lemma~\ref{lem:lexi}. 
 %where~$e$ is the isomorphism~$\rmC_A(E)\to\rmC_A(X)$ described by the triple~$(E,E,\id)$, 
%and on morphisms by taking an expansion to the corresponding morphism of~$I/\rmC_A(X)$. 
If~$E'\ge E$ is an expansion of~$E$, we have a diagram
\begin{equation}\label{diag:1}
\xymatrix@R-1pc{E'\ar@{}[d]|{\IV} \ar[rr]^-\cong_-{\hat \la_E} \ar@/^2pc/[rrrr]^-\cong_-{\la_{E'}} && F \ar@{}[d]|{\IV} \ar[rr]^-\cong_-{\la_F} && [e']\\
 E \ar[rr]^-\cong_-{\la_E} &&  [e]& & %\\
%[r] &  & &&
}\end{equation}
where triangle commutes by Lemma~\ref{lem:lexi}. We  define~$\Lambda$ on the inequality~\hbox{$E\leqslant E'$} in $\calE[r]$ to be the morphism of~$\Exp_A$ defined by the pair~\hbox{$(F,\la_F)\colon\rmC_A[e]\rar \rmC_A[e']$} for~$F=\hat\la_E(E')$ as in (\ref{diag:1}). %the expansion of~$[e]$ corresponding to~$E'$ along~$\hat\la_E$. 
The fact that the diagram 
$$\xymatrix{\rmC_A[e]\ar[rr]^-{I(F,\la_F)} &&\rmC_A[e'] \\
& \ar[ul]^{(E,[e],\la_E)} \rmC_A[r] \ar[ur]_{(E',[e'],\la_{E'})}& }$$ 
commutes is the commutativity of the triangle in (\ref{diag:1}). 

We now define a functor~\hbox{$\Pi\colon \rmC_A[r]\setminus I\to\calE[r]$} in the other direction: Given an object~$(\rmC_A[s],f)$ in the fiber, with~$f\colon\rmC_A[r]\to\rmC_A[s]$ an isomorphism of Cantor algebras, let~\hbox{$(E,G,\mu)$} be the minimal representative of~$f$. We set~$\Pi (\rmC_A[s],f)=E$. On morphisms, we have no choice as the target category is a poset, but we have to check that, given a morphism~$(F,\lambda)\colon(\rmC_A[s],f)\to(\rmC_A[s'],f')$ in the fiber, the expansion~$\Pi(\rmC_A[s'],f')$ of~$[r]$ is an expansion of~$\Pi(\rmC_A[s],f)$. This follows from the fact that~\hbox{$f'= I(F,\lambda)\circ f$}, implying that also~\hbox{$f^{-1}= f'^{-1}\circ I(F,\lambda)$}  and the fact that~$I(F,\lambda)$ comes from $\Exp_A$:  
%minimality of~$\Pi(\rmC_A(Y),y)$ as a presentation of~$y$ and the fact that~$\Pi(\rmC_A(Z),z)$ is involved in a presentation of~$y$ coming from its description as a composition~\hbox{$y=z\circ I(E,\lambda)$}.
if~$f'$ has minimal presentations~$(E',G',\mu)$, then $f'^{-1}$ can be represented by~$(G',E',\mu'^{-1})$, and the composition~$f'^{-1}\circ I(F,\lambda)$ can be computed using the diagram: 
$$\xymatrix@R-1pc{ \hat F \ar@{}[d]|{\IV}\ar[rr]^-{\hat\la}&&   G' \ar@{}[d]|{\IV} \ar[rr]^-{\mu'^{-1}} && E' \ar@{}[d]|{\IV}\\
 F \ar@{}[d]|{\IV}\ar[rr]^-{\la} &&  [s']& & [r] \\
[s] &  & &&
}$$
%for~$\hat G=\hat F$ a common expansion of~$G$ and~$F$, and~$H=\hat \la(\hat F)$, and 
with the triple~$(\hat F,E',\mu'^{-1}\circ \hat\la)$ representing the composition.  
%Given that~$(F,G,\mu)$ was the least representative of~$f$,  the expansion of~$E$ corresponding to~$F'$ under~$\la$, with~$\hat\la:E'\to F'$ the induced map.
Given that this composition is equal to~$f^{-1}$, which has minimal presentation~$(G,E,\mu^{-1})$, we must have that $E'\ge E$, as required. 
%necessarily have that~$\hat E$ is an expansion of~$E'$, and~$H$ is an expansion of~$G'$. Now as~$G'$ is itself an expansion of~$[s']$, we see that in fact~$H\ge G'\ge [s']$ and hence we also have, pulling back along~$\la\circ \mu$, that~$\hat E\ge E'\ge E$, yielding in particular the required inequality~$E\le E'$. 
%and~$G'$ is an expansion of~$G$, the latter giving the needed inequality. 

The composition~$\Pi \Lambda$ is the identity: it is enough to check this on objects as~$\calE[r]$ is a poset, and~$\Lambda$ takes an expansion~$E$ to the object~$(\rmC_A[e],(E,[e],\la_E))$, itself taken back to~$E$ by~$\Pi$. On the other hand, the composition~$\Lambda\Pi$ takes an object~$(\rmC_A[s],f)$ with~$f$ minimally represented by~$(E,G,\mu)$ to~$(\rmC_A[e],f_E)$ for~$e=|E|$ and~$f_E$ represented par~$(E,[e],\la_E)$. Now note that the diagram in~$\Can_A^\x$
$$\xymatrix{\rmC_A[s]\ar[rr]^-{(G,[e],\la_E\circ\mu^{-1})} &&\rmC_A[e] \\
& \ar[ul]^{(E,G,\mu)} \rmC_A[r]  \ar[ur]_{(E,[e],\la_E)}& }$$ 
commutes, where we have written the morphisms in terms of representatives. 
Thus~$(G,\la_E\circ\mu^{-1})$ defines a morphism in~$\rmC_A[r]\setminus I$ from~$(\rmC_A[s],f)$ to~$\Lambda\Pi(\rmC_A[s],f)$. We check now that these morphisms assemble to a natural transformation between the identity functor on~$\rmC_A[r]\setminus I$ and the composition~$\Lambda\Pi$. Indeed, consider a morphism~$(F,\la)\colon(\rmC_A[s],f)\to (\rmC_A[s'],f')$ where~$f$ has minimal representative~$(E,G,\mu)$ and~$f'$ has minimal representative~$(E',G', \mu')$. We need to check that the following square commutes in~$\rmC_A[r]\setminus I$: 
$$\xymatrix{(\rmC_A[s],f) \ar@{<--}[drr]_{f}\ar[dd]_{(F,\la)} \ar[rrrr]^-{(G,\la_E\circ\mu^{-1})} &&&& (\rmC_A[e],f_E)=\Lambda\Pi(\rmC_A[s],f) \ar@<-4ex>[dd]^{\Lambda\Pi(F,\la)} \ar@{<--}[dll]^{f_E} \\
&& \rmC_A[r] && \\
(\rmC_A[s'],f') \ar@{<--}[urr]^{f'}\ar[rrrr]^-{(G',\la_{E'}\circ\mu'^{-1})} &&&& (\rmC_A[e'],f_{E'}) =\Lambda\Pi(\rmC_A[s],f) \ar@{<--}[ull]_{f_{E'}}
}$$
Now this square commutes in~$\Exp_A$ if and only if its image in~$\bCan_A^\x$ commutes, as~$\Exp_A$ is a subcategory of~$\bCan_A^\x$. But in~$\bCan_A^\x$, each triangle in the diagram commutes by the fact that the maps in the square are maps in the fiber~$\rmC_A[r]\setminus I$. As all the maps are isomorphisms of Cantor algebras, it follows that the outside square commutes, as needed. 
\end{proof}

%%%

\subsection{Level expansions}\label{levexpsec}

The category~$\Exp_A$ has morphisms given by pairs consisting of an expansion and a bijection. We will now decompose the expansions into simpler types of expansions, which we call {\em level expansions}, that can be described in terms of subsets. We will construct a category~$\Lev_A$, which we will show is equivalent to~$\Exp_A$, where the morphisms will now be given by level expansions and bijections.

\begin{Def}
An expansion~$E$ of~$X$ is called a {\em level expansion} if there exists a subset~$P$ of~$X$ such that~\hbox{$E=(P\x A)\cup Q$} as a subset of~$\rmC_A^+(X)$, where~$Q=X\setminus P$ is the complement of~$P$ in~$X$. 
%partition~$X=P\sqcup Q$ with~$P\cap Q=\emptyset$ so that~$E$ results from~$X$ by expanding the elements in~$P$ once each while leaving the complement~$Q$ as it is. In particular, there is a canonical bijection~$E\cong(A\times P)\sqcup Q$. Note that there is at most one such partition.
% for some subset~$E\subset X$, with~$(E\times A)\sqcup (X\setminus E)$ seen as a subset of 
%$X\lambdangle A\rangle$  by identifying the pair~$(e,a)$ with~$ea\in X\lambdangle A\rangle$ for each~$e\in E$. 
\end{Def}

Note that level expansions do not ``compose'' in the sense that if~$E$ is a level expansion of~$X$ and~$F$ a level expansion of~$E$, then~$F$, while still an expansion of~$X$, need not be a level expansion of~$X$. For that reason, we do not get a category of level expansions analogous to~$\Exp_A$ right away. To define the category~$\Lev_A$, we will start with a semi-simplicial set of level expansions, and then pass to its poset of simplices.

% old definition of L(X):
%
%\begin{definition}\label{def:L}
%Given a set~$X$, we let~$\rmL(X)$ denote the simplicial complex whose vertices are the expansions~$E$ of~$X$, and where a set of~$p+1$ expansions~$E_0,\dots,E_p$ forms a~$p$--simplex if, after re-ordering them, all of them are level expansions of~$E_0$, and there exists pairwise disjoint subsets~$P_1,\dots,P_p\subseteq E_0$ such that~$E_i$ is also the level expansion of~$E_{i-1}$ along~$P_i$ for each~$i=1,\dots,p$. 
%\end{definition}

\begin{definition}\label{def:L}
Given a finite set~$X$, we define~$\rmL(X)$ to be the semi-simplicial set with~$0$--simplices the expansions of~$X$ and  with~$p$--simplices the sequences 
$$E_0< E_1<E_2<\dots< E_p$$
in the poset~$\calE(X)$ satisfying that each~$E_i$ is a level expansion of the smallest set~$E_0$. 
The~$i$--th face map forgets~$E_i$: 
%there is a semi-simplicial set whose vertices are the expansions~$E$ of~$X$, and whose~$p$--simplices are the~$(p+1)$--tuples~$(E_0,\dots,E_p)$ of expansions~$E_j$ of~$X$ such that there exist pairwise disjoint subsets~$P_1,\dots,P_p\subseteq E_0$ with~$E_j$ the level expansion of~$E_{j-1}$ along~$P_j$ for each~\hbox{$j=1,\dots,p$}. (Then all of the expansions~$E_j$ are level expansions of~$E_0$, and a~$p$--simplex can be identified with~$(E_0,P)$, where~$E_0$ is an expansion of~$X$, and~$P=(P_1,\dots,P_p)$ is a~$p$--tuple of pairwise disjoint subsets of~$E_0$.) The face maps~$\partial_j$, for~$j=0,\dots,p$, are given by omitting the entry~$E_j$ from the~$(p+1)$--tuple~$(E_0,\dots,E_p)$, that is
\[
d_i(E_0<\dots<E_p)= \ E_0<\dots<\widehat{E}_i<\dots<E_p.
\]
%We define~$\rmL(X)$ as the associated simplicial complex, with the same vertices, the expansions of~$X$, and where a set of~$p+1$ vertices forms a~$p$--simplex if and only if there is a~$p$--simplex in the semi-simplicial set that has these as its vertices. A pair of expansions of~$X$ forms an edge in~$\rmL(X)$ if and only if one is a level expansion of the other. 
\end{definition}

Note that also the differential $d_0$ makes sense. In fact,  if~$E_i=(Q_i\x A)\cup E_0\minus Q_i$ and~$E_j=(Q_j\x A)\cup E_0\minus Q_j$ are both level expansions of~$E_0$ with~$E_j>E_i$, then one necessarily has that~$Q_i\subset Q_j$ and~$E_j$ is the level expansion of~$E_i$ as well. 
%along~$Q_j\minus Q_i$. So in a simplex~$E_0< E_1<E_2<\dots< E_p$ of~$\rmL(X)$, the set~$E_j$ is a level expansion of~$E_i$ whenever~$j>i$. 
%In fact,~$\rmL(X)$ is isomorphic to the following semi-simplicial set: 
This leads to the following alternative description of $\rmL(X)$: 

\begin{definition}\label{def:L'}
Given a finite set~$X$, we define~$\rmL'(X)$ to be the semi-simplicial set with~$0$--simplices the expansions of~$X$ and  with~$p$--simplices the tuples
$$ (E,(P_1,\dots,P_p))$$
with~$E$ an expansion of~$X$ and~$(P_1,\dots,P_p)$ a collection of disjoint non-empty subsets of~$E$. The face maps are defined by
\[
d_i(E,(P_1,\dots,P_p)) =\left\{\begin{array}{ll} \big((P_1\x A)\cup (E\minus P_1), (P_2,\dots,P_p)\big) & i=0\\
 (E,(P_1,\dots,P_i\cup P_{i+1},\dots,P_p)) & 0< i <p\\
 (E,(P_1,\dots,P_{p-1})) & i=p.
\end{array}\right.
\]
\end{definition}

\begin{lem}\label{lem:corr}
The map~$\Phi\colon\rmL'(X)\to \rmL(X)$ that takes~$(E,(P_1,\dots,P_p))$ to~$E<E_1<\dots<E_p$ with the sets~$E_i$ defined by~\hbox{$E_i=(P_{[1,i]}\x A)\cup (E\minus P_{[1,i]})$} for~$P_{[1,i]}=P_1\cup\dots\cup P_i$, is an isomorphism of semi-simplicial sets. 
\end{lem}

\begin{proof}
The map~$\Phi$ is well-defined since each~$E_i$ is by definition a level expansion of~$E=E_0$, and the fact that~$E_i<E_{i+1}$ follows from the fact that~$P_{[1,i]}\subset P_{[1,i+1]}$; in fact,~$E_{i+1}$ is the level expansion of~$E_i$ along~\hbox{$P_{i+1}=P_{[1,i+1]}\minus P_{[1,i]}$}. The face maps in~$\rmL'(X)$ were defined to make this map simplicial:~$d_0$ corresponds to replacing~$E$ by~$E_1$, forgetting~$E=E_0$, and when~$0<i<p$, the map~$d_i$ corresponds to forgetting~$E_i$, and~$d_p$ forgets~$E_p$. 
Injectivity is immediate. For surjectivity, note that if~$E_0<\dots<E_p$ is a simplex of~$\rmL(X)$, then we must have that~$E_i=(Q_{i}\x A)\cup(E_0\minus Q_{i})$ for each~$i$, as~$E_i$ is a level expansion of~$E_0$, and the fact that~$E_i<E_{i+1}$ imposes that~$Q_{i}\subset Q_{i+1}$ as we have seen above. Setting~$P_i=Q_{i}\minus Q_{i-1}$ for each~$i$ gives a simplex~$ (E,(P_1,\dots,P_p))$ of~$\rmL'(X)$ mapping to~$E_0<\dots<E_p$, showing that~$\Phi$ is also surjective. 
\end{proof}

This lemma is in some way crucial, as it is the point where we switch from talking about expansions and morphisms of Cantor algebras, to solely talking about sets and subsets. 
Both ways of thinking of the simplices of~$\rmL(X)$ will be useful throughout the rest of the section. 

Note that any expansion~$E$ of a set~$X$ can be factored as a ``composition'' of level expansions, in the sense that one can always find expansions~$E_1,\dots,E_k$ of~$X$ such that 
$$X=E_0< E_1< \dots<  E_k<  E_{k+1}=E$$
in the poset~$\calE(X)$ of expansions of~$X$, in such a way that each~$E_i$ is a level expansion of~$E_{i-1}$. 
(There is even a canonical such factorization using the height filtration of~$\rmC_A^+(X)$, but we will not use it.) 
Given that~$\calE(X)$ is contractible, we get that~$\rmL(X)$ is connected for every finite set~$X$. The following result shows that the realization of~$\rmL(X)\cong \rmL'(X)$, which is a subspace of the nerve of~$\calE(X)$, is in fact,  like~$\calE(X)$,  contractible.

\begin{rem}\label{rem:semisimple}
Semi-simplicial sets have a realization, defined just like the ``thick'' realization of simplicial sets  (see eg.~\cite{EbertRW}). Given a poset~$P$, its nerve is a simplicial set, which can be considered as a semi-simplicial set by forgetting the degeneracies. This semi-simplicial set has~$q$--simplices the sequences~\hbox{$p_0\le \dots\le p_q$} in the poset~$P$, with the face map~$d_i$ forgetting the~$i$--th element. Alternatively, we can associate a smaller semi-simplicial  set to~$P$, 
that has~$q$--simplices the sequences~$p_0< \dots< p_q$ in the poset~$P$. Now the classical nerve of~$P$ can be recovered from this smaller nerve by freely adding all the degeneracies. The classical nerve as simplicial set or as semi-simplicial set (forgetting the degeneracies), and the smaller nerve using only strict inequalities, all have homotopy equivalent realizations (see eg.~\cite[Lem.~1.7, 1.8]{EbertRW}). In particular, if a poset has a least or greatest element, its realization is contractible, using whichever of these three possible realizations. 
%if~$L$ is a semi-simplicial set with~$L_p$ its set of~$p$--simplices, then its realization is the space~$\bigsqcup_{p\ge 0} L_p\x \Delta^p/\sim$, where~$(d_ix,t)\sim (x,d^it)$
\end{rem}

\begin{prop}\label{expcomplex}
For all finite sets~$X$, the semi-simplicial set~$\rmL(X)\cong \rmL'(X)$ is contractible. 
\end{prop}

\begin{proof}
If~$Y$ is an expansion of~$X$, we define~$\rmL(X,Y)$ to be the full subcomplex of~$\rmL(X)$ whose vertices are expansions~$E$ of~$X$ admitting~$Y$ as an expansion. As any finite collection of expansions of~$X$ admits a common expansion (by repeated use of Lemma~\ref{lem:common}), compactness of the spheres implies that every homotopy class can be represented by a map into some~\hbox{$\rmL(X,Y)\subseteq\rmL(X)$}. It is therefore sufficient to show that the complexes~$\rmL(X,Y)$ are contractible.

For an expansion~$E$ of~$X$, we define its {\em rank} to be~$rk(E)=(|E|-|X|)/(n-1)$. This is the number of simple expansions (expanding a single element~$x$ once) needed to obtain~$E$ from~$X$, and thus also the rank of~$E$ in the poset~$\calE(X)$. Suppose that the expansion~$Y$ has rank~$r$. Let~$F_i\rmL(X,Y)$ denote the full subcomplex of~$\rmL(X,Y)$ on the vertices of rank at least~$i$. This defines a descending filtration of~$\rmL(X,Y)$. 
\[
\rmL(X,Y)=F_0\rmL(X,Y)\supseteq F_1\rmL(X,Y)\supseteq F_2\rmL(X,Y)\supseteq\dots\supseteq F_r\rmL(X,Y)=\{Y\}
\]
For each~$i<r$,  the complex~$F_i\rmL(X,Y)$ is obtained  from~$F_{i+1}\rmL(X,Y)$ by attaching cones on the vertices~$E$ of rank~$i$ along their links, because no two such vertices are part of the same simplex. Now~$E_0<\dots<E_p$ is in~$\Link(E)\cap F_{i+1}(X,Y)$ if and only if 
$$E<E_0<\dots<E_p\le Y,$$
that is if and only if~$E_0<\dots<E_p$ is a sequence of non-trivial level expansions of~$E$ admitting~$Y$ as an expansion. 
Let $\calE_{\textrm Lev}(E,Y)$ denote the set of non-trivial level expansions of~$E$ admitting $Y$ as an expansion. Consider $\calE_{\textrm Lev}(E,Y)$ as a subposet of $\calE(E)$. We claim  that~\hbox{$\Link(E)\cap F_{i+1}(X,Y)$} is isomorphic to its associated (small) semi-simplicial set (in the sense of in Remark~\ref{rem:semisimple}). 
%Moreover,  the non-trivial level expansions of~$E$ admitting~$Y$ as an expansion form a subposet of~$\calE(E)$, and~$\Link(E)\cap F_{i+1}(X,Y)$  is precisely the thin nerve of that poset~(in the sense of Remark~\ref{rem:semisimple}): 
Indeed, if~$E_0,\dots,E_p\in \calE_{\textrm Lev}(E,Y)$ are level expansions of~$E$ such that $E_0<\dots<E_p$ is a simplex in the nerve of the poset $\calE$, and thus also by definition in the nerve of  $\calE_{\textrm Lev}(E,Y)$, then $E<E_0<\dots<E_p\le Y$ and we exactly have a simplex in~\hbox{$\Link(E)\cap F_{i+1}(X,Y)$}.

Our last step in the proof is to show that $\calE_{\textrm Lev}(E,Y)$ has a greatest element, and is hence contractible. 
Under the isomorphism of  Lemma~\ref{lem:corr}, $\calE_{\textrm Lev}(E,Y)$ is isomorphic to the poset~$\calP(E,Y)$ of $\emptyset\neq P\subset E$ such that~\hbox{$E(P)=(P\x A)\cup (E\minus P)$} admits~$Y$ as an expansion, where the poset structure is now inclusion.  Let~$\hat P=\bigcup_{P\in \calP(E,Y)} P$ be the union of all such~$P$'s. The check that $\hat P\in \calP(E,Y)$, we need to check that the expansion~$E(\hat P)=(\hat P\x A)\cup (E\minus \hat P)$ still admits $Y$ as an expansion. This follows from the fact that we can check this component-wise: writing~$Y=\cup_{e\in E}Y_e$ for~$Y_e=Y\cap \rmC_A^+(e)\subset \rmC_A^+[E]$, the condition that~$Y\ge E(\hat P)$ is equivalent to~$Y_e\in \rmC_A^+(e\x A)\subset \rmC_A^+(e)$ for every~$e\in \hat P$, which holds by the definition of~$\hat P$ as every such $e$ lies in an element $P$ of $\calP(E,Y)$. Hence $\hat P$ is a greatest element for $\calP(E,Y)$, showing that    $\calP(E,Y)$, and thus also~$\calE_{\textrm Lev}(E,Y)$, is contractible, as needed. The result follows by induction. 
\end{proof}

Let~$\calL(X)$ denote the poset of simplices of~$\rmL'(X)$, so the elements of~$\calL(X)$ are the simplices of~$\rmL'(X)$, and the morphisms are the inclusions among them. 
Explicitly, the elements of~$\calL(X)$ are the  tuples~$ (E,(P_1,\dots,P_p))$ with~$E$ an expansion of~$X$ and~$(P_1,\dots,P_p)$, for some~$p\ge 0$, a sequence of disjoint non-empty subsets of~$E$. 
To describe the poset structure, recall from Definition~\ref{def:L'} that the face map~$d_0$ in~$L'(X)$ takes a level expansion along~$P_1$, then for~$0<i<p$, the face map~$d_i$ takes the union of the neighbouring subsets~$P_i$ and~$P_{i+1}$, and finally the last face map~$d_p$ forgets the last subset~$P_p$. 
Writing
\[
P_{[i,j]}:=P_i\cup \dots\cup P_j
\]
for~$i\le j$, we thus have that
\[
(F,(Q_1,\dots,Q_q))\le (E,(P_1,\dots,P_p))
\]
if there exists~$0\le p_0,p_1\le p$ and~$\kappa_1,\dots,\kappa_q\ge 1$ with 
$p_0+\sum_{j=1}^q\kappa_j+p_1=p$
such that 
$$F=(P_{[1,p_0]}\x A)\cup (E\minus P_{[1,p_0]})$$ is the level expansion of~$E$ along the first~$p_0$ subsets and 
each~$Q_j$ is a union
%The poset~$\calL(X)$ can be identified with the poset of pairs~$(E,P)$, with~$E$ an expansion of~$X$ and~\hbox{$P=(P_1,\dots,P_p)$} a~$p$--tuple of pairwise disjoint subsets of~$E$, for some~$p\geqslant0$. We have~\hbox{$(E,P)\leqslant(F,Q)$} if there is an~$i$ such that~
%$F$  the level expansion of~$E$ along the union~$P_1\sqcup\dots\sqcup P_i$ for some~$i\ge 0$ (that is~$F=E_i$ of the above description), and 
\[
%Q_j=P_{i+\kappa(1)+\dots+\kappa(j-1)+1}\sqcup\dots\sqcup P_{i+\kappa(1)+\dots+\kappa(j)}
Q_j=P_{[p_0+\kappa_1+\dots+\kappa_{j-1}+1,p_0+\kappa_1+\dots+\kappa_j]}.
\]
In particular 
$Q_1\cup\dots \cup Q_q=P_{[p_0+1,p-p_1]}$. (The numbers~$p_0$ and~$p_1$ count the number of~$0$--th and last face maps that have been applied.)

Note that the order relation in the poset~$\calL[r]$ is opposite to the order relation in the poset~$\calE[r]$ in the sense that if~\hbox{$ ([r],(P_1,\dots,P_p))\ge(F,(Q_1,\dots,Q_q)$} in~$\calL[r]$, then~$F\ge [r]$ in~$\calE[r]$. This is why the order relation of~$\calL[r]$ will occur in the reversed direction in the definition of~$\Lev_A$ below. 

%%%

We have already seen that a bijection~$\lambda\colon X\to Y$ induces an isomorphism~$\hat\la\colon\calE(X)\to\calE(Y)$ of posets (see Notation~\ref{not:hatla}). This in turn induces a map~$\hat\lambda\colon \rmL(X)\to \rmL(Y)$ of semi-simplicial sets, and hence also  a map~$\hat\lambda\colon\calL(X)\to \calL(Y)$ of posets of simplices, which we will denote by~$\hat\la$ again. Explicitly, this last map takes~$ (E,(P_1,\dots,P_p))$ to~~$(\hat\la E,(\hat\la P_1,\dots,\hat\la P_p))$.
%, where, abusing notation as in Notation~\ref{not:hatla}, we use~$\hat\la$ to denote all the maps induced by~$\la$. 
%$(\hat{la}(E),(Q_1,\dots,Q_q))$
%As this maps comes from a map between the semi-simplicial sets and preserves the cardinality of expansions, in the subsets description of the elements of~$\calL(X)$, this map has the property that  if 
%$(E,(P_1,\dots,P_p))$ is mapped to~$(F,(Q_1,\dots,Q_q))$ by~$\calL(\lambda)$, then~$p=q$ and~$\lambda$ also induces a bijection~\hbox{$\lambda:E\to F$} with the property that~$\lambda(P_i)=Q_i$ for each~$i$. %, which we also write as~$\lambda(P)=Q$. 
Note also that for any expansion $E$ of $X$, the poset $\calL(E)$ identifies with the subposet of $\calL(X)$ with objects the tuples $(F,(P_1,\dots,P_p))$ such that $F$ is an expansion of $E$.

We are now ready to construct the category~$\Lev_A$ from the posets~$\calL(X)$ in a way similar to the one in which we constructed~$\Exp_A$ from the posets~$\calE(X)$. 
As in the case of~$\bCan_A^\x$ and~$\Exp_A$, we will restrict to objects build from the non-zero natural numbers. %in order to get a strict monoidal structure on the category. 
%, where we consider all objects~$(E,P)$ of the poset~$\calL(X)$ for all sets~$X$, and where morphisms are generated by bijections and ``level expansions,'' now interpreted as the poset structure of~$\calL(X)$ for some~$X$. 

\begin{definition}
The category~$\Lev_A$ has as objects the sequences~$([r],(P_1,\dots,P_p))$ for~$r>0$ and~$(P_1,\dots,P_p)$ a~$p$--tuple of~$p\geqslant0$ disjoint non-empty subsets of~$[r]$. (In particular, 
$([r],(P_1,\dots,P_p))\in \calL[r]$.) 
The morphisms
\[
\phi\colon([r],(P_1,\dots,P_p)) \ \rar \  ([s],(Q_1,\dots,Q_q))
\]
are given by tuples~$\phi=(F,(Q'_1,\dots,Q'_q),\lambda)$ where the element $(F,(Q'_{1},\dots,Q'_{q}))\in \calL[r]$ is such that the relation~$([r],(P_1,\dots,P_p))\geqslant (F,(Q'_{1},\dots,Q'_{q}))$  holds in~$\calL([r])$
and~\hbox{$\la\colon F\to [s]$} is a bijection such that~\hbox{$\la(Q'_{i})=Q_{i}$} for each~$i$. (In particular, $F$ is a level expansion of $[r]$ along some $P_i$'s and each $Q_i'$ is a union of $P_i$'s.) 
The composition of~$\phi$ with a morphism~\hbox{$\psi=(G,(R'_1,\dots,R'_u),\mu)\colon ([s],(Q_1,\dots,Q_q))\to ([t],(R_1,\dots,R_u))$}  is defined by 
$$(G,(R'_1,\dots,R'_u),\mu) \circ (F,(Q'_1,\dots,Q'_q),\lambda) = (\hat\la^{-1}G,(\hat\la^{-1}R'_1,\dots,\hat\la^{-1}R'_q),\mu\circ\hat\lambda)$$
where 
$$\xymatrix@R-1pc{\hat\la^{-1}(G) \ar[rr]^-{\hat \la} \ar@{}[d]|\IV && G \ar[rr]^-{\mu} \ar@{}[d]|\IV && [t]\\
F \ar[rr]^-{\la} \ar@{}[d]|\IV && [s]&&\\
[r] &&&& 
}$$
and where we use that~$\la$ induces an isomorphism~$\hat \la\colon\calL(F)\to \calL[s]$ and that~$\calL(F)$ is naturally a subposet of~$\calL[r]$, so that $ ([s],(Q_1,\dots,Q_q)) \ge (G,(R'_1,\dots,R'_u)$ in $\calL[s]$ gives 
\[
([r],(P_1,\dots,P_p))\ge  (F,(Q'_1,\dots,Q'_q)\ge  (\hat\la^{-1}G,(\hat\la^{-1}R'_1,\dots,\hat\la^{-1}R'_q))
\]
in~$\calL[r]$. (In particular, $\hat\la^{-1}(G)$ is necessarily again a level expansion of $[r]$ along some $P_i$'s and each $\hat\la^{-1}(R'_i)$ is a union of $P_i$'s.) 
\end{definition}

We see here that it is crucial for the composition that we use the poset structure of $\calL[r]$ and not the simpler notion of level expansion to define the morphisms in $\Lev_A$, as level expansions in general do not compose to level expansions.

Note that part of the data of a morphism~$\phi=(F,(Q'_1,\dots,Q'_q),\lambda)$ from~$([r],(P_1\dots,P_p))$ to~$([s],(Q_1,\dots,Q_q))$ in~$\Lev_A$ is an expansion~$F$ of~$[r]$ and a bijection~$\la\colon F\stackrel{\cong}\rar [s]$. In fact, there is a functor
$$J\colon \Lev_A\to \Exp_A$$
defined on objects and morphisms as a forgetful map: 
$$J([r],(P_1\dots,P_p)):=\rmC_A[r] \ \ \ \textrm{and}\ \ \ J(F,(Q'_1,\dots,Q'_q),\lambda)=(F,\la).$$
This is compatible with composition as can be seen directly from the definition of composition in~$\Lev_A$ above. 

\begin{proposition}\label{prop:levexp}
The  functor
\[
J\colon\Lev_A\to\Exp_A   %\,,\,([r],(P_1,\dots,P_p))\longmapsto \rmC_A[r]
\]
defined above induces an equivalence~$|\Lev_A|\simeq |\Exp_A|$ on classifying spaces.
\end{proposition}

\begin{proof}
We will show that the fiber of~$J$ under~$\rmC_A[r]$ is homotopy equivalent to~$\calL[r]$. As~$\calL[r]$ is the poset of simplices of~$\rmL'[r]$, the result will then follow from Proposition~\ref{expcomplex}, which says that~$\rmL'[r]$ is contractible. 

\noindent {\em Description of the fibers.} 
The objects of the fiber~$\rmC_A[r]\setminus J$ can be written as tuples 
$$([s],(P_1,\dots,P_p), [s] \stackrel{\ \la}\lar E)$$ 
for~$([s],(P_1,\dots,P_p))$ an object of the category~$\Lev_A$, the set~$E$ an expansion of~$[r]$ and~$\lambda\colon E\stackrel{\cong}\rar [s]$ a bijection defining a morphism 
$(E,\la)\colon \rmC_A[r]\to \rmC_A[s]=J([s],(P_1,\dots,P_p))$ in~$\Exp_A$. 
A morphism 
$$\phi\colon([s],(P_1,\dots,P_p), [s]\stackrel{\ \la}\lar E)\ \rar \ ([t],(Q_1\dots,Q_q), [t]\stackrel{\ \mu}\lar F)$$ 
in the fiber is given by a morphism
\[
\phi=(G,(Q_1',\dots,Q'_q),G\stackrel{\kappa}\to [t])\colon([s],(P_1,\dots,P_p))\ \rar \ ([t],(Q_1,\dots,Q_q))
\]
in the category~$\Lev_A$, with in particular~$G$ an expansion of~$[s]$ and~$\kappa$ a bijection (where we have spelled out the source and target of~$\kappa$ for clarity), such that the diagram 
$$\xymatrix{J([s],(P_1,\dots,P_p))=\rmC_A[s] \ar[rr]^-{J(\phi)=(G,\kappa)}&& \rmC_A[t]=J([t],(Q_1\dots,Q_q)) \\
&\rmC_A[r] \ar[ul]+<10ex,-2ex>^{(E,\la)} \ar[ur]+<-10ex,-2ex>_{(F,\mu)}&}$$
commutes in~$\Exp_A$. 
%~\hbox{$(F',\mu)=J(\phi)\circ(E',\lambda)\colon E\cong\rar [s]$} as morphisms~$\rmC_A(X)\to\rmC_A(F)$ in~$\Exp_A$.

Just as in the proof of Proposition~\ref{prop:expcan}, we define a pair of functors~$\Lambda\colon\calL[r]^{\op}\longleftrightarrow \rmC_A[r]\backslash J\colon\Pi$, where the order of~$\calL[r]$ is reversed. 

\noindent {\em Definition of~$\Lambda\colon \calL[r]^{\op}\longrightarrow \rmC_A[r]\backslash J$.} 
Define~$\Lambda$ on objects by 
\begin{gather*}
%\Pi(E,P,\lambda\colon E'\cong E)=(E',\lambda^{-1}P)\\
\Lambda(E,(P_1,\dots,P_p))=([e],(\la_E P_1,\dots,\la_EP_p),[e]\stackrel{\ \la_E}\lar E)
\end{gather*}
for~$e=|E|$ and~$\la_E\colon E\to [e]$ the map of Lemma~\ref{lem:lexi} given by the lexicographic ordering of~$E$. 
The functor~$\Lambda$ is given on morphisms as follows: if~$(E,(P_1\dots,P_p))\ge (E',(Q_1,\dots,Q_q))$ in~$\calL[r]$, we have that~$E'$ is an expansion of~$E$ along~$P_1\cup \dots\cup P_i$ for some~$i$ and 
we get a diagram of expansions and isomorphisms 
$$\xymatrix@R-1pc{[e'] & \ar[l]_-{\la_{E'}} E' \ar@{}[d]|{\IV} \ar[r]^-{\hat\la_E} & \hat\la_EE' \ar@{}[d]|{\IV}\\
& E \ar[r]^-{\la_E} & [e].
}$$
We define a morphism 
$$\phi\colon([e],(\la_EP_1,\dots,\la_EP_p),[e]\stackrel{\ \la_E}\lar E) \rar ([e'],(\la_{E'}Q_1,\dots,\la_{E'}Q_q),[e']\stackrel{\ \la_{E'}}\lar E')$$
in the fiber by setting~$\phi=(\hat\la_EE',\hat\la_E Q_1,\dots,\hat\la_EQ_q,\hat\la_EE'\xrightarrow{\la_{E'}\hat\la_E^{-1}} [e'])$. 
%So~$\phi$ is the image under~$\la_E$ of the poset morphism~$(E,(P_1\dots,P_p))\leqslant(E',(Q_1,\dots,Q_q))$, 
This makes sense because 
$$([e],(\la_EP_1,\dots,\la_EP_p))\ \ge \ (\hat\la_EE',(\hat\la_E Q_1,\dots,\hat\la_EQ_q))$$
in~$\calL[e]$, as it is 
the image under~$\hat \la_E\colon \calL(E)\subset \calL[r]\to \calL[e]$ of~$(E,(P_1\dots,P_p))\geqslant(E',(Q_1,\dots,Q_q))$ in~$\calL[r]$. 
Also, we have~$\la_{E'}\hat\la_E^{-1}(\hat\la_EQ_i)=\la_{E'}Q_i$ for each~$i$, so~$\phi$ defines a morphism in~$\Lev_A$ from~$([e],(\la_EP_1,\dots,\la_EP_p))$ to~$([e'],(\la_{E'}Q_1,\dots,\la_{E'}Q_q))$.  
%Finally, using Lemma~\ref{lem:lexi}, we have that~$\la_{E'}\circ \hat\la_E^{-1}:\hat\la_E(E')\to [e']$ identifies with the map induced by the lexicographic ordering. This gives the commutativity of 
Finally, we see that the diagram 
$$\xymatrix{\rmC_A[e] \ar[rr]^-{(\hat\la_EE',\la_{E'}\hat\la_E^{-1})} && \rmC_A[e']\\
&\rmC_A[r], \ar[ul]^{(E,\la_E)}\ar[ur]_{(E',\la_{E'})}&
}$$
commutes in~$\Exp_A$, 
showing that~$\phi$ indeed defines a morphism in the fiber. Functoriality follows from the following computation: 
if~$(E,(P_1\dots,P_p))\ge (E',(Q_1,\dots,Q_q))\ge (E'',(R_1,\dots,R_u))$ in~$\calL[r]$, the composition of the images of these morphisms under~$\Lambda$ is 
\begin{align*}
\big(\hat\la_{E'}E'',(\hat\la_{E'}R_1,&\dots,\hat\la_{E'}R_u),\la_{E''}\hat\la_{E'}^{-1}\big)\circ \big(\hat\la_EE',(\hat\la_EQ_1,\dots,\hat\la_EQ_q),\la_{E'}\hat\la_E^{-1}\big)\\
=&\big((\hat\la_{E'}\hat\la_E^{-1})^{-1}\hat\la_{E'}E'',((\hat\la_{E'}\hat\la_E^{-1})^{-1}\hat\la_{E'}R_1,\dots,(\hat\la_{E'}\hat\la_E^{-1})^{-1}\hat\la_{E'}R_u),(\la_{E''}\hat\la_{E'}^{-1})\circ (\hat\la_{E'}\hat\la_E^{-1})\big)\\
=&\big(\hat\la_{E}E'',(\hat\la_{E}R_1,\dots,\hat\la_{E}R_u),\la_{E''}\hat\la_E^{-1}\big), 
\end{align*}
which is equal to image under~$\Lambda$ of the composition. 
%this, together with~$\id_F$, defines a morphism~$\phi\colon(E,P)\to(F,Q)$ in~$\Lev_A$ such that~$(F,\id_F)=J(\phi)\circ(E,\id_E)$ in~$\Exp_A$. By what has been explained above, this~$\phi$ defines a morphism~$\Lambda(E,P)\to\Lambda(F,Q)$ in the fiber~$X\backslash J$.

\noindent {\em Definition of~$\Pi\colon \rmC_A[r]\backslash J\rar \calL[r]^{\op}$.} 
Define the functor~$\Pi$ on objects by 
\begin{gather*}
\Pi([s],(P_1,\dots,P_p),[s]\stackrel{\ \la}\lar E)=(E,(\lambda^{-1}P_1,\dots,\la^{-1}P_p))
%\Lambda(E,P)=(E,P,\id\colon E=E).
\end{gather*}
which makes sense as~$E$ is an expansion of~$[r]$ and~$\la$ a bijection. As the target of~$\Pi$ is a poset, we are left to check that the existence of a morphism in the source gives the appropriate inequality in the target. 
Given a morphism
$$\phi=(G,(Q_1',\dots,Q'_q),G\stackrel{\kappa}\to [t])\colon([s],(P_1,\dots.P_p),[s]\stackrel{\ \la}\lar E)\ \rar\ ([t],(Q_1,\dots,Q_q),[t]\stackrel{\ \mu}\lar F)$$ 
in the fiber~$\rmC_A[r]\backslash J$, 
%with~$\phi=(G,(Q_1',\dots,Q'_q),G\stackrel{\kappa}\to [t])\colon([s],(P_1,\dots,P_p)) \rar ([t],(Q_1,\dots,Q_q))$ the underlying morphism in~$\Lev_A$, 
we have $ ([s],(P_1,\dots.P_p))\ge (G,(Q_1',\dots,Q'_q))$
in~$\calL[s]$ and~$\kappa\colon G\to [t]$ an isomorphism taking~$Q_i'$ to~$Q_{i}$ for each~$i$.
Pulling back along the bijection~$\la\colon E\to [s]$ we get  that
\[
(E,(\lambda^{-1}P_1,\dots,\la^{-1}P_p))\ge  (\hat\la^{-1}G,(\hat\la^{-1}Q_1',\dots,\hat\la^{-1}Q'_q))
\]
in the poset~$\calL(E)$, which can be identified with a subposet of~$\calL[r]$. 
We claim that the right hand side is equal to~\hbox{$(F,(\mu^{-1}Q_1,\dots,\mu^{-1}Q_q))$}, which  will give the required inequality. 
Indeed, because $\phi$ is a morphism in the fiber, we have~\hbox{$(G,\kappa)\circ (E,\la)=(F,\mu)$}. Computing the left hand side in 
\[
\xymatrix@R-1pc{\hat\la^{-1}G \ar[r]^-{\hat\la} \ar@{}[d]|{\IV} & G \ar[r]^-{\kappa} \ar@{}[d]|{\IV} & [t]\\
E \ar[r]^-{\la} \ar@{}[d]|{\IV} & [s] & \\
[r],&&}
\]
shows that~$\hat\la^{-1}G=F$ and~$\kappa\circ \hat\la=\mu$. The claim then follows from the fact that~$Q'_i=\kappa^{-1}Q_i$. 
%\begin{align*}
% (\hat\la^{-1}G,(\hat\la^{-1}(Q_1'),\dots,\hat\la^{-1}(Q'_q))) &=  (\hat\la^{-1}G,(\hat\la^{-1}\circ \kappa^{-1}(Q_1),\dots,\hat\la^{-1}\circ \kappa^{-1}(Q_q))) \\
%&=(F,(\mu^{-1}(Q_1),\dots,\mu^{-1}(Q_q)))
%\end{align*}

\noindent {\em The functors~$\Lambda$ and~$\Pi$ define an equivalence.} 
The composition~$\Pi\Lambda$ is the identity on~$\calL^{\op}[r]$: it is enough to check this on objects as~$\calL^{\op}[r]$ is a poset. An object~$(E,(P_1,\dots,P_p))\in\calL^{\op}[r]$ is mapped by~$\Lambda$ to~$([e],(\la_EP_1,\dots,\la_EP_p),[e]\stackrel{\la_E}\lar E)$, which is 
mapped back to the original tuple by~$\Pi$. On the other hand, we have
\begin{align*}
\Lambda\Pi([s],(P_1,\dots,P_p),[s]\stackrel{\ \la}\lar E) %&=([e],(\la_E\lambda^{-1}P_1,\dots,\la_E\la^{-1}P_p),[e]\stackrel{\ \la_E}\lar E)\\
&=([s],(\la_E\lambda^{-1}P_1,\dots,\la_E\la^{-1}P_p),[s]\stackrel{\ \la_E}\lar E)
%(E,P,\id\colon E'\cong E),
\end{align*}
where we have used that~$e=s$ as~$E$ and~$[s]$ have the same cardinality. These objects of~$\rmC_A[r]\minus J$ are not equal in general, but 
$$\eta=([s],(P_1,\dots,P_p),[s]\xrightarrow{\la_E\la^{-1}} [s])$$ is a morphism in~$\Lev_A$ from~$([s],(P_1,\dots,P_p))$ to~$([s],(\la_E\lambda^{-1}P_1,\dots,\la_E\la^{-1}P_p)$, which induces a morphism in the fiber from~$([s],(P_1,\dots,P_p),[s]\stackrel{\ \la}\lar E)$ to its image by~$\Lambda\Pi$ by 
 the commutativity of the diagram 
$$\xymatrix{\rmC_A[s]\ar[rr]^-{([s],\la_E\la^{-1})} && \rmC_A[s] \\
&\rmC_A[r] \ar[ur]_{(E,\la_E)} \ar[ul]^{(E,\la)} & 
}$$
in~$\Exp_A$.  We are left to check that these morphisms~$\eta$ fit together to define a natural transformation from the identity to~$\Lambda\Pi$. Consider a morphism 
$$\phi=(G,(Q_1',\dots,Q'_q),G\stackrel{\kappa}\to [t])\colon([s],(P_1,\dots.P_p), [s]\stackrel{\ \la}\lar E)\ \rar\ ([t],(Q_1,\dots,Q_q),[t]\stackrel{\ \mu}\lar F)$$
in the fiber. Then~$\Pi(\phi)$ just remembers that the images under~$\Pi$ of the source and target are comparable in~$\calL[r]$, and 
\begin{align*}\Lambda\Pi(\phi) %&=(\hat\la_E\hat\la^{-1}(G),(\hat\la_E\hat\la^{-1}(Q_1'),\dots,\hat\la_E\hat\la^{-1}(Q_q')),\hat\la_E\hat\la^{-1}(G)\stackrel{\la_F\circ\hat\la_E^{-1}}\rar [s]) \\
&=(\hat\la_EF,(\hat\la_E\mu^{-1}Q_1,\dots,\hat\la_E\mu^{-1}Q_q),\hat\la_EF\xrightarrow{\la_F\hat \la_E^{-1}} [t]).
\end{align*}
We need to check that 
\begin{equation}\label{diag:2}
\xymatrix{ ([s],(P_1,\dots.P_p)) \ar[d]_{([s],(P_1,\dots,P_p),[s]\xrightarrow{\la_E\la^{-1}} [s])} \ar[rrr]^-{\phi=(G,(Q_1',\dots,Q'_q),G\stackrel{\kappa}\to [t])} 
&&& ([t],(Q_1,\dots,Q_q)) \ar[d]^{([t],(Q_1,\dots,Q_q),[t]\xrightarrow{\la_F\mu^{-1}} [t])}\\
([s],(\la_E\lambda^{-1}P_1,\dots,\la_E\la^{-1}P_p)) \ar[rrr]^-{\Lambda\Pi(\phi)} &&& ([t],(\la_F\mu^{-1}Q_1,\dots,\la_F\mu^{-1}Q_q)}
\end{equation}
commutes in~$\Lev_A$. 
Now because all the morphisms in the diagram define morphisms in the fiber $\rmC_A[r]\minus J$, we get a diagram 
$$\xymatrix@R-1pc{\rmC_A[s]\ar[dd]_{([s],\la_E\la^{-1})} \ar[rr]^-{(G,\kappa)} & & \rmC_A[t] \ar[dd]^{([t],\la_F\mu^{-1})}\\
& \rmC_A[r]  \ar[ul] \ar[ur] \ar[dl] \ar[dr] & \\
\rmC_A[s] \ar[rr]_-{(\hat\la_EF,\la_F\hat\la_E^{-1})} && \rmC_A[t]}$$
of commuting triangles in $\Exp_A$, from which it follows that the outer square commutes  as $\Exp_A$ is a subcategory of the groupoid $\bCan_A^\x$. 
From this, we get that 
~$G=\hat\la(F)$ and   $\la_F\mu^{-1}\kappa=\la_F\hat\la^{-1}$. 
Now going along the top of diagram (\ref{diag:2}), the maps compose to $(G,(Q_1',\dots,Q_q'),G\xrightarrow{\la_F\mu^{-1}\kappa}[t])$ while going along the bottom they compose to 
$((\hat\la F,(\la\mu^{-1}Q_1,\dots,\la\mu^{-1}Q_q),\hat\la F\xrightarrow{\la_F\hat\la^{-1}} [t])$. And these are equal by the previous computation, using also that $Q_i=\kappa(Q_i')$ for each $i$.  
This completes the proof.
\end{proof}

%%%

\subsection{Thomason's homotopy colimit}\label{def:Sigma_A}

We will now recall Thomason's explicit model for the homotopy colimit~\eqref{eq:defTho}. Let
\[
\Si_A\colon\bSet^\x \to \bSet^\x
\]
be the functor defined on objects by~$\Si_A[r]=[rn]$ and on morphisms~$\la\colon[r]\to [r]$ by the diagram 
\begin{equation}\label{diag:SiA}
\xymatrix@R-1pc{[rn] \ar@{-->}[rr]^-{\Si_A(\la)} && [rn]\\
[r]\x A=E \ar[rr]^-{\hat\la=\la\x A} \ar@{}[d]|{\IV} \ar[u]^{\la_E}&& E=[r]\x A \ar@{}[d]|{\IV} \ar[u]_{\la_E}\\
[r] \ar[rr]^-{\la} && [r], }
\end{equation}
that is setting~$\Si_A(\la)=\la_E\circ \hat\la\circ \la_E^{-1}$, where~$\la_E$ is the map given by lexicographic ordering of Lemma~\ref{lem:lexi}. 
Explicitly~$\la_E$ orders the elements of~$[r]\x A$ as 
$$(1,a_1),\dots,(1,a_n),\dots, (r,a_1),\dots,(r,a_n),$$
so~$\Si_A(\la)$ is the well-known ``block'' version of~$\la$, permuting blocks of size~$n$. 
Using this block permutation description, one sees that~$\Si_A$ is a strict symmetric monoidal functor. 

Recall that 
$$\Tho_A=\hocolim_\bbN(\bSet^\x,\Si_A)$$
is obtained as a homotopy colimit in the category of non-unital permutative categories for~$\bbN$ here considered as the category with one object~$*$ associated to the monoid~$(\bbN,+)$ and 
$$(\bSet^\x,\Si_A)\colon \bbN\rar \mathbf{Perm}$$ 
denoting the functor from~$\bbN$ to the category of non-unital 
permutative categories taking the unique object to~$\bSet^\x$ and a morphism~$k\in \bbN$ to~$\Si_A^k$, the~$k$--th iterate of the functor~$\Si_A$.  
Here we use the explicit model given by Thomason in~\cite[Construction 3.22]{Thomason} for this homotopy colimit.
The objects of~$\Tho_A$ are the tuples~$([r_1],\dots,[r_m])$ with~$m>0$ and each~$r_i>0$.
Morphisms are given as tuples
\[
(\psi,\mu,f_1,\dots,f_{n})\colon ([r_1],\dots,[r_{m}])\rar ([s_1],\dots,[s_{n}])
\] 
where~$\psi\colon [m]\to [n]$ is a surjection,~$\mu\colon [m]\to\bbN$ is a function, and for each $j=1,\dots,n$,
\[
f_j\colon [\sum_{i\in\psi^{-1}(j)} r_in^{\mu(i)}]\longrightarrow [s_j]
\]
is a bijection. 
%where~$[m]=\{1,\dots,m\}$ as before and
% where we have already used the notation
%for
%\[
%E(\psi,\mu)_j=\bigsqcup_{i\in\psi^{-1}(j)} [r_i]\times A^{\mu(i)}.
%\] 
%(Here of course we will soon identify the source of~$f_j$ with an expansion of the set~$[\sum_{i\in \psi^{-1}(j)}r_i]$.) 
The composition of such morphisms is defined by composition of the surjections, corresponding~``multi-additions'' of the values~$\mu(i)$ and composition of the bijections~$f_j$, appropriately multiplied with powers of~$A$: if~$(\psi,\mu,f_1,\dots,f_{n})$ is as above and~$(\psi',\mu',f'_1,\dots,f'_{p})\colon ([s_1],\dots,[s_{n}])\rar ([t_1],\dots,[t_{p}])$, then the composition~$(\psi',\mu',f'_1,\dots,f'_{p})\circ (\psi,\mu,f_1,\dots,f_{n})$ 
is given by
\[
(\psi'\circ \psi,
\mu+(\mu'\circ\psi),
f'_1\circ\left(\bigoplus_{j\in \psi'^{-1}(1)}\Si_A^{\mu'(j)}f_j\right)\circ\s_1,
\dots,
f'_p\circ\left(\bigoplus_{j\in \psi'^{-1}(p)}\Si_A^{\mu'(j)}f_j\right)\circ \s_p
)
\]
for~$\s_k$ the block permutation induced by the permutation reordering~$(\psi'\psi)^{-1}(k)$ as~$\psi^{-1}(j_1),\dots,\psi^{-1}(j_q)$ for~$\{j_1<\dots<j_q\}=(\psi')^{-1}(k)$. 
The monoidal structure is given on objects by juxtaposition: 
\[
([r_1],\dots,[r_{m}])\oplus ([s_1],\dots,[s_{n}])=([r_1],\dots,[r_{m}],[s_1],\dots,[s_{n}])
\]
%The unit is the empty sequence~$()$ 
and likewise on morphisms, 
and the symmetry is the map~$(\s_{m,n},{\bf 0},\id,\dots,\id)$ is induced by the symmetry in~$\bSet^\x$, with~${\bf 0}\colon[m+n]\to \bbN$ denoting the zero map.

The goal of this section is to prove the following

\begin{theorem}\label{thm:LevTh}
There is a functor
$$H\colon\Lev_A\rar \Tho_A$$
which induces a homotopy equivalence on classifying spaces. 
\end{theorem}

\begin{proof}[Part 1 of the proof of Theorem~\ref{thm:LevTh}: Definition of the functor.]
The functor~$H\colon\Lev_A\rar \Tho_A$ is defined on objects~$([r],(P_1,\dots,P_p))$ of~$\Lev_A$ with~$r_i=|P_i|$  by setting 
\[
H([r],(P_1,\dots,P_p))=\left\{\begin{array}{ll}([r_1],\dots,[r_p],[r-\sum_ir_i]) & \textrm{if}\ \sum_i r_i\neq r\\
([r_1],\dots,[r_p]) & \textrm{if}\ \sum_i r_i= r.
\end{array}\right.
\]
Consider now a morphism 
$$\phi\colon([r],(P_1,\dots,P_p)) \ge (E,(Q'_1,\dots,Q'_q))\ \stackrel{\la}\rar\ ([s],(Q_1,\dots,Q_q))$$ in~$\Lev_A$.
% is given by a tuple~$\phi=(E,(Q'_1,\dots,Q'_q),E\stackrel{\ \la}\rar [s])$ with 
%$([r],(P_1,\dots,P_p))\ge (E,(Q'_1,\dots,Q'_q))$ in~$\calL[r]$ and~$\la$ an isomorphism  between~$(E,(Q'_1,\dots,Q'_q))$ and~$([s],(Q_1,\dots,Q_q))$. 
%To define~$H$ on morphisms, we first associate to the poset morphisms of~$\calL[r]$ objects of the same form as the morphisms in the category~$\Tho_A$: 
The relation~$([r],(P_1,\dots,P_p))\ge (E,(Q'_1,\dots,Q'_q))$ in~$\calL[r]$ gives that 
there exists~$0\le p_0,p_1\le p$ and~$\kappa_1,\dots,\kappa_q\ge 1$ such that~$p_0+\sum_j \kappa_j+p_1=p$, with~$E$ the level expansion of~$[r]$ along~$P_1\cup \dots\cup P_{p_0}$ and
\[
Q'_j=P_{p_0+\kappa(1)+\dots+\kappa(j-1)+1}\sqcup\dots\sqcup P_{p_0+\kappa(1)+\dots+\kappa(j)}.
\]
The object~$H([r],(P_1,\dots,P_p))$ will be a sequence of~$p'=p$ or~$p+1$ elements depending on whether~\hbox{$\cup_\ell P_\ell=[r]$} or not. Likewise, the image~$H([s],(Q_1,\dots,Q_q))$ will be a sequence of~$q'=q$ or~$q+1$ elements. 
We define~$H(\phi)=(\psi,\mu,f_1,\dots,f_{q'})$ for~$\psi,\mu$, and the functions~$f_j$ are defined as follows. Let 
\[
\psi(i)=
\begin{cases}
j 	& \textrm{if }P_i\subset Q'_j \\
q+1 & \textrm{if }P_i\not\subset\cup_j Q'_j\textrm{ or }i=p+1.
\end{cases}
\]
This gives a well-defined surjective map~$\psi\colon[p']\to [q']$: any~$j\in [q]$ with~$j\le q$ is in the image of~$\psi$ because any~$Q_j$ is a union of~$P_i$'s, and if~$q'=q+1$, then we must either have~$p'=p+1$ and~$\psi(p+1)=q+1$ or~$p'=p$ with~$p_0\ge 1$ or~$p_1\ge 1$, in which case either~$\psi(1)=q+1$ or~$\psi(p)=q+1$.    
%$$\psi(l)=\left\{\begin{array}{ll} q+1& 1\le l\le i\\
%j & i+\kappa_1+\dots+\kappa_{j-1}+1\le l\le i+\kappa_1+\dots+\kappa_{j}\\
%q+1 & p-k+1\le l\le p+1
%\end{array}\right.$$
We define~$\mu\colon[p']\to \bbN$ by 
$$\mu(i)=\left\{\begin{array}{ll} 1 & 1\le i\le p_0\\
0 & \textrm{else}. 
\end{array}\right.$$
Finally, for~$1\le j\le q$, let~$f_j\colon[\sum_{i\in \psi^{-1}(j)}r_i]\to [\sum_{i\in \psi^{-1}(j)}r_i]=[s_j]$ for $s_j=|Q_j|$, be the bijection making the following diagram commute: 
\begin{equation}\label{equ:Hmor}
\xymatrix@R-2pc{&[\sum_{i\in \psi^{-1}(j)}r_i]  \ar[rr]^-{f_j} && [\sum_{i\in \psi^{-1}(j)}r_i]=[s_j]\\
\bigsqcup_{i\in \psi^{-1}(j)}[r_i]\ar[ur]^-\oplus & &&  \\
&\bigsqcup_{i\in \psi^{-1}(j)}P_i  \ar[ul]^-{\sqcup_i\la_{P_i}} & \ar[l]\bigcup_{i\in \psi^{-1}(j)}P_i = Q'_j \ar[r]^-{\la|_{Q'_j}} & Q_j \ar[uu]_-{\la_{Q_j}} 
}\end{equation}
for~$\la_{P_i}$ and~$\la_{Q_j}$ the lexicographic maps of Lemma~\ref{lem:lexi}. (In this case, each $\mu(i)=0$.)  %and~$\hat\la$ the restriction of~$\la\colon E\to [s]$ to~$Q'_j$.  
If~$q'=q+1$, we set $f_{q+1}$ 
%\[
%f_{q+1}\colon[\sum_{i\in \psi^{-1}(q+1)}\!\!\! r_i n^{\mu(i)}]=[\sum_{1\le i\le p_0}\!\! r_i n + \sum_{p-p_1\le i\le p}\!\! r_i \,+\,(r-\!\!\sum_{1\le i\le p}\!\! r_i)]\to [\sum_{i\in \psi^{-1}(q+1)}\!\!\! r_i n^{\mu(i)}]=[s-\sum_{1\le j\le q}\!\! s_j]
%\] 
to be the bijection making the following diagram commute: 
$$\xymatrix{[\sum_{i\in \psi^{-1}(q+1)}r_i n^{\mu(i)}]  \ar[rr]^-{f_{q+1}} && [\sum_{i\in \psi^{-1}(q+1)}r_i n^{\mu(i)}]\\
\big(\bigsqcup_{1\le i\le p_0}P_i\x A \big)\,\sqcup\,\big(\bigsqcup_{p-p_1< i\le p}P_i\big)\,  \sqcup \,\big([r]\minus \cup_i P_i\big) \ar[u]^{\sqcup(\textrm{lexi})}%^{\sqcup_i\la_{P_i\x A^{\mu(i)}}} 
& \ar[l] E\minus (\cup_jQ'_j) \ar[r]^-{\la} & [s]\minus(\cup_jQ_j) \ar[u]_{(\textrm{lexi})}%_{\la_{Q_j}} 
}$$
with the vertical maps defined as in the previous case.

We check that~$H$ respects composition. Consider a composition in~$\Lev_A$: 
$$\xymatrix@R-1pc{([r],(P_1,\dots,P_p)) \ge   (E,(Q_1',\dots,Q_q')) %\ar@{}@<-6ex>[d]|{\IV} 
\ar[r]^-\la & ([s],(Q_1,\dots,Q_q)) \ge  (G,(R_1',\dots,R_u')) \ar[r]^-\kappa  %\ar@<6ex>[d]^\kappa 
&  ([t],(R_1,\dots,R_u))  %\\
%(\hat\la^{-1}G,(\hat\la^{-1}R_1',\dots,\hat\la^{-1}R_u')) \ \ \ \ \ \ \ \ar[r]^-{\kappa\circ\hat\la} & \ \ \ \  \ \ \ ([t],(R_1,\dots,R_u))  
}$$
with image under~$H$ the composition 
$$\xymatrix{([r_1],\dots,[r_{p'}]) \ar[r]^-{(\psi,\mu,{\bf f})} & ([s_1],\dots,[s_{q'}]) \ar[r]^-{(\psi',\mu',{\bf f'})} & ([t_1],\dots,[t_{u'}]).}$$
By definition of the composition in~$\Tho_A$, these maps compose to
\[
\big(\psi'\circ \psi,\mu+(\mu'\circ\psi),f'_1\circ (\oplus_{i\in \psi'^{-1}(1)}\Si_A^{\mu'(i)} f_i)\circ\s_1,\dots,f'_{u'}\circ (\oplus_{i\in \psi'^{-1}(u')}\Si_A^{\mu'(i)}f_i)\circ \s_{u'}\big).
\]
The fact that the value of the functor~$H$ on the composed morphism  
$$\xymatrix{([r],(P_1,\dots,P_p)) \ge (\hat\la^{-1}G,(\hat\la^{-1}R_1',\dots,\hat\la^{-1}R_u')) \ar[r]^-{\kappa\circ\hat\la} &  ([t],(R_1,\dots,R_u))  
}$$
agrees on the first component follows from the fact that~\hbox{$P_i\subset \hat\la^{-1}(R_k)$} if and only if~$P_i\subset Q'_j$ and~$\la(Q'_j)=Q_j\subset R_k$. 
On the second component, it is because~$\mu$ records the~$P_i$'s that are expanded by the first map and~$\mu'$ records those that are expanded by the second map. For the last component, assuming first  $1\le k\le u$,  with~$(\psi')^{-1}(k)=\{j_1<\dots<j_m\}$, we have $\mu+(\mu'\circ\psi)(i)=0$ for each $i\in (\psi'\psi)^{-1}(k)$, and the diagram 
$$\xymatrix{[\sum_{i\in (\psi'\psi)^{-1}(k)}r_i]\stackrel{\s_k}\cong  [\sum_{\ell=1}^m(\sum_{i\in (\psi)^{-1}(j_\ell)}r_i)] \ar[r]^-{\oplus_\ell f_\ell} & [\sum_{\ell=1}^m(\sum_{i\in (\psi)^{-1}(j_\ell)}r_i)] \ar[r]^-{f'_k} & [\sum_{i\in (\psi'\psi)^{-1}(k)}r_i]\\
\bigsqcup_{i\in (\psi'\psi)^{-1}(k)}P_i \stackrel{\s_k}\cong \bigsqcup_{\ell=1}^m(\bigsqcup_{i\in (\psi)^{-1}(j_\ell)}P_i) \lar \bigsqcup_{\ell=1}^mQ'_{j_\ell}  \ar[r]^-{\hat \la}
 \ar@<-6ex>[u]^{\sqcup(\textrm{lexi})} \ar@<13ex>[u]^{\sqcup(\textrm{lexi})}
 & \bigsqcup_{\ell=1}^mQ_{j_\ell} \ar[u]^{\sqcup(\textrm{lexi})} \lar R'_k\ar[r]^-\kappa & R_k \ar[u]^{(\textrm{lexi})}\\
& \bigcup_{i\in (\psi'\psi)^{-1}(k)}P_i=\hat\la^{-1}R'_k \ar[ul] \ar[ul]+<-17ex,-2ex> \ar[ul]+<20ex,-2ex> \ar[u]+<5ex,-2ex>^-{\hat \la}  \ar[ur]_{\kappa\circ\hat\la}& 
}$$
%with $\bar \mu:=\mu+(\mu'\circ\psi)$, and similarly for
commutes by definition of $H$ on morphisms for the squares, and with the first three diagonal maps at the bottom of the diagram  the canonical maps coming from the fact that the target is in each case a decomposition of the source into disjoint subsets of it. The top row is the map $f'_{u'}\circ (\oplus_{i\in \psi'^{-1}(u')}\Si_A^{\mu'(i)}f_i)\circ \s_{u'}$, but is also equal to corresponding component of the value of $H$ on the composition by the commutativity of the outer diagram. 
The case~$k=u+1$, when relevant, is similar.   
\end{proof}

The proof that the functor~$H$ induces a homotopy equivalence on classifying spaces will be analogous to that of the functor~$\Lev_A\to \Exp_A$: we will compare the homotopy fibers to  a poset of simplices of a semi-simplicial set analogous to~$\rmL(X)$. %Given a tuple of positive numbers~$(r_1,\dots,r_k)$, we will now define a semi-simplicial set~$\rmL(r_1,\dots,r_k)$. 
We start by introducing the relevant semi-simplicial set. 

Let~$\bbN=\{0,1,2,\dots\}$ denote the set of  natural numbers and  
$(\bbN^k,\le)$ the poset of~$k$--tuples of natural numbers, with the poset structure defined by~$(m_1,\dots,m_k)\le (n_1,\dots,n_k)$ if~$m_i\le n_i$ for each~$i$. 

\begin{definition}
For~$k> 0$, let~$\rmM(k)$ denote the semi-simplicial subset of the nerve of~$(\bbN^{k},\le)$ with 
%with 0--simplices the elements of~$\bbN^{k}$, and
the same vertices and where a~$p$--simplex is a sequence 
\[
v_0= (n_1^0,\dots,n^0_{k})<\dots <  v_p=(n^p_1,\dots,n^p_{k})
\]
in~$(\bbN^{k},\le)$ such that~$v_i -v_0\in \{0,1\}^k$ for all~$1\le i\le p$, i.e.~$n_j^i-n_j^0\in\{0,1\}$ for all~$j$.  
%If~$k=0$,~$\rmL()$ is a point, denoted~$()$. 
\end{definition}

Note that the condition~$v_i -v_0\in \{0,1\}^k$ for all~$1\le i\le p$ is equivalent to the condition~$v_i -v_j\in \{0,1\}^k$ for all~$0\le j\le  i\le p$ given that~$v_0<\dots<v_p$ in the poset~$(\bbN^k,\le)$. 

For any a tuple of positive numbers~$(r_1,\dots,r_k)$, we can identify $\rmM(k)$ with a semi-simplicial subset of~$\rmL[\sum r_i]$: to a vertex~$v=(n_1,\dots,n_k)$ of~$\rmM(k)$, we can associate the set  
$$[r_1]\x A^{n_1}\sqcup \dots \sqcup [r_k]\x A^{n_k}.$$
which we think of as an expansion of the set~$[r]=[r_1]\oplus \dots\oplus [r_k]$. Under this identification, the condition for a sequence~\hbox{$(n_1^0,\dots,n^0_k)<\dots <(n^p_1,\dots,n^p_k)$} to form a~$p$--simplex in~$\rmM(k)$ 
is precisely that the corresponding expansions of~$[r]$ form a~$p$--simplex in~$\rmL[r]$. But note that here we are only considering level expansions of a specific form, allowing only expansions along the subsets~$[r_i]\times A^{n_i}$, or unions of such. 

\begin{lem}\label{lem:Lk}
For any~$k> 0$, the semi-simplicial set~$\rmM(k)$ is contractible. 
\end{lem}

\begin{proof}
%It is enough consider the case where~$k\ge 1$ and each~$r_i\ge 1$. 
%This proof is very similar to the proof of the contractibility of~$\rmL(X)$. 
We consider the restriction of the rank filtration of the poset~$(\bbN^k,\le)$ to~$\rmM(k)$: the {\em rank} of a vertex~$(n_1,\dots,n_k)$ is~$\sum_in_i$. Let~$F_i\rmM(k)$ be the semi-simplicial subset on the vertices of rank at most~$i$.  We have 
$$\{(0,\dots,0)\}=F_0\rmM(k)\subset F_1\rmM(k)\subset \dots \subset F_i\rmM(k)\subset \dots\subset \rmM(k).$$
By compactness, any map from a sphere into~$\rmM(k)$ will land in a finite filtration, so it is enough to show that~$F_i\rmM(k)$ is contractible for each~$i$. Clearly~$F_0\rmM(k)$ is contractible, as it is just a point. Assuming that we have proved that~$F_i\rmM(k)$ is contractible, we will show that~$F_{i+1}\rmM(k)$ is also contractible. As vertices of any simplex in 
$\rmM(k)$ necessarily have distinct rank, we see that~$F_{i+1}\rmM(k)$ is obtained from~$F_i\rmM(k)$ by adding a cone on each vertex~$v$ of rank equal to~$i+1$, attached to 
$F_i\rmM(k)$ along~$\Link(v)\cap F_i\rmM(k)$. By definition, a simplex~$w_0<\dots<w_q$ of~$\rmM(k)$ lies in~$\Link(v)$ precisely if~$w_0<\dots<w_p<v<w_{p+1}<\dots<w_q$ is a simplex of~$\rmM(k)$. As~$v$ has rank~$i+1$, we see that only~$w_0,\dots,w_p$ have rank at most~$i$, and hence that 
$\Link(v)\cap F_i\rmM(k)$ identifies with the semi-simplical subset~$\rmM(k)_{<v}\subset \rmM(k)$ of simplices~$w_0<\dots<w_p$ such that~\hbox{$w_0<\dots<w_p<v$} is also a simplex of~$\rmM(k)$. 
% Indeed, any such~$w$ necessarily has rank at most~$i$ and any simplex~$w_0<\dots<w_p$ in~$\rmM(k)_{<v}$  lies in the link of~$v$. And conversally, if 
%$w_0<\dots<w_p<v<w_{p+1}<\dots<w_q$ is a simplex of~$\rmM(k)$ with the~$w_i$'s in~$F_i\rmM(k)$, then we must have~$p=q$ and~$w_0<\dots<w_p$ a simplex of~$\rmM(k)_{<v}$.

So we are left to show that each~$\rmM(k)_{<v}$ is contractible. Note first that the semi-simplicial set~$\rmM(k)_{<v}$ is the~(small) nerve of a poset: for vertices~$w,w'$ of~$\rmM(k)_{<v}$, set~\hbox{$w\prec w'$} if~$w<w'$ is a~$1$--simplex in~$\rmM(k)_{<v}$. This defines a poset structure on~$\rmM(k)_{<v}$  as~$w\prec w'\prec w''$ implies that~$w<w''$ is also a~$1$--simplex in $\rmM(k)_{<v}$, because~$w<v$ and~$w''<v$ are simplices in~$\rmM(k)$, which means that~\hbox{$v-w\in \{0,1\}^k$} and~\hbox{$0\le w''_i-w_i\le v_i-w_i\le 1$} is also necessarily in~$\{0,1\}$ for each~$i$. 
%Contractibility then follows from the fact that this is the nerve of a poset with a least element: 
This allows us to finish the proof, because that poset has a least element: indeed, if~$v=(n_1,\dots,n_k)$, then~\hbox{$w=(n_1,\dots,n_k)-(\epsilon_1,\dots,\epsilon_k)$} 
with $\epsilon_i=\min\{1,n_i\}$ 
%\[
%\epsilon_i =\left\{ \begin{array}{ll}1 & \textrm{if}\ n_i\ge 1\\
% 0 & \textrm{if}\ n_i=0.
%\end{array}\right. 
%\]
is the least element.
\end{proof}

As in the case of~$\rmL[r]$, we will need to pass to the poset of simplices of~$\rmM(k)$. We start by giving a description of~$\rmM(k)$ that is more suitable for our needs. The idea is, just as in the case of $\rmL(X)$,  that we can encode the information of a simplex
\[
(n_1^0,\dots,n^0_{k})<\dots <  (n^p_1,\dots,n^p_{k})
\]
by remembering the first tuple of elements~$ (n_1^0,\dots,n^0_{k})$, and then remembering for each~$1\le i\le p$ which set of indices~$J_i\subset \{1,\dots,k\}$ was raised by 1 when going from 
$ (n_1^{i-1},\dots,n^{i-1}_{k})$ to~$ (n_1^i,\dots,n^i_{k})$, that is for which index~$\ell$ the difference~$n_\ell^{i}-n_\ell^{i-1}$ is equal to 1. This idea formalizes to the following:  

\begin{lem}\label{lem:LL'k}
%Let~$\calL(r_1,\dots,r_k)$ denote the poset of simplices of~$\rmM(k)$. This poset is isomorphic to the poset 
The semi-simplicial set~$\rmM(k)$ is isomorphic to the semi-simplicial set 
$\rmM'(k)$ whose~$p$--simplices are the tuples 
$$((n_1,\dots,n_k),J_1,\dots,J_p)$$
for~$(n_1,\dots,n_k)\in \bbN^k$ and~$J_1,\dots,J_p$ a collection of~$p\ge 0$ disjoint non-empty subsets of~$[k]=\{1,\dots,k\}$, 
%with poset structure defined by 
%$$((m_1,\dots,m_k),I_1,\dots,I_q)\le ((n_1,\dots,n_k),J_1,\dots,J_p)$$
%if one can get~$((m_1,\dots,m_k),I_1,\dots,I_q)$ from~$ ((n_1,\dots,n_k),J_1,\dots,J_p)$ by a composition of the following operations: 
and with the face maps defined by 
%\begin{enumerate}
%\item~$d_0 ((n_1,\dots,n_k),J_1,\dots,J_p)=((n_1+\epsilon_1,\dots,n_k+\epsilon_k),J_2,\dots,J_p)$ for~$\epsilon_j=1$ if~$j\in J_1$ and~$0$ otherwise;
%\item~$d_i ((n_1,\dots,n_k),J_1,\dots,J_p)=((n_1,\dots,n_k),J_1,\dots,J_i\cup J_{i+1},\dots ,J_p)$ for~$1\le i\le p$; 
%\item~$d_p ((n_1,\dots,n_k),J_1,\dots,J_p)=((n_1,\dots,n_k),J_1,\dots,J_{p-1})$.
%\end{enumerate}
\[
d_i ((n_1,\dots,n_k),J_1,\dots,J_p)=
\begin{cases}
((n_1+\epsilon_1,\dots,n_k+\epsilon_k),J_2,\dots,J_p) & i=0 \\
((n_1,\dots,n_k),J_1,\dots,J_i\cup J_{i+1},\dots ,J_p) & 0<i<p\\ 
((n_1,\dots,n_k),J_1,\dots,J_{p-1}) & i=p,
\end{cases}
\]
where $\epsilon_j=1$ if $j\in J_1$ and $0$ otherwise.
%there exists~$0\le i,k\le p$ and~$\kappa_1,\dots,\kappa_q$ such that 
%\begin{enumerate}
%\item~$n_j=m_j+\epsilon_j$ with~$\epsilon_j=1$ if~$j\in\cup_{1\le l\le i}I_i$ and~$\epsilon_j=0$ otherwise, 
%\item~$J_j=I_{i+\kappa_1+\dots \kappa_{j-1}+1}\cup\dots\cup I_{i+\kappa_1+\dots \kappa_j}$
%\item~$i+\kappa_1+\dots+\kappa_q=p-k$. 
%\end{enumerate}
\end{lem}

\begin{proof}
Define a map~$F\colon\rmM(k)\to \rmM'(k)$ on~$p$--simplices by setting 
$$F\big((n_1^0,\dots,n^0_{k})<\dots <  (n^p_1,\dots,n^p_{k})\big)=\big((n_1^0,\dots,n^0_{k}),J_1,\dots,J_p\big)$$
for~$J_i\subset \{1,\dots,k\}$ the set of indices of~$(n_1^{i},\dots,n^{i}_{k})- (n_1^{i-1},\dots,n^{i-1}_{k})$ that are equal to 1. By definition of~$\rmM(k)$, this difference lies in~$\{0,1\}^k\minus \{0\}^k$ so~$J_i$ is non-empty.
Also the subsets~$J_i$ must be pairwise disjoint since the existence of $\ell\in J_i\cap J_j$ for $i<j$ would contradict that $n_\ell^j-n_\ell^0\in\{0,1\}$. 
%as each difference~$(n_1^{j},\dots,n^{j}_{k})- (n_1^{i},\dots,n^{i}_{k})$ for any~$j>i$ likewise   lies in~$\{0,1\}^k\minus\{0\}^k$. 
Hence~$F$ has image inside~$\rmM'(k)$. The face maps on the latter semi-simplicial set are defined so that~$F$ is semi-simplicial. An inverse is obtained by taking 
a tuple~$\big(v=(n_1,\dots,n_{k}),J_1,\dots,J_p\big)$ to the simplex~$v<v+\langle J_1\rangle<\dots<v+\langle J_p\rangle$ with~$\langle J_i\rangle=(\epsilon_1,\dots,\epsilon_k)$ for~$\epsilon_j=1$ if~$\epsilon_j\in J_i$ and~$0$ otherwise. 
%As these cover the requirements for being a simplex of~$\rmL(r_1,\dots,r_k)$, we see that~$F$ is well-defined and bijective as a map of sets. We are left to check that it respects the ordering. We have defined the ordering on~$\calL'(r_1,\dots,r_k)$ so that this holds, each operation~$d_i$ corresponding to the~$i$--th face map in~$\rmL(r_1,\dots,r_k)$.  \nnote{enough?} 
\end{proof}

\begin{proof}[Part 2 of the proof of Theorem~\ref{thm:LevTh}: The functor~$H\colon\Lev_A\to \Tho_A$ is an equivalence]
We will again consider the homotopy fibers of the functor~$H$. We will show that they are homotopy equivalent to  
the poset~$\calM(k)$ of simplices of~$\rmM'(k)$, 
from which the result will follow as a consequence of Lemma~\ref{lem:Lk} and Lemma~\ref{lem:LL'k}. 
For any object~$([r_1],\dots,[r_k])$ in~$\Tho_A$, we will define functors
$$\Lambda\colon\calM(k)^{\op} \ \longleftrightarrow \ ([r_1],\dots,[r_k])\minus H\colon\Pi$$
and show that they define an equivalence. 
We start by defining~$\Lambda$. 

\noindent
{\em Definition of~$\Lambda$.} To a~$p$--simplex~$((n_1,\dots,n_k),J_1,\dots,J_p)$ of the semi-simplicial set~$\rmM'(k)$ we 
can associate the set~\hbox{$E=[r_1]\x A^{n_1}\sqcup \dots\sqcup [r_k]\x A^{n_k}$}, considered as an expansion of $[\sum_ir_i]$ via the sum $[r_1]\sqcup\dots\sqcup [r_k]\xrightarrow{\oplus} [\sum_ir_i]$, and  canonically identified with~$[\sum_i r_in^{n_i}]$ via~$\la_E$. 
%with its corresponding lexicographic ordering bijection~$\la_E:E\to [\sum r_in_i]$, 
For $1\le j\le p$, let 
\[
P_j=\la_E(\bigsqcup_{i\in J_j}[r_i]\x A^{n_i})\subset [\sum_{1\le i\le k} r_in^{n_i}]
\]
be the subset of $ [\sum_ir_in^{n_i}]$ corresponding to~$J_j$. 
We define~$\Lambda$ on objects by setting
$$\Lambda((n_1,\dots,n_k),J_1,\dots,J_p)= (([\sum_{1\le i\le k} r_in^{n_i}],P_1,\dots,P_p),(\psi,\mu,\id))$$
with
\[
(\psi,\mu,\bid)\colon([r_1],\dots,[r_k])\rar H([\sum_{1\le i\le k} r_in^{n_i}],(P_1,\dots,P_p))=([\sum_{i\in J_1}r_in^{n_i}],\dots,[\sum_{i\in J_p}r_in^{n_i}],([\sum_{i\notin\cup_j J_j}r_in^{n_i}]))
\]
defined by  
setting~$\psi(i)=j$ if~$i\in J_j$, and~$\psi(i)=p+1$ if~$i\notin\cup_jJ_j$, and~$\mu(i)=n_i$, and with~$\bid=(\id,\dots,\id)$ denoting a $p'$--tuple of identity maps; 
here, and in the rest of the proof, the last component of the rightmost term above is dropped if it is an empty set. And as before, we use $p'$ to denote $p$ or $p+1$ accordingly.

To define~$\Lambda$ on morphisms, note first that 
\[
X=((m_1,\dots,m_k),I_1,\dots,I_p)\ge ((n_1,\dots,n_k),J_1,\dots,J_q)=Y
\]
in~$\calM(k)$ if there are~$0\le p_0,p_1\le p$ and~$\kappa_1,\dots,\kappa_q\ge 1$ with 
$p_0+\sum \kappa_j+p_1=p$ such that  
$n_i=m_i+\epsilon_i$ for~$\epsilon_i=1$ if~$i\in I_1\cup\dots\cup I_{p_0}$ and~$0$ otherwise, and~$J_j=I_{p_0+\kappa_1+\dots+\kappa_{j-1}+1}\cup\dots\cup I_{p_0+\kappa_1+\dots+\kappa_{j}}$. 
To such an inequality, we associate the morphism 
$$\alpha\colon\Lambda(X)=([\sum_{1\le i\le k}r_in^{m_i}],(P_1,\dots,P_p)) \ge (E,(Q_1',\dots,Q_q')) \stackrel{\la_E}\rar ([\sum_{1\le i\le k}r_in^{n_i}],(Q_1,\dots,Q_q))=\Lambda(Y)$$
 in~$\Lev_A$ with~$E$ the level expansion of~$[\sum_ir_in^{m_i}]$ along the subset~$P_1\cup\dots\cup P_{p_0}$, and~$Q'_j=\cup_{i\in J_j}P_i$. 
For this to make sense, we need that $\la_E(Q'_j)=Q_j$, which holds by the compatibility of the lexicographic orderings~(Lemma~\ref{lem:lexi}). 
Analyzing diagram~\eqref{equ:Hmor} in the definition of $H$, we see that~\hbox{$H(\al)=(\phi,\nu,\s_1,\dots,\s_{q'})$} with 
$$\s_j\colon \sum_{i\in \psi^{-1}(j)}\sum_{\ell\in I_i}r_\ell n^{m_\ell} \rar  \sum_{i\in J_j}r_\ell n^{m_\ell}$$
if $j\le q$ reordering the summands (corresponding to the natural map $\sqcup_{i\in \psi^{-1}(j)}P_i\to Q'_j$ as subsets of $[\sum_ir_in^{m_i}]$), and similarly for $j=q+1$ if relevant. And the diagram 
$$\xymatrix{([r_1],\dots,[r_k])\ar[rr]^-{(\psi,\mu,\bid)} \ar[drr]_-{(\psi',\mu',\bid)} &&  ([\sum_{i\in I_1}r_in^{m_i}],\dots,[\sum_{i\in I_p}r_in^{m_i}],([\sum_{i\notin\cup_{j=1}^p I_j}r_in^{m_i}])) \ar[d]^{H(\alpha)=(\phi,\nu,{\bf\s})} \\
&&  ([\sum_{i\in J_1}r_in^{n_i}],\dots,[\sum_{i\in J_q}r_in^{n_i}],([\sum_{i\notin\cup_{j=1}^q I_j}r_in^{n_i}]))}$$
commutes, as~$\phi\circ\psi=\psi'$, both surjections keeping track of which component each expanded~$[r_i]$ goes to,  and~$\mu'=\mu+(\nu\circ\psi)$ as~$\nu$ records which additional expansions occurred from~$\mu$ to~$\mu'$, and because the maps $\s_j$ precisely take care of the necessary reordering of the summands.  
Hence we have defined a morphism in the fiber from~$\Lambda(X)$ to~$\Lambda(Y)$. Left is to check that this assignment is compatible with composition, which is a direct computation.

\noindent
{\em Definition of~$\Pi$.} We define~$\Pi$ on objects of the fiber as follows. Let~$X=\big(([s],(P_1,\dots,P_p)),(\psi,\mu,{\bf f})\big)$ be an object of the fiber, with 
$$([r_1],\dots,[r_k]) \ \stackrel{(\psi,\mu,{\bf f})}\rar\ H([s],(P_1,\dots,P_p))= ([s_1],\dots,[s_p],[s-\sum_js_j])= ([s_1],\dots,[s_{p'}])$$
for~${\bf f}=(f_1,\dots,f_{p'})$ with~$f_j\colon\sqcup_{i\in \psi^{-1}(j)}[r_in^{\mu(i)}] \ \stackrel{\cong}\rar \ [s_j]$, where we set $[s_{p+1}]:=[s-\sum_js_j]$ if the latter set is non-empty. 
We define 
$$\Pi(X)=((\mu(1),\dots,\mu(k)),\psi^{-1}(1),\dots,\psi^{-1}(p)).$$ 
This is well-defined as the sets $[s_1],\dots,[s_p]$ are non-empty. 
Given a morphism~$X\to Y$ in the fiber, we need to check that we have the relation~$\Pi(X)\ge \Pi(Y)$ in $\calM(k)$. If~\hbox{$X=\big((\psi,\mu,{\bf f}),([s],(P_1,\dots,P_p))\big)$} and 
$Y=\big((\psi',\mu',{\bf f'}),([s'],(Q_1,\dots,Q_q))\big)$, such a morphism is given by a morphism in~$\Lev_A$
$$\al\colon([s],(P_1,\dots,P_p))\ge (E,(Q_1',\dots,Q_q')) \ \stackrel{\la}\rar\ ([s'],(Q_1,\dots,Q_q))$$
with~$E$ the level expansion of~$[s]$ along~$P_1\cup\dots\cup P_{p_0}$ and each~$Q_j'$ a union of~$P_i'$, and~$\la\colon E\to [s']$ a bijection, such that the diagram 
\begin{equation}\label{equ:morfib}
\xymatrix{([r_1],\dots,[r_k])\ar[rr]^-{(\psi,\mu,{\bf f})} \ar[drr]_-{(\psi',\mu',{\bf f'})} && ([s_1],\dots,[s_{p'}]) \ar[d]^{H(\al)}\\
&& ([s'_1],\dots,[s'_{q'}])
}\end{equation}
commutes. 
Because~$\al$ is a level expansion along the first~$p_0\ge 0$ subsets~$P_i$, we get that~$\mu'=\mu+\epsilon$ for~\hbox{$\epsilon_i=1$} if~\hbox{$i\in \psi^{-1}\{1,\dots,p_0\}$} and $\epsilon_i=0$ otherwise. Also, each~$\psi'^{-1}(j)$ is the union of the~$\psi^{-1}(i)$ for those $i$ for which~$Q'_j$ is the union of the~$P_i$'s. This precisely gives that 
$$\Pi(X)=((\mu(1),\dots,\mu(k)),\psi^{-1}(1),\dots,\psi^{-1}(p))\ge ((\mu'(1),\dots,\mu'(k)),\psi'^{-1}(1),\dots,\psi'^{-1}(q))=\Pi(Y)$$ in~$\calM(k)$, as required. 

\noindent
{\em The functors define an equivalence.} The composition~$\Pi\Lambda$ is the identity: As $\calM(k)^{\op}$ is a poset, it is enough to check this on objects, where 
\[
\Pi\Lambda ((n_1,\dots,n_k),J_1,\dots,J_q)=((\mu(1),\dots,\mu(k)),\psi^{-1}(1),\dots,\psi^{-1}(p))
\]
for the function~$\mu$ defined by setting~$\mu(i)=n_i$ and~$\psi$ defined by setting~$\psi(i)=j$ if~$i\in J_j$. On the other hand, the composition~$\Lambda\Pi$ is not the identity, but we will see that, just like in the previous cases,  the information forgotten by that composition can be used to define a natural transformation between~$\Lambda\Pi$ and the identity functor. Indeed, given an object~\hbox{$X=\big(([s],(P_1,\dots,P_p)),(\psi,\mu,{\bf f})\big)$} in the fiber, we have~\hbox{$\Lambda\Pi (X)=\big(([s],(\bar P_1,\dots,\bar P_p)),(\psi,\mu,\bid)\big)$} with~$s=\sum_ir_in^{\mu(i)}$ and~$\bar P_j$  the subset of~$[s]$ corresponding to 
the inclusion~\hbox{$[\sum_{i\in \psi^{-1}(j)}r_in^{\mu(i)}]\subset [\sum_ir_in^\mu(i)]=[s]$} coming from the inclusion of indexing sets $\psi^{-1}(j)\subset [k]$. 
Consider the following bijection
\begin{align*}
\underline f\colon [s]\rar P_1\sqcup\dots\sqcup P_p\sqcup [s\minus \cup_iP_i] &\xrightarrow{\sqcup_j\la_{P_j}}  [s_1]\sqcup\dots\sqcup [s_{p'}] \\
&\xrightarrow{\sqcup_jf_{j}}  [s_1]\sqcup\dots\sqcup [s_{p'}] \\
%&= [\textstyle\sum_{i\in \psi^{-1}(1)}r_in^{\mu(i)}]\sqcup\dots\sqcup [\textstyle\sum_{i\in \psi^{-1}(p')}r_in^{\mu(i)}]
& \xrightarrow{\sqcup_j\la_{\bar P_j}^{-1}} \bar P_1\sqcup\dots\sqcup \bar P_p\sqcup [s\minus \cup_i\bar P_i] \rar  [\textstyle\sum_{i}r_in^{\mu(i)}]=[s]
\end{align*}
where the first map splits $[s]$ into components according to the $P_j$'s and the last map  reassembles the components using the inclusion of the sets $\bar P_i$ instead. 
As $\underline f(P_j)=\bar P_j$, it defines a morphism $([s],(P_1,\dots,P_p),\underline f)$ in~$\Lev_A$ from $([s],(P_1,\dots,P_p))$ to~\hbox{$([s],(\bar P_1,\dots,\bar P_p))$}. 
Moreover, applying~$H$ to this morphism gives a commutative diagram in~$\Tho_A$
\[
\xymatrix{([r_1],\dots,[r_k])\ar[rr]^-{(\psi,\mu,\id)} \ar[drr]_-{(\psi,\mu,f_1,\dots,f_{p'})} &&([\sum_{i\in \psi^{-1}(1)}r_in^{n_i}],\dots,[\sum_{i\in \psi^{-1}(p')}r_in^{n_i}]) %([\sum_{i\notin\cup_j J_j}r_in^{n_i}])) %\ar@{=}[r] 
\ar[d]^{H([s],(P_1,\dots,P_p),\underline f)=(\id,{\bf 0},f_1,\dots,f_{p'})}\\ %& ([s_1],\dots,[s_p],[s-\sum_js_j]) \\
&& ([s_1],\dots,[s_p],[s-\sum_js_j]). %& 
}
\]
%and the maps~$f$ assemble to the map~$$[\sum_i r_i\mu(i)\ \stackrel{\la_E^{-1}}\rar \ \sqcup_{i}[r_i]\x A^{\mu(i)} \ \stackrel{\underline f}\rar \ [s]$$
%for~$\underline f$ the map of (\ref{equ:barf}), which takes each~$\tilde P_j$ to~$P_j$, 
Hence we have defined a morphism~$\eta$ in the fiber from~$\Lambda\Pi(X)$ back to~$X$. 
We are left to check that these morphisms together form a natural transformation. 
Given~$X=\big(([s],(P_1,\dots,P_p)),(\psi,\mu,{\bf f})\big)$ and~\hbox{$Y=\big(([s'],(Q_1,\dots,Q_q)),(\psi',\mu',{\bf f'})\big)$}, a  morphism~$X\to Y$ in the fiber is defined by a map 
$$\al\colon([s],(P_1,\dots,P_p))\ge (E,(Q_1',\dots,Q_q')) \ \stackrel{\la}\rar\ ([s'],(Q_1,\dots,Q_q))$$
in~$\Lev_A$ making the diagram (\ref{equ:morfib}) commute. We need to check that 
\begin{equation}\label{diag:3}
\xymatrix{ ([s],(P_1,\dots,P_p))\ge (E,(Q_1',\dots,Q_q'))\ar[r]^-\la\ar@<-7ex>[d]_{\underline f} & ([s'],(Q_1,\dots,Q_q))\ar[d]^{\underline f'}\\
 ([s],(P_1,\dots,P_p)) \ge (\bar E,(\bar Q_1',\dots,\bar Q_q'))\ar[r]^-{\la_{\bar E}} & ([s'],(Q_1,\dots,Q_q))
}\end{equation}
commutes in~$\Lev_A$, where we have suppressed the trivial inequalities in the vertical morphisms and where the bottom row is $\Lambda\Pi(\al)$. 
Now going along the top of the diagram, the morphisms compose to give $(E,(Q_1',\dots,Q_q'),\underline f'\circ \la)$ whereas going along the bottom, they compose to the morphism
$(\underline{\hat f}^{-1}\bar E,(\underline{\hat f}^{-1}\bar Q_1',\dots,\underline{\hat f}^{-1}\bar Q_q'),\la_E\circ \underline f)$. So we are left to check that $\underline f'\circ \la=\la_E\circ \underline f$ and that $\underline{\hat f}$ takes $(E,(Q_1',\dots,Q_q'))$ to 
$(\bar E,(\bar Q_1',\dots,\bar Q_q'))$. The latter fact follows from the fact that $\underline f$ takes $([s],(P_1,\dots,P_p))$ to  $([s],(\bar P_1,\dots,\bar P_p))$ and that the fact that $E$ and the $Q_j'$'s are build from $[s]$ and the $P_i$'s in exactly the same way as $\bar E$ from $[s]$ and the $\bar P_i$'s. Finally, the fact that the morphisms in (\ref{diag:3}) are morphisms in the fiber, gives a commutative diagram 
in $\Tho_A$: 
$$\xymatrix@R-1pc{([s_1],\dots,[s_{p}],[s-\sum_is_i]]) \ar[dd]_{H(\eta_X)}\ar[rr]^{H(\al)} & & ([s'_1],\dots,[s'_{q}],[s'-\sum_js'_j])\ar[dd]^{H(\eta_Y)} \\
& ([r_1],\dots,[r_k]) \ar[ul]\ar[ur]\ar[dl]\ar[dr] & \\
([s_1],\dots,[s_{p}],[s-\sum_is_i]]) \ar[rr]^{H(\Lambda\Pi\al)} & & ([s'_1],\dots,[s'_{q}],[s'-\sum_js'_j]) }$$
Now the maps in this diagram assemble to give a commutative diagram of bijections on the set $[s']$ where $s'=|E|=|\bar E|=\sum_ir_in^{\mu'(i)}$, where the left vertical map is induced by $\underline f$, the right vertical  map is induced by $\underline f'$, the top one by $\la$ and the bottom one by $\la_E$. Because this is now a commutative diagram of invertible maps, we can conclude that the outer square commutes, which precisely  gives the equality $\underline f'\circ \la=\la_E\circ \underline f$. 
\end{proof}

\subsection{From the Thomason construction to Cantor algebras}\label{sec:square}

There is a canonical functor~$\Tho_A$ to~$\bCan_A^\x$ coming from the universal property of the homotopy colimit~$\Tho_A$ and the defining structure of Cantor algebras. We give it here explicitly and check that it is symmetric monoidal using the universal property. Finally, we prove that it is compatible with the zig-zag of functors we have so far defined between the two categories. 

Let
\[
F\colon\Tho_A\ \rar\ \bCan_A^\x
\]
be the functor defined on objects by 
$$F([r_1],\dots,[r_m])=\rmC_A[r_1+\dots+r_m]$$
and on morphisms by taking 
$$(\psi,\mu,{\bf f})\colon([r_1],\dots,[r_m])\rar ([s_1],\dots,[s_n]) \ \ \textrm{with}\ \ f_j\colon\bigsqcup_{i\in\psi^{-1}(j)}[r_in^{\mu(i)}]\stackrel{\cong}\rar [s_j]$$
to the isomorphism of Cantor algebras represented by 
$$(E,[s_1+\dots+s_n],\underline f)\colon\rmC_A[r_1+\dots+r_m]\ \rar \ \rmC_A[s_1+\dots+s_n]$$
for~$E=\bigsqcup_i r_i\x A^{\mu(i)}$ considered as an expansion of~$[r_1]\oplus\dots\oplus[r_k]=[r_1+\dots+r_m]$ and with 
$$\underline f\colon E=\bigsqcup_{1\le i\le m}r_i\x A^{\mu(i)} \ \stackrel{\cong}\rar\ \bigsqcup_{1\le j\le n}\ \bigsqcup_{i\in \psi^{-1}(j)} r_i\x A^{\mu(i)}\ \xrightarrow{\sqcup_j f_j}\ \bigsqcup_{1\le j\le n} [s_j]\xrightarrow{\oplus} [s_1+\dots+s_n]$$
where the first map reorders the factors. 

\begin{prop}
The functor $F\colon\Tho_A\ \rar\ \bCan_A^\x$ is symmetric monoidal. 
\end{prop}

\begin{proof}
The functor $F$ is the functor induced, via the universal property of the homotopy colimit $\Tho_A$ \cite[Prop 3.21]{Thomason}, from the free Cantor algebra functor 
$\rmC_A\colon \bSet^\x\to \bCan_A^\x$ and the natural transformation $\eta\colon \rmC_A\to \rmC_A\circ \Si_A$ defined on objects as the morphism $\eta_r\colon \rmC_A[r]\to \rmC_A[rn]$ represented by the triple 
$(E=[r]\x A,[rn],\la_E)$. Indeed, for the naturality, we need to check that for any bijection $\la\colon [r]\to [r]$,  the diagram 
$$\xymatrix{\rmC_A[r] \ar[rr]^-{\rmC_A(\la)} \ar[d]_{\eta_r} && \rmC_A[r] \ar[d]^{\eta_r}\\
\rmC_A[rn] \ar[rr]^-{\rmC_A\Si_A(\la)} && \rmC_A[rn] }$$
commutes, which precisely follows from the commutativity of diagram (\ref{diag:SiA}) defining $\Si_A(\la)$. 
(One can also check directly that $F$ is a symmetric monoidal functor. In fact it is strict symmetric monoidal.) 
%For the monoidality, recall that 
%\[
%([r_1],\dots,[r_m])\oplus  ([s_1],\dots,[s_n])= ([r_1],\dots,[r_m],[s_1],\dots,[s_n])
%\]
%while 
%\[
%\rmC_A[r_1+\dots+r_m]\oplus \rmC_A[s_1+\dots+s_n]= \rmC_A[r_1+\dots+r_m+s_1+\dots+s_n].
%\]
%This shows that we have~$F(X\oplus Y)=F(X)\oplus F(Y)$ for objects~$X=([r_1],\dots,[r_m])$ and~$Y=([s_1],\dots,[s_n])$ in the category~$\Tho_A$. And if~$(\psi,\mu,f)$ and~$(\psi',\mu',f')$ are two morphisms of~$\Tho_A$, one sees likewise that~\hbox{$F((\psi,\mu,f)\oplus (\psi',\mu',f'))=F(\psi,\mu,f)\oplus F(\psi',\mu',f')$}, as in both cases, the maps and expansions are just sat ``side-by-side'' in the sum. 
%The symmetry~$\s_{X,Y}\colon X\oplus Y\to Y\oplus X$ in~$\Tho_A$ is the map
%$$(\s_{m,n},{\bf 0},\id)\colon([r_1],\dots,[r_m],[s_1],\dots,[s_n])  \ \rar \ ([s_1],\dots,[s_n],[r_1],\dots,[r_m])$$
%induced by the block permutation~$\s_{m,n}\colon[m+n]\to [n+m]$. Now the map  
%$$F(\s_{X,Y})\colon\rmC_A[r_1+\dots+r_m+s_1+\dots+s_n]\ \rar\ \rmC_A[s_1+\dots+s_n+r_1+\dots+r_m]$$
%is then precisely the map induced by the (larger) block permutation
%\[
%\s_{\sum r_i,\sum s_i}\colon[r_1+\dots+r_m+s_1+\dots+s_n]\to [s_1+\dots+s_n+r_1+\dots+r_m]
%\]
%and hence coincides with the symmetry of~$\bCan_A^\x$.  
\end{proof} 

Note that~$F$ has image in the subcategory~$I\colon\Exp_A\inc \bCan_A^\x$. We will denote by~$F_0\colon\Tho_A\to \Exp_A$ the functor considered with~$\Exp_A$ as target category. 

We have constructed a diagram of functors 
\[\xymatrix{\Tho_A \ar[rr]^{F}\ar@{-->}[rrd]^{F_0} & &  \bCan_A^\times\\
\Lev_A \ar[u]^{H}_\sim\ar[rr]_J^\sim &&  \ar[u]_I^\sim \Exp_A }
\]
and we are left to check that the diagram commutes up to homotopy, which will follow once we have shown that the bottom triangle commutes up to homotopy. 

\begin{prop}
Consider the diagram 
$$\xymatrix{\Lev_A \ar[dr]_{J} \ar[r]^-H & \Tho_A \ar[d]^{F_0}\\
& \Exp_A}$$
of categories. There is a natural transformation~$\eta\colon J \to F_0\circ H$. In particular, the corresponding diagram of classifying spaces commutes up to homotopy. 
\end{prop}

\begin{proof}
Given an object~$([r],(P_1,\dots,P_p))$ of~$\Lev_A$, we have that
$$J([r],(P_1,\dots,P_p))=\rmC_A[r]=F_0\circ H([r],(P_1,\dots,P_p)).$$ Indeed,~$J$ forgets the subsets~$(P_1,\dots,P_p)$, while~$H$ takes 
$([r],(P_1,\dots,P_p))$ to~$([r_1],\dots,[r_p],[r-\sum_i r_i])$ in~$\Tho_A$ for~$r_i=|P_i|$ (dropping the last component if it is zero), which is then taken by~$F_0$ to the Cantor algebra~\hbox{$\rmC_A[\sum_i r_i+r-\sum_i r_i]=\rmC_A[r]$}. However we see that the subsets~$P_i$ have been ``permuted'' in the composition~$F_0\circ H$, which will affect the value of the composed functor~$F_0\circ H$ on morphisms. The natural transformation~$\eta$ is given by that permutation: 
Let~$\eta([r],(P_1,\dots,P_p))\colon\rmC_A[r] \ \rar\ \rmC_A[r]$ be induced by the bijection  
%\begin{multline*}
$$[r]\rar P_1\sqcup\dots\sqcup P_p\sqcup\left([r]\minus\sqcup_iP_i\right) \stackrel{(\textrm{lexi})}\rar [r_1]\sqcup\dots\sqcup [r_p]\sqcup [r-\sum_i r_i]\xrightarrow{\oplus} [\sum r_i + (r-\sum r_i)]=[r].$$
%\end{multline*}
% where the first map takes each subset~$P_i$ of~$[r]$, as well as the complement of~$\sqcup_i P_i$, to itself, and the second map is the canonical order preserving map. 
We need to check that this is natural with respect to the morphisms of~$\Lev_A$. It is enough to check naturality for each of the generating morphisms~$d_0,d_i,d_p$ and~$\la$.
Recall that~$d_0$ expands along~$P_0$, that~$d_i$ takes the union of~$P_i$ and~$P_{i+1}$, and that~$d_P$ forgets~$P_p$. We see that~$\eta$ is precisely the map that ensures that we apply these operations to the ``same'' subsets of~$[r]$: indeed,~$J$ considers each~$P_i$ as the original~$P_i\subset [r]$, while~$F_0\circ H$ considers~$P_i$ as the subset of~\hbox{$P_1\sqcup\dots\sqcup P_p\sqcup ([r]\minus \sqcup_iP_i)\cong [r]$}, and~$\eta$ precisely takes the one version of~$P_i$ to the other. 
Likewise, applying a bijection~$\la$ will be compatible by the definition of the functors involved: 
Given a bijection~$\la\colon[r]\to [r]$ inducing a morphism~\hbox{$\la\colon([r],(P_1,\dots,P_p))\to ([r],(\la P_1,\dots,\la P_p))$} in~$\Lev_A$, we get a commutative diagram 
$$\xymatrix{\rmC[r] \ar[r]^{J(\la)=\la} \ar[d]_\eta^\cong& \rmC[r] \ar[d]^\eta_\cong\\
 \big(\rmC_A[r_1]\oplus\dots\oplus \rmC_A[r_p]\big)\oplus \rmC_A[r-\sum_i r_i] \ar@{=}[d] \ar[r]^{\oplus_i \la_i} &\big(\rmC_A[r_1]\oplus\dots\oplus \rmC_A[r_p]\big)\oplus \rmC_A[r-\sum_i r_i] \ar@{=}[d]\\
\rmC_A[r] \ar[r]^{F_0\circ H(\la)} &  \rmC_A[r]
}$$
where~$\la_i$ is the map induced by~$\la$ on each~$P_i$ (resp. on~$[r]\minus\sqcup_iP_i$), through their canonical identification with~$[r_i]$ via the lexicographic ordering map~$\la_{P_i}$. 
\end{proof}

% !TEX root = allinone_main_v2.tex

\section{Homotopy theory and homology}\label{sec:main}

In this section, we identify the spectrum~$\bbK(\Tho_A)$ from the preceding section with the Moore spectrum for~$\bbZ/(n-1)$, where~$n$ is the cardinality of the set~$A$ as above.  
We will use this and the relationship between the category~$\Tho_A$
and the Higman--Thompson groups to give, in
Section~\ref{sec:explicit}, concrete computations of the homology of
the Higman--Thompson groups.

%\Mnote{Yes.
%\begin{remark}
%A Cantor algebra of type~$1$ is a set~$X$ together with a bijection~$f\colon X\to X$. A model for the free Cantor algebra of type~$1$ on a set~$S$ is~$X=\bbZ\times S$ with~$f(n,x)=(n+1,x)$. If~$S=[r]$, then its automorphism group is the wreath product~$\bbZ\wr\Sigma_r$. For these, homological stability is known, and also the stable homology, for instance from Segal's generalization of the BPQ theorem: the algebraic K-theory spectrum is~$\bbS\vee\Sigma\bbS$ with
%\[
%\Omega^\infty_0(\bbS\vee\Sigma\bbS)=Q_0(\rmS^0)\times Q(\rmS^1).
%\]
%This agrees with our computation, because the cofiber of~$1-1=0$ on~$\bbS$ is also~$\bbS\vee\Sigma\bbS$.
%\end{remark}
%}

\subsection{Entry of the Moore spectra}\label{sec:Moore}

For any integer~$n\geqslant 1$ let~$\bbM_n$ denote the homotopy cofiber of the multiplication by~$n$ map on the sphere spectrum~$\bbS$, so that there is a homotopy cofiber sequence as follows.
\begin{equation}\label{eq:cofib}
\bbS\overset{n}{\longrightarrow}\bbS\longrightarrow\bbM_n
\end{equation}
The spectrum~$\bbM_n$ is the Moore spectrum for~$\bbZ/n$, also known as the {\em mod~$n$ Moore spectrum}.

\begin{theorem}\label{thm:hocolim_is_Moore}
Let~$A$ be a finite set of cardinality~$n\geqslant 2$. There is an equivalence of spectra
\[
\bbK(\Tho_A)\simeq\bbM_{n-1}.
\]
\end{theorem}

\begin{proof} Recall from Section~\ref{def:Sigma_A} that the category~$\Tho_A$ is defined from the diagram of categories on~$\bbN$ which takes the unique object to the category~$\bSet^\x$ and the arrow~$k$ to the functor~$\Sigma_A^k$ that corresponds to taking the product with~$A^k$. 
By the Barratt--Priddy--Quillen theorem \cite{BP72}, we have that~$\bbK(\bSet^\x):=\bbK(\Set^\x)\simeq\bbS$ and  the functor~$\Sigma_A$  induces multiplication with the cardinality~$n$ of~$A$ on~$\bbS$. 
We apply Thomason's formula~\eqref{eq:thmTho} to our definition~\eqref{eq:defTho} and get
\[
\bbK(\Tho_A)
=
\bbK(\hocolim_{\bbN}(\Set^\times,\Sigma_A))
\simeq
\hocolim_{\bbN}(\bbK(\Set^\times),\bbK(\Sigma_A))
\simeq
\hocolim_{\bbN}(\bbS,n).
\]

There is a spectral sequence
\[
\rmE^2_{p,q}=\rmH_p(\bbN;(\pi_q\bbS,n))\ \Longrightarrow\ \pi_*\hocolim_{\bbN}(\bbS,n)
\]
that computes the homotopy groups of the homotopy colimit, see~\cite[Sec.~3]{Thomason}. Here~$(\pi_q\bbS,n)$ is the diagram of abelian groups on~$\bbN$ that takes the object to~$\pi_q\bbS$ and the morphism~$1$ to multiplication by~$n$. 
The monoid ring~$\bbZ[\bbN]\cong\bbZ[\rmT]$ is polynomial on one generator~$\rmT$, so that we can use the standard Koszul resolution 
\[
\bbZ[\rmT]\stackrel{\rmT-1}\rar \bbZ[\rmT]\rar\bbZ
\] 
to compute the~$\rmE^2$ page. As the resolution has length one, the spectral sequence degenerates at~$\rmE^2$, and yields a long exact sequence 
\begin{equation}\label{eq:les}
\cdots
\xrightarrow{\hspace*{2em}}\pi_*\bbS
\overset{n-1}{\xrightarrow{\hspace*{2em}}}\pi_*\bbS
\xrightarrow{\hspace*{2em}}\pi_*\hocolim_{\bbN}(\bbS,n)
\xrightarrow{\hspace*{2em}}\cdots.
\end{equation}
Let~$\epsilon\in \pi_0\hocolim_{\bbN}(\bbS,n)=[\,\bbS,\hocolim_{\bbN}(\bbS,n)\,]$ be the image of~$1\in \pi_0\bbS=[\,\bbS,\bbS\,]$ in the long exact sequence. 
From the long exact sequence we get that~$(n-1)\epsilon=0$. Therefore, this map factors through the Moore spectrum to give a map
\[
\overline\epsilon\colon\bbM_{n-1}\longrightarrow\hocolim_{\bbN}(\bbS,n).
\]
Comparison of the long exact sequence obtained from the homotopy cofibration sequence~\eqref{eq:cofib} with the long exact sequence~\eqref{eq:les} shows that~$\overline\epsilon$ induces an isomorphism on homotopy groups, and thus that it is an equivalence. 
\end{proof}

Together with Corollary~\ref{cor:KCan=KTho} the proposition gives the following result.

\begin{corollary}\label{cor:KCan=M}
For any finite set~$A$ of cardinality~$n\geqslant 2$ there is an equivalence of spectra
\[
\bbK(\Can_A^\times)\simeq\bbM_{n-1}.
\]
\end{corollary}

\begin{remark}\label{rem:components}
It is instructive to work out the implications of Corollary~\ref{cor:KCan=M} on the level of components. The cofiber sequence \eqref{eq:cofib} shows that the group~$\pi_0\bbM_{n-1}$ is the cokernel of the multiplication by~$n-1$ map on the group~$\pi_0\bbS=\bbZ$, so that it is cyclic of order~\hbox{$n-1$}. On the other hand the abelian group~\hbox{$\pi_0\Omega^\infty\bbK(\Can_A^\times)=\pi_0\rmK(\Can_A^\times)$} is the group completion of the abelian monoid~$\pi_0|\Can_A^\times|$. The latter can be identified with~$\{0,1,\dots,n-1\}$ as a set, where an integer~$r$ corresponds to the free Cantor algebra~$\rmC_{A}[r]$ of type~$A$ on~$r$ generators. The monoid structure is dictated by~$\rmC_{A}[r]\oplus\rmC_{A}[s]=\rmC_A[r+s]$ and~$\rmC_A[r+(n-1)]\cong\rmC_A[r]$ if~$r\geqslant1$. Note that the neutral element~$0$ is the only invertible element in this monoid. The group completion of this monoid is~$\bbZ/(n-1)$ by the theorem, but this can of course also be worked out by hand: Once~$1$ is inverted, the element~$n-1$ is identified with~$0$ and we immediately obtain the group~$\bbZ/(n-1)$. 
\end{remark}

\subsection{The (stable) homology of the Higman--Thompson groups}

If~$\bbX$ is a spectrum, its associated infinite loop space~$\Omega^\infty\bbX$ is a group-like~$\rmE_\infty$--space: the monoid of components is a group. 
It follows that all its components are homotopy equivalent, and in the following, we denote by~$\Omega_0^\infty\bbX$ the component corresponding to the zero element~$0\in\pi_0\bbX$. 
We will now relate the homology of the zeroth component~$\Omega_0^\infty\bbK(\Can_A)$ to that of the Higman--Thompson groups. 
Together with  Corollary~\ref{cor:KCan=M} this will yield the main identification we are after, namely the isomorphism of the homology of the Higman--Thompson groups with that of~$\Omega_0^\infty\bbM_{n-1}$.
%We will now pass from equivalences of spaces and spectra to maps that are merely homology isomorphisms. The slogan is: K-theory computes stable homology.

Recall from Section~\ref{sec:stabsubsec} the stabilization homomorphism~$s_r\colon\rmV_{n,r}\to \rmV_{n,r+1}$ and let
\[
\rmV_{n,\infty}=\colim(\rmV_{n,1}\longrightarrow \rmV_{n,2}\longrightarrow \cdots)
\] 
denote the associated {\em stable group}. The homology of~$\rmV_{n,\infty}$ is usually called the {\em stable homology} of the groups~$\rmV_{n,r}$ for~$r\geqslant 1$, but 
a direct consequence of our stability theorem, Theorem~\ref{thm:HS},  is that
\[
\rmH_*(\rmV_{n,r})\cong \rmH_*(\rmV_{n,\infty}).
\]
So in the case at hand, {\em all} the homology is stable and hence it is enough to identify the homology of~$\rmV_{n,\infty}$. 

Given a monoid~$M$, one can form its bar construction~$\rmB M$ and the loop space~$\Omega\rmB M$ is a group-like space, known as the {\em group completion of~$M$}. The group completion theorem of McDuff and Segal~\cite{McDuff-Segal}  identifies in good cases the  homology of~$\Omega\rmB M$ with that of its ``stable part.'' We will use this theorem to compute the homology of~$\Omega^\infty\bbK(\Can_A^\times)\simeq \Omega\rmB|\Can_A^\x|$, the group completion of the~$\rmE_\infty$ monoid~$|\Can_A^\x|$. 

\begin{theorem}\label{thm:gpcomp}
Let~$A=\{a_1,\dots,a_n\}$ with~$n\geqslant2$. There is a map 
\[
\rmB\rmV_{n,\infty}\longrightarrow\Omega_0^\infty\bbK(\Can_A^\times)
\]
which induces an isomorphism in homology 
%\[
%\rmH_*(\rmV_{n,\infty})
%\cong
%\rmH_*(\Omega_0^\infty\bbK(\Can_A^\times))
%\]
with all systems of local coefficients on~$\Omega_0^\infty\bbK(\Can_A^\times)$. 
\end{theorem}

%\begin{theorem}\label{gpcomp}
%$H_*(\Omega B |Cantor_n|)\cong  H_*(\bbZ/_{n-1}\x BV_{n,\infty})$ 
%\end{theorem}

\begin{proof}
We apply the group completion theorem to the monoid~$M=|\Can_A^\x|$. More precisely, we will 
use Theorem 1.1 in \cite{Randal-Williams}, which makes explicit the relevant result in \cite{McDuff-Segal}. 
The monoid~$M$ is homotopy commutative, as it is the classifying space of a permutative category. 
It has components indexed by~$0,\dots,n-1$, forming a monoid in the way describe in Remark~\ref{rem:components}. To apply Theorem 1.1 of  \cite{Randal-Williams}, we use the constant sequence of elements of~$M$ given by~$\rmC_A[1],\rmC_A[1],\dots$. 
We need to check that for every~$m\in M$, the component of~$m$ in~$M$ is a right factor of the component of some finite sum~$\rmC_A[1]\oplus \dots\oplus \rmC_A[1]$, which is obvious as every component is reached this way except for the zero component which is a right factor of any such sum. 

Form  the colimit 
\[
M_\infty=\colim\big(|\Can^\x_A|\stackrel{\oplus \rmC_A[1]}\rar |\Can^\x_A|\stackrel{\oplus \rmC_A[1]}\rar\cdots).
\]
Theorem 1.1 in \cite{Randal-Williams} says that there is a homology isomorphism
\[
M_\infty\longrightarrow\Omega\rmB|\Can^\x_A|\simeq \Omega^\infty\bbK(\Can_A)
\]
with respect to all local coefficient systems on the target.  
%We have already computed that~$\Omega B|\Can_A|$ has components~$\bbZ/_{n-1}$, so we are left to identify the homology of its components, which we do using~$M_\infty$ and the equivalence just computed. 
%We are left to identify the~$0$th component of~$M_\infty$. Under the equivalence of groupoids~$$\Can^\x_A\simeq \coprod{0\leqslant r\leqslant n-1}{V_{n,r}},$$  
Now the zeroth component of~$M_\infty$ can be identified with the colimit on classifying space of the maps 
\[
\rmV_{n,0}\longrightarrow\rmV_{n,1}\stackrel{s_1}\longrightarrow\cdots\longrightarrow \rmV_{n,n-1} \stackrel{s_{n-1}}\longrightarrow\rmV_{n,n} \stackrel{s_n}\longrightarrow\rmV_{n,n+1}\longrightarrow\cdots
\]
of groups , and this colimit of groups is the group~$\rmV_{n,\infty}$ in the stability statement. The result follows. 
\end{proof}

%\begin{proof}
%This follows from the agreement of various models for algebraic K-theory. Thomason~\cite[p.~1653]{Thomason} argues that there is an equivalence
%\[
%\Omega_0^\infty\bbK(\Mod_R^\times)\simeq\BGL_\infty(R)^+,
%\]
%where~$\Mod_R^\times$ is the symmetric monoidal category of finitely generated projective~$R$-modules and their isomorphisms over a given ring~$R$, and~$\BGL_\infty(R)^+$ is the Quillen plus construction on the classifying space of the infinite general linear group. The same arguments give an equivalence
%\[
%\Omega_0^\infty\bbK(\Can_a^\times)\simeq\rmB\rmV_{a,\infty}^+,
%\]
%where~$\rmV_{a,\infty}$ is infinite Higman--Thompson group. Since the Quillen plus construction~$\rmB\rmV_{a,\infty}\to\rmB\rmV_{a,\infty}^+$ is a homology equivalence, we are done.
%\end{proof}

%\begin{corollary}
%For all~$n\geqslant2$ there is a homology  isomorphism
%\[
%\rmH_*(\rmV_{n,\infty})
%\cong
%\rmH_*(\Omega_0^\infty\bbM_{n-1}).
%\]
%\end{corollary}

%%%

%\subsection{The unstable homology of the Higman--Thompson groups}

We are now ready to prove the main result of this text.

%\begin{theorem}\label{thm:main}
%For all~$n\geqslant2$ and all~$r\geqslant1$ there is a map  isomorphisms
%\[
%\rmH_*(\rmV_{n,r})
%\cong
%\rmH_*(\Omega_0^\infty\bbM_{n-1})
%\]
%in homology.
%\end{theorem}

%Note that the right hand side does not depend on~$r$. The proof reflects that.

\begin{proof}[Proof of Theorem~\ref{thm:main}]
%From Section~\ref{sec:Cantor} we have isomorphisms~$\rmV_{n,r}\cong\rmV_{n,r+(n-1)}$ of groups. Therefore, we only need to compute~$\rmH_*(\rmV_{n,r})$ for arbitrarily large~$r$. However, 
As a consequence of our stability result, Theorem~\ref{thm:HS}, for all~$r\geqslant 1$, the stabilization map~\hbox{$s_r\colon\rmV_{n,r}\to \rmV_{n,r+1}$} induces an isomorphism in homology with coefficients in any~$\rmH_1(\rmV_{n,\infty})$--module. Note that, in particular, we have an isomorphism~$\rmH_1(\rmV_{n,r})\cong\rmH_1(\rmV_{n,\infty})$, so that~$\rmH_1(\rmV_{n,\infty})$--modules are the same as~$\rmH_1(\rmV_{n,r})$--modules.  It follows that 
the map~$\rmB\rmV_{n,r}\to \rmB\rmV_{n,\infty}$ induces an isomorphism in homology with abelian coefficients for all~$r\geqslant 1$. 
%$\rmH_*(\rmV_{n,r})\cong\rmH_*(\rmV_{n,\infty})$ for all~$r\geqslant 1$. 
Theorem~\ref{thm:gpcomp} gives that there is a map~\hbox{$\rmB\rmV_{n,\infty}\to \Omega^\infty_0(\Can_A^\x)$} which induces an isomorphism in homology with all local coefficients on the target. (Note that~\hbox{$\pi_1 \Omega^\infty_0(\Can_A^\x)\cong\rmH_1( \Omega^\infty_0(\Can_A^\x))\cong\rmH_1(\rmV_{n,\infty})$}, so these are again~$\rmH_1(\rmV_{n,r})$--modules as above.)
%identified the homology of~$\rmV_{n,\infty}$ to that of~$\Omega^\infty_0(\Can_A^\x)$ which, 
By Corollary~\ref{cor:KCan=M}, we have a homotopy equivalence~$\Omega^\infty_0\bbK(\Can_A^\x)\simeq \Omega_0^\infty\bbM_{n-1}$, proving the result.
\end{proof}

%In the following Section~\ref{sec:explicit} we will give explicit consequences of Theorem~\ref{thm:main}.

%%%

\section{Computational consequences}\label{sec:explicit}

We will in this section explain how Theorem~\ref{thm:main} can be used to give explicit consequences for the homology of the Higman--Thompson groups. Concretely, we will compute the abelianizations and Schur multipliers of the Higman--Thompson groups directly from Theorem~\ref{thm:main}, and we will also completely decide which of the groups are integrally or rationally acyclic. We recover old results with new methods, as well as prove new results. In particular, we prove the acyclicity of the Thompson group~$\rmV$. 

The results in this section are based on computations of the homology groups of the infinite loop space of the Moore spectrum~$\bbM_{n}$ with classical methods from homotopy theory. Given a spectrum~$\bbX$, the stable homotopy groups~$\pi_*\bbX$ agree with the (unstable) homotopy groups~$\pi_*\Omega^\infty\bbX$ of the underlying infinite loop space~$\Omega^\infty\bbX$. The situation is different for homology, however. The homology of the Moore spectrum~$\bbM_n$ is, up to a shift, the homology of the mod~$n$ Moore space:~$\rmH_0\bbM_n\cong\bbZ/n$ and~$\rmH_d\bbM_n=0$ for~$d\not=0$. In contrast, the homology of the underlying infinite loop space~$\Omega^\infty\bbM_n$ is more difficult to compute. We will here give some partial computations of these homology groups. Further computations can be obtained by working harder. 

%%%

\subsection{Abelianizations and Schur multipliers}

%In general, it can be a very laborious endeavor to compute the homology of an infinite loop space, even (or perhaps especially) if the homology of the spectrum is very simple, as for the Moore spectra. We will here content ourselves with the computation of~$\rmH_1$ and~$\rmH_2$ in order to confirm and extend the known result from the literature.

In this section, we compute~$\rmH_1$ and~$\rmH_2$ as well as the first non-trivial homology group of~$\rmV_{n,r}$  by computing these groups for~$\Omega_0^\infty\bbM_{n-1}$. We confirm and extend the known results from the literature.

\begin{proposition}\label{prop:H1}
For all~$n\geqslant 2$ and~$r\geqslant1$ there are isomorphisms
\[
\rmH_1(\rmV_{n,r})\cong
\begin{cases}
0 		& n\text{ even}\\
\bbZ/2 	& n\text{ odd.}\\
\end{cases}
\]
\end{proposition}

\begin{proposition}\label{prop:non-trivial}
For all~$n\geqslant 3$, we have that
\[
\rmH_d(\rmV_{n,r})\cong
\begin{cases}
0 & 0<d<2p-3\\
\bbZ/p & d=2p-3
\end{cases}
\]
for~$p$ the smallest prime dividing~$n-1$, and 
\[
\rmH_{2q-3}(\rmV_{n,r})\neq 0
\]
for~$q$ any prime dividing~$n-1$.
%\begin{center}
%\begin{tabular}{ll}
%$\rmH_{0<*<2p-3}(\rmV_{n,r})=0$ & for~$p$ the smallest prime dividing~$n-1$\\
%$\rmH_{2p-3}(\rmV_{n,r})\cong \bbZ/p$ & for~$p$ the smallest prime dividing~$n-1$\\
%$\rmH_{2q-3}(\rmV_{n,r})\neq 0$ & for~$q$ any prime dividing~$n-1$. 
%\end{tabular}
%\end{center}
\end{proposition}

Proposition~\ref{prop:H1} is essentially the case~$p=2$ of Proposition~\ref{prop:non-trivial}. We prove both propositions together. 

\begin{proof}[Proof of Propositions~\ref{prop:H1} and~\ref{prop:non-trivial}]
%For~$n\geqslant 3$, the connected cover of~$\bbM_{n-1}$ is non-contractible. If~$d>0$ is the first dimension with a non-trivial homotopy group, then~$\rmH_d\rmV_{n,r}\not=0$ by Hurewicz' theorem.
By Theorem~\ref{thm:main}, it is equivalent to compute these homology groups for~$\Omega^\infty_0\bbM_{n-1}$. 
The space~$\Omega^\infty_0\bbM_{n-1}$ is a connected infinite loop space, and so
\[
\rmH_1\Omega^\infty_0\bbM_{n-1}
\cong\pi_1\Omega^\infty_0\bbM_{n-1}
\cong\pi_1\bbM_{n-1}.
\]
As~$\pi_1\bbS\cong \bbZ/2$, the latter can easily be computed from the cofiber sequence~\eqref{eq:cofib} of spectra to be the cokernel of multiplication by~$n-1$ on~$\pi_1\bbS=\bbZ/2$, which proves the first proposition.

For the second proposition, we assume that~$n\geqslant 3$ so that~$n-1\geqslant 2$. Let~$p$ be a prime. The~$p$-parts of the homotopy groups of the sphere spectrum are zero between dimension~$0$ and~\hbox{$2p-3$} and then~\hbox{$\pi_{2p-3}(\bbS)\otimes \bbZ_{(p)}\cong \bbZ/p$}. 
Now multiplication by~$n-1$ is an isomorphism on~$\bbZ/p$ if and only if~$p$ does not divide~$n-1$, and if it does, the map is zero. The result follows using the same long exact sequence now for homotopy groups with coefficients and Hurewicz's theorem.   
\end{proof}

%\begin{remark}\label{rem:non-trivial} The proof above extends to compute the first non-trivial homology of~$\rmV_{n,r}$ at any prime~$p$, and in particular its first non-trivial homology group. Indeed, the~$p$-parts of the homotopy groups of the sphere spectrum vanish below dimension~\hbox{$2p-3$}, where they are~$\bbZ/p$. Consequently, the same vanishing holds for the~$p$-parts of the homotopy groups of the Moore spectrum, and in dimension~\hbox{$2p-3$} we see the cokernel of multiplication by~$n-1$ on~$\bbZ/p$. It follows that, when~$n\geqslant3$, 
%\begin{equation}\label{eq:not_acylic}
%\rmH_{2p-3}(\rmV_{n,r})\not=0
%\end{equation}
%for any prime~$p$ that divides~$n-1$. And if~$p$ is the smallest prime dividing~$n-1$, we have that~$\rmH_{2p-3}(\rmV_{n,r})=\bbZ/p$ is the first non-trivial homology group of~$\rmV_{n,r}$.  
%The above proposition gives the case~$p=2$. \nnote{I would prefer to have this remark as the result, with the proposition as a particular case, specially as we use it later.}
%\end{remark}

The following result recovers and extends the computation of Kapoudjian~\cite{Kapoudjian}, who worked out the case~$r=1$ by entirely different methods.

\begin{proposition}\label{prop:H2}
For all~$n\geqslant 2$ and~$r\geqslant1$ there are group isomorphisms
%\[
%\rmH_1\rmV_{n,r}\cong
%\begin{cases}
%0 		& n\text{ even}\\
%\bbZ/2 	& n\text{ odd}\\
%\end{cases}
%\]
%and
\[
\rmH_2(\rmV_{n,r})\cong
\begin{cases}
0						& n\text{ even}\\
\bbZ/4 					& n\equiv 3\text{ mod }4\\
\bbZ/2\oplus\bbZ/2 		& n\equiv 1\text{ mod }4.
\end{cases}
\]
\end{proposition}

\begin{proof}
Again, by Theorem~\ref{thm:main}, it is equivalent to compute these homology groups for~$\Omega^\infty_0\bbM_{n-1}$. 
Let us first tick off the case when~$n$ is even, so that~$n-1$ is odd. 
%Then the cofibration sequence~\eqref{eq:cofib} shows that the space~$\Omega^\infty_0\bbM_{n-1}$ is~$2$-connected, so that~$\rmH_1$ and~$\rmH_2$ vanish.
By Proposition~\ref{prop:non-trivial}, when~$n\geqslant 3$ with~$n-1$ odd, the first possible non-trivial homology group of~$\rmV_{n,r}$ is in degree~$2\cdot 3-3=3$ (if~$3$ divides~$n-1$). In particular~$\rmH_2$ vanishes. 
If~$n=2$, the group vanishes because multiplication by~$n-1$ is the identity, making~$\bbM_1$ the trivial spectrum. (See also the more general Theorem~\ref{thm:Vacyclic}.) 

Let us now assume that~$n$ is odd, and write~$X=\Omega^\infty_0\bbM_{n-1}$. We have that~$\pi_1X\cong \bbZ/2$. Consider the Postnikov truncation~$X_2$, with the same first and second homotopy groups as~$X$. We have~\hbox{$\rmH_2X\cong\rmH_2X_2$}. As~$X$ is an infinite loop space, its first~$k$-invariant vanishes (see~\cite{Arlettaz}), so that 
\[
X_2\simeq \rmK(\pi_1X,1)\times \rmK(\pi_2X,2).
\] 
%where~$X{[1]}\cong \rmK(\pi_1X,1)\cong \rmK(\bbZ/2,1)$ and~$X{[2]}\cong \rmK(\pi_2X,2)$. 
Now the K\"unneth theorem gives
\[
\rmH_2X_2 \cong  \rmH_2\rmK(\pi_1X,1) \oplus \rmH_2\rmK(\pi_2X,2) \cong  \rmH_2\rmK(\pi_2X,2),
\]
given that~$ \pi_1X\cong \bbZ/2$ has trivial~$\rmH_2$. As~$\rmH_2\rmK(\pi_2X,2)\cong \pi_2X$, we obtain that 
\[
\rmH_2\Omega^\infty_0\bbM_{n-1}\cong \pi_2\Omega^\infty_0\bbM_{n-1}\cong \pi_2\bbM_{n-1},
\]
%and let~$\widetilde\Omega^\infty_0\bbM_{n-1}$ denote the universal (i.e.~$2$-fold) covering. 
%The spectral sequence of the covering shows that the only contributions to~$\rmH_2\Omega^\infty_0\bbM_{n-1}$ are the following three groups. Firstly, there is the group
%\[
%\rmH_2(\bbZ/2;\rmH_0\widetilde\Omega^\infty_0\bbM_{n-1})
%\cong\rmH_2(\bbZ/2;\bbZ)=0,
%\]
%since the action on~$\rmH_0$ is trivial. Secondly, there is the group~$\rmH_1(\bbZ/2;\rmH_1\widetilde\Omega^\infty_0\bbM_{n-1})$, and this vanishes since~$\widetilde\Omega^\infty_0\bbM_{n-1}$ is simply-connected (by construction). Thirdly and finally, we need to consider the cokernel of the~$\mathrm{d}_3$ differential that hits the group~$\rmH_0(\bbZ/2;\rmH_2\widetilde\Omega^\infty_0\bbM_{n-1})$. This differential is given by the first~$k$-invariant, and that is always trivial for infinite loop spaces, \nnote{Markus promised some reference or extra explanation here :)}  so that this cokernel is
%\[
%\rmH_0(\bbZ/2;\rmH_2\widetilde\Omega^\infty_0\bbM_{n-1})
%=\rmH_2\widetilde\Omega^\infty_0\bbM_{n-1},
%\]
%the latter because the actions of the fundamental groups are trivial for h-spaces, so in particular for infinite loop spaces. In summary, we have seen that
%\[
%\rmH_2\Omega^\infty_0\bbM_{n-1}
%\cong
%\rmH_2\widetilde\Omega^\infty_0\bbM_{n-1},
%\]
%and it remains to compute the latter. This can be done as follows. The computation
%\[
%\rmH_2\widetilde\Omega^\infty_0\bbM_{n-1}
%\cong\pi_2\widetilde\Omega^\infty_0\bbM_{n-1}
%\cong\pi_2\Omega^\infty_0\bbM_{n-1}
%\cong\pi_2\bbM_{n-1}
%\]
which reduces the question to stable homotopy theory as above. 
When~$n$ is odd, the homotopy cofibre sequence~\eqref{eq:cofib} only shows that this group is of order~$4$. We need a further computation to identify the group. 

Let us assume that~$n\equiv 3$ mod~$4$. Then~$n-1$ is even but not divisible by~$4$, and we can use the cofibre sequence 
\[
\bbM_{(n-1)/2}\longrightarrow
\bbM_{n-1}\longrightarrow
\bbM_2
\]
from the octahedral axiom to see that~$\pi_2\bbM_{n-1}\cong\pi_2\bbM_2$, and the latter group is known to be cyclic of order~$4$, see~\cite[Thm.~3.2]{Mukai}.

%One way to see this is to note that the unit~$\bbS\to\KO$ of the real topological K-theory spectrum~$\KO$ induces an isomorphism
%\[
%\pi_2\bbM_n\cong\KO_2\bbM_n\cong\KO^0\Sigma\bbM_n
%\]
%for any~$n$ because~$\pi_i\bbS\to \KO_i\bbS$ is an isomorphism for~$i=1,2$. Now~$\KO^0\Sigma\bbM_2$ is the~$0$-th reduced real~K-theory of the real projective plane, which is cyclic of order~$4$. 

Lastly, if~$n\equiv 1$ mod~$4$, then~$n-1$ is divisible by~$4$. The factorization~$4k=2k\cdot 2$ gives a map~$j\colon\bbM_2\to\bbM_{4k}$ that has even degree on the bottom cell and odd degree on the top cell. Analyzing the diagram
\[
\xymatrix{
\cdots&\KO^0(\Sigma\bbS)\cong \bbZ/2\ar[l]&\KO^0(\Sigma\bbM_2)\cong \bbZ/4\ar[l]&\KO^0(\Sigma^2\bbS)\cong \bbZ/2\ar[l]&\cdots\ar[l]\\
\cdots&\KO^0(\Sigma\bbS) \cong \bbZ/2\ar[l]\ar[u]_0&\KO^0(\Sigma\bbM_{4k})\cong \ ?\ \ar[l]\ar[u]_{j^*}&\KO^0(\Sigma^2\bbS) \cong \bbZ/2\ar[l]\ar[u]_\cong&\cdots\ar[l]\\
}
\]
shows that~$j^*$ cannot be zero or epi, so that~$\KO^0\Sigma\bbM_{4k}$ must split into~$\bbZ/2$ summands. 
\end{proof}

%%%

\subsection{Acyclicity results}

We now deduce global results about the homology of the groups~$\rmV_{n,r}$. 

The following result has been suggested by Brown~\cite[Sec.~6]{Brown:Suggestion}.

\begin{theorem}\label{thm:Vacyclic}
For all~$r\geqslant 1$, the Thompson group~$\rmV\cong\rmV_{2,r}$ is integrally acyclic: 
\[
\rmH_d(\rmV)=\rmH_d(\rmV_{2,r})=0
\]
for all~$d\not=0$.
\end{theorem}

\begin{proof}
For~$n=2$ multiplication by~$n-1=1$ is homotopic to the identity, so that it is a self-equivalence of the sphere spectrum. Then the homotopy cofiber, the Moore spectrum~$\bbM_1$, is contractible, and the homology of the infinite loop space vanishes.
\end{proof}

\vbox{\begin{theorem}\label{thm:Qacyclic}
For all~$n\geqslant 3$ and~$r\geqslant 1$, the group~$\rmV_{n,r}$ is rationally but not integrally acyclic: 
\[
\rmH_d(\rmV_{n,r})\otimes\bbQ=0
\]
for all~$d\not=0$, but 
\[
\rmH_{2p-3}(\rmV_{n,r})\neq 0
\]
for any prime~$p$ such that~$p$ divides~$n-1$.
\end{theorem}}

\begin{proof}
For~$n\geqslant2$ multiplication by~$n-1$ is a rational equivalence, so that the Moore spectrum~$\bbM_{n-1}$ is rationally contractible,  and the rational homology groups vanish.
This proves the first part of the statement. The second part of the statement is given by Proposition~\ref{prop:non-trivial}.
\end{proof}

\begin{remark}
In the case~$n=2$, rationally acyclicity of the Thompson group has earlier been shown by Brown~\cite[Thm.~4]{Brown:Suggestion}, where the author also indicated that his proof can be adapted to prove the case~$n\geqslant 3$. The case~$n=2$, still only rationally, was later reproved by Farley in~\cite{Farley}. 
%In view of the integral acyclicity of the groups in that case, see our Theorem~\ref{thm:Vacyclic}, that turns out to be arguably the least interesting case. 
\end{remark}

%%%

We end by mentioning a consequence of our work for the commutator subgroups. When~$n$ is odd, Proposition~\ref{prop:H1} implies that the commutator subgroup~$\rmV_{n,r}^+$ of~$\rmV_{n,r}$ is an index-two subgroup. Let~$\tilde\Omega^\infty_0\bbK(\Can_A^\x)$ denote the universal cover of the space~$\Omega^\infty_0\bbK(\Can_A^\x)$. Shapiro's Lemma, Theorem~\ref{thm:main} and Theorem~\ref{thm:HS} applied to the twisted coefficients~\hbox{$M=\bbZ\rmH_1(\rmV_{n,r})\cong \bbZ\rmH_1(\Omega^\infty_0\bbK(\Can_A^\x))$} give the following:  

\begin{corollary}\label{cor:commutator}
There are homology isomorphisms 
\[
\rmH_*(\rmV_{n,r}^+)\cong \rmH_*(\rmV_{n,\infty}^+)\cong \rmH_*(\tilde\Omega^\infty_0\bbK(\Can_A^\x)).
\]
In particular, the groups~$\rmV_{n,r}^+$ and~$\rmV_{n,\infty}^+$ are not acyclic when~$n$ is odd. 
\end{corollary}
 
\begin{proof}
See Sections 3.1 and 3.2 of \cite{Randal-Williams+Wahl}. 
\end{proof}

%If~$n$ is an odd integer, the isomorphism~$\rmH_1(\rmV_{n,\infty})\cong\bbZ/2$ from Proposition~\ref{prop:H1} defines a unique normal subgroup~$\rmV_{n,\infty}^+$ of~$\rmV_{n,\infty}$ of index~$2$. That notation seems standard and the~$+$ does {\em not} refer to the Quillen plus construction. \nnote{This is the commutator subgroup, right? stability also holds for these guys but the twisted stability. Why do you want to consider the stable guy here? Also, the group completion holds with abelian coefficients  and we do also get a computation as the homology of the universal cover of~$\Omega_0^\infty\bbM_{n-1}$. Should we write this up?}

%\begin{proposition}
%For all odd~$n$ the group~$\rmV_{n,\infty}^+$ is not acyclic.
%\end{proposition}

%\begin{proof}
%The spectral sequence
%\[
%\rmH_*(\bbZ/2;\rmH_*(\rmV_{n,\infty}^+;\bbZ[1/2]))\Longrightarrow\rmH_*(\rmV_{a,\infty};\bbZ[1/2]))
%\]
%degenerates to give an isomorphism
%\[
%(\rmH_*(\rmV_{n,\infty}^+;\bbZ[1/2]))_{\bbZ/2}\cong\rmH_*(\rmV_{n,\infty};\bbZ[1/2]).
%\]
%Since the right hand side is not trivial, the left hand side is not either.
%\end{proof}

This answers a question of Sergiescu.

%%%

%%%

\section*{Acknowledgment}

This research has been supported by the Danish National Research
Foundation through the Centre for Symmetry and Deformation~(DNRF92) in
Copenhagen. Parts of this paper were conceived while the authors were
visiting the Hausdorff Research Institute for Mathematics (HIM) in
Bonn, and the paper was revised during a visit at the Newton Institute in Cambridge. We thank both institutes for their support. The authors would like to thank Dustin Clausen and Oscar Randal-Williams for pointing out gaps in early versions of this paper, and the referee for a report that helped us improving the paper. The first author would also like to thank Ricardo Andrade, Ken Brown, Bj\o{}rn Dundas, Magdalena Musat, Martin Palmer, and Vlad Sergiescu for conversations related to the subject of this paper. 

%%%

% !TEX root = allinone_main_v2.tex

%%%

\vfill

Department of Mathematical Sciences, 
NTNU Norwegian University of Science and Technology,
7491 Trondheim,
NORWAY\\
\href{mailto:markus.szymik@ntnu.no}{markus.szymik@ntnu.no}

Department of Mathematical Sciences, 
University of Copenhagen, 
Universitets\-parken 5, 
2100 Copenhagen,
DENMARK\\
\href{mailto:wahl@math.ku.dk}{wahl@math.ku.dk}

%%%

\end{document}